\documentclass{memo-l}
\usepackage[active]{srcltx}

\usepackage{amssymb,amsmath,amsfonts,bbm,pifont,upgreek,bbold,accents,ulem,xcolor
} 
\usepackage[colorlinks=true]{hyperref}
\hypersetup{urlcolor=blue, citecolor=red}
\usepackage[T1]{fontenc}

\renewcommand{\em}{\sl} 

\font\teneufm=eufm10
\font\seveneufm=eufm7
\font\fiveeufm=eufm5
\newfam\eufmfam
\textfont\eufmfam=\teneufm
\scriptfont\eufmfam=\seveneufm
\scriptscriptfont\eufmfam=\fiveeufm
\def\eufm#1{{\fam\eufmfam\relax#1}}

\newcommand\beq[1]{ \begin{equation}\label{#1} }
\newcommand{\eeq}{ \end{equation} }
\newcommand\beqa[1]{ \begin{eqnarray}\label{#1} }
\newcommand{\eeqa}{ \end{eqnarray} }
\newcommand{\beqano}{ \begin{eqnarray*} }
\newcommand{\eeqano}{ \end{eqnarray*} }
\newcommand\arr[1]{\left\{\begin{array}{l}#1\end{array}\right.}

\newcommand\equ[1]{{\rm(\ref{#1})}}

\expandafter\chardef\csname pre amssym.def
at\endcsname=\the\catcode`\@
\catcode`\@=11
\def\undefine#1{\let#1\undefined}
\def\newsymbol#1#2#3#4#5{\let\next@\relax
 \ifnum#2=\@ne\let\next@\msafam@\else
 \ifnum#2=\tw@\let\next@\msbfam@\fi\fi
 \mathchardef#1="#3\next@#4#5}
\def\mathhexbox@#1#2#3{\relax
 \ifmmode\mathpalette{}{\m@th\mathchar"#1#2#3}%
 \else\leavevmode\hbox{$\m@th\mathchar"#1#2#3$}\fi}
\def\hexnumber@#1{\ifcase#1 0\or 1\or 2\or 3\or 4\or 5\or 6\or 7\or
8\or
 9\or A\or B\or C\or D\or E\or F\fi}
\ifcase\@ptsize
 \font\tenmsb=msbm10
 \font\sevenmsb=msbm7
 \font\fivemsb=msbm5
\or
 \font\tenmsb=msbm10 scaled \magstephalf
 \font\sevenmsb=msbm7 scaled \magstephalf
 \font\fivemsb=msbm5  scaled \magstephalf
\or
 \font\tenmsb=msbm10 scaled \magstep1
 \font\sevenmsb=msbm7 scaled \magstep1
 \font\fivemsb=msbm5 scaled \magstep1
\fi
\newfam\msbfam
\textfont\msbfam=\tenmsb
\scriptfont\msbfam=\sevenmsb
\scriptscriptfont\msbfam=\fivemsb
\edef\msbfam@{\hexnumber@\msbfam}
\def\Bbb#1{\fam\msbfam\relax#1}
\def\widehat#1{\setboxz@h{$\m@th#1$}%
 \ifdim\wdz@>\tw@ em\mathaccent"0\msbfam@5B{#1}%
 \else\mathaccent"0362{#1}\fi}
\def\widetilde#1{\setboxz@h{$\m@th#1$}%
 \ifdim\wdz@>\tw@ em\mathaccent"0\msbfam@5D{#1}%
 \else\mathaccent"0365{#1}\fi}

\def\RIfM@{\relax\ifmmode}
\def\nonmatherr@#1{\errmessage{\string#1\space allowed only in math mode}}
\def\Bbb{\RIfM@\expandafter\Bbb@\else
 \expandafter\nonmatherr@\expandafter\Bbb\fi}
\def\Bbb@#1{{\Bbb@@{#1}}}
\def\Bbb@@#1{\fam\msbfam\relax#1}
\def\setboxz@h{\setbox\z@\hbox}
\def\wdz@{\wd\z@}
\catcode`\@=\csname pre amssym.def at\endcsname

\newcommand{\nl}{{\smallskip\noindent}}
\newcommand{\noi}{{\noindent}}

\newcommand{\dst}{\displaystyle}
\newcommand\ovl[1]{ \overline {#1} }
\renewcommand{\Im}{{\,\rm Im\,}}
\renewcommand{\Re}{{\,\rm Re\,}}

\renewcommand{\a}{\alpha}
\renewcommand{\b}{\beta}
\newcommand{\G}{\Gamma}
\newcommand{\g}{\gamma}
\newcommand{\D}{\Delta}
\renewcommand{\d}{\delta}
\renewcommand{\k}{\kappa}
\renewcommand{\L}{\Lambda}

\newcommand{\m}{\mu}
\newcommand{\n}{\nu}
\renewcommand{\P}{\Pi}
\newcommand{\p}{\pi}
\renewcommand{\r}{\rho}
\newcommand{\s}{\sigma}
\renewcommand{\t}{\tau}
\newcommand{\f}{\varphi}
\renewcommand{\O}{\Omega}
\renewcommand{\o}{\omega}
\newcommand{\cA}{{\rm A}}
\newcommand{\cB}{{\rm B}}
\newcommand{\cE}{{\rm E}}
\newcommand{\cR}{{\rm R}}
\newcommand{\cK}{{\rm K}}
\newcommand{\cC}{{\rm C}}
\newcommand{\cD}{{\rm D}}
\newcommand{\cN}{{\rm N}}
\newcommand{\cP}{{\rm P}}
\newcommand{\cJ}{{\rm J}}
\newcommand{\cT}{{\rm T}}
\newcommand{\cS}{{\rm S}}
\newcommand{\torus}{{\Bbb T}}
\renewcommand{\natural}{{\Bbb N}}
\newcommand{\real}{{\Bbb R}}
\newcommand{\integer}{{\Bbb Z}}
\newcommand{\complex}{{\Bbb C}}

\newcommand\ppo{{(0)}}
\newcommand\ppu{{(1)}}
\newcommand\ppd{{(2)}}
\newcommand\ppt{{(3)}}

\newcommand\ppn{{(n)}}
\newcommand\ppm{{(m)}}
\newcommand\pph{{(h)}}
\newcommand\ppi{{(i)}}
\newcommand\ppj{{(j)}}
\newcommand\ppk{{(k)}}
\newcommand\id{{\, \rm id \,}}
\newcommand\fg{{\g_1}}
\newcommand\tk{{ k }}
\newcommand\fk{{ k_1 }}
\newcommand\sk{{ k_2 }}
\newcommand\td{{n}}
\newcommand\meas{{\, \rm meas\,}}
\newcommand\pertnorm{{E}}
\newcommand\KAM{{\hat E}}

\newcommand{\ie}{{i.e.}}
\newcommand{\eg}{{e.g.}}

\newcommand{\cDel}{{\cD e\ell}}

\newtheorem{theorem}{Theorem}[chapter]
\newtheorem{lemma}[theorem]{Lemma}

\theoremstyle{definition}
\newtheorem{definition}[theorem]{Definition}

\theoremstyle{remark}
\newtheorem{remark}[theorem]{Remark}

\theoremstyle{proposition}
\newtheorem{proposition}[theorem]{Proposition}

\theoremstyle{corollary}
\newtheorem{corollary}[theorem]{Corollary}

\numberwithin{section}{chapter}
\numberwithin{equation}{chapter}

\makeindex

\begin{document}

\frontmatter

\title{Perihelia Reduction and  Global Kolmogorov Tori in the Planetary  Problem}


\author{Gabriella Pinzari}
\address{Dipartimento di Matematica ed Applicazioni ``R. Caccioppoli'', Universit\`a di Napoli ``Federico II", Monte Sant'Angelo -- Via Cinthia I-80126 Napoli (Italy)}
\email{gabriella.pinzari@unina.it}
\thanks{This research has been financially supported {partly (up to February 28, 2016)} by ERC Ideas-Project 306414 ``Hamiltonian PDEs and small divisor problems: a dynamical systems approach'' {and partly (since March 1, 2016) by the ERC Project 677793 ``Stable and Chaotic Motions in the Planetary Problem.''}}

\date{February 21, 2015}

\subjclass[2010]{Primary 34C20, 70F10,  37J10, 37J15, 37J40;
Secondary 34D10,  70F07, 70F15, 37J25, 37J35}

\keywords{Canonical coordinates; Jacobi's reduction; Deprit's reduction; Perihelia reduction; Symmetries; Quasi-periodic motions; Arnold's theorem on the stability of planetary motions.}

\dedicatory{
}

\maketitle

\tableofcontents

\begin{abstract}
We prove the existence of an almost full measure set of $(3n-2)$-dimensional {quasi-periodic} motions in the planetary problem with $(1+n)$ masses, with eccentricities arbitrarily close to the Levi--Civita limiting value and relatively high inclinations. {This extends previous results, where smallness of eccentricities and inclinations was assumed. The question had been previously considered by {V.~I.~Arnold} \cite[Ch III, \S 1, n. 6, p. 128]{arnold63} in the 60s, for the particular case of the planar three-body problem, where, due to the limited number of degrees of freedom, it was enough to use the invariance of the system by the SO(3) group}.\\
    The proof exploits nice parity properties of a new set of coordinates for the planetary problem, which reduces completely the number of degrees of freedom for the system (in particular, its degeneracy due to rotations) and, moreover, is well fitted to its reflection invariance. It allows the explicit construction of {an associated close to be integrable system,} replacing Birkhoff normal form, common tool of previous literature.
\end{abstract}

\mainmatter

 \chapter{Background and results} 

In recent years, substantial progress on a statement by Vladimir Igorevich Arnold concerning the stability of the planetary system has been achieved \cite{kolmogorov54, arnold63, laskarR95, robutel95, herman09, fejoz04, pinzari-th09, chierchiaPi11b}.

\nl It sounds as follows.

\vskip.1in \noi {\sl ``For the majority of initial conditions under which the instantaneous orbits of the planets are close to circles lying in a single plane, perturbation of the planets on one another produces, in the course of an infinite interval of time, little change on these orbits provided the masses of the planets are sufficiently small. {\rm [\,\ldots]} In particular {\rm[\,\ldots]} in the n-body problem there exists a set of initial conditions having a positive Lebesgue measure and such that, if the initial positions and velocities of the bodies belong to this set, the distances of the bodies from each other will remain perpetually bounded.''} {\cite[Chapter III, p. 125]{arnold63}}.

\vskip.1in \noi Solving the differential equations of the motions of the planetary problem, {\ie,} $n$ planets interacting among themselves and with a star via gravity is, for $n\ge 2$, a problem with ancient roots. {This story goes back to Sir Isaac Newton -- who brilliantly solved the case of two bodies and then, tackling the {analogous} one for three bodies, soon realized the necessity of turning to a ``perturbative'' study (except for naming it a ``head ache problem'') -- passed through investigations by eminent mathematicians like Delaunay, Lagrange, the prize publicly announced by {King} Oscar II of Sweden and Norway and awarded to Henri Poincar\'e, but its ``solution'' is nowadays open.} Chaotic and stable regions may coexist \cite{arnold63, FGKR14, DKRS15}.

\nl The question {received} a new mathematical description, and a strong modern endorsement, \mbox{after} A.~N.~Kolmogorov announced, at the International Congress of Mathematicians of 1954 {in} Amsterdam, what is now almost unanimously considered the most important result of the last century for dynamical systems: {the} theorem of conservation of the invariant torus. This breakthrough result, next {enriched by} substantial contributions by J. Moser and V. I. Arnold himself \cite{kolmogorov54, moser1962, arnold63c}, states that for a generic Hamiltonian system close to an integrable one, {\ie, a system of the form
\[ {\rm H}(I,\f)={\rm h}(I)+\m\, f(I,\f)\qquad (I,\f)\in B\times \torus^N\quad B\subset \real^N\ \quad \torus:=\real/(2\p\integer)\quad \m\ll1, \]}
the major part of unperturbed motions survives, after a small perturbation is switched on, provided suitable ``non-degeneracy'' conditions are verified by {the} ``unperturbed part'' ${\rm h}$. Moreover, the theory provides precise arithmetic (``diophantine'') properties to be verified by the ``unperturbed frequencies'' $\o_*=\partial{\rm h}(I_*)$, {so that} they will be preserved in the full system.\\
In 1962, V. I. Arnold, extending Kolmogorov's ideas, and looking for an application to the planetary problem, at the International Congress of Mathematicians of Stockholm, announced the theorem of stability of planetary motions quoted above. In 1965 Kolmogorov and Arnold were {awarded} the Lenin Prize for their studies on the stability of the planetary problem -- but the story was not finished there.

\nl In order to introduce the results of this paper, we highlight basic facts of this story and its continuation, referring the reader to \cite{fejoz13, chierchia13, pinzari13, chierchiaPi14, pinzari14} for more notices.

\nl The planetary problem is close to the integrable problem of $n$ uncoupled two-body problems, where each planet interacts separately with the sun. The mutual interactions among planets are regarded as a perturbing function, the smallness of which is ruled by the planets' masses. However, as a perturbed system, the planetary problem has a limiting degeneracy. Its associated integrable system (the two-body problem) is ``super-integrable'': it has more integrals than degrees of freedom. {At a technical level, the limiting degeneracy is {exhibited by} the disappearance of degrees of freedom in the unperturbed part. Therefore, continuing the unperturbed motions to a {\sl positive measure set} of quasi-periodic trajectories might, in general, be not possible, in absence of further informations on the perturbing function.}

\nl Arnold found, for the planetary problem, a brilliant solution to the problem of the limiting degeneracy. This {led} him to add {to the} assumptions and assertions that are proper of perturbation theories (\eg, ``the masses of the planets are sufficiently small'', ``set of initial conditions having a positive Lebesgue measure'', ``the distances \ldots will remain perpetually bounded'')
a further requirement of smallness of eccentricities and inclinations of the unperturbed Keplerian ellipses (``the instantaneous orbits of the planets are close to circles lying in a single plane''). Let's summarize Arnold's ideas.

 \noi Choosing, as Arnold did, Poincar\'e coordinates \cite{Poincare:1892} (see, also \cite[Ch. III, \S 2]{arnold63}, or, \eg, \cite{chierchiaPi11c, fejoz13b}), the system takes the usual {close to be integrable} form $${\rm H}_{\cP oi}={\rm h}_{\rm Kep}+\m\, f_{\cP oi},$$ where $\m$ is a small parameter related to the planetary masses, but the unperturbed ``Keplerian'' part ${\rm h}_{\rm Kep}(\L)$ depends on only $n$ action variables $\L=(\L_1,\cdots,\L_n)$ (related to the semi-major axes of the instantaneous Keplerian ellipses), out of an overall of $3n$ degrees of freedom. The perturbing function, $f_{\cP oi}$, on the other hand, depends on all the coordinates: the actions $\L$, their conjugated angles $\ell=(\ell_1,\cdots, \ell_n)$ (proportional to the areas of the elliptic sectors spanned by the planets), and, moreover, on some other coordinates $({\rm p}, {\rm q})=({\rm p}_1,\cdots, {\rm p}_{2n},{\rm q}_1,\cdots, {\rm q}_{2n} ) $, $4n$-dimensional, related to those (``secular'') quantities (eccentricities, inclinations, nodes and perihelia of the ellipses) that in the unperturbed problem stay fixed, and for this reason do not appear in ${\rm h}_{\rm Kep}$.

\nl It is of great help that the averaged perturbing function (with respect to the angles $\ell$) $\ovl{f_{\cP oi}}(\L, {\rm p}, {\rm q})$ enjoys several parities in the coordinates $({\rm p}, {\rm q})$, geometrically related to its invariance by rotations and reflections with respect to the coordinate planes. The ``secular origin'' $({\rm p}, {\rm q})=0$, corresponding to all the planets moving on co-centric circles in the same plane, turns out to be an elliptic equilibrium point for the averaged perturbing function, for any value of $\L$.

\nl Arnold brilliantly argued to exploit this circumstance to his purpose. By Birkhoff theory, one might think to switch to another set of canonical coordinates $(\L,\widetilde\ell, \widetilde{\rm p}, \widetilde{\rm q})$, analogous to Poincar\'e's coordinates, possibly defined only for $(\widetilde{\rm p}, \widetilde{\rm q})$ in a small neighborhood of radius $\varepsilon$ around the origin, such that the Hamiltonian of the system, or, more precisely, its $\widetilde\ell$-averaged (``secular'') perturbing function $\ovl{f_{\cB ir}}$, takes a ``normalized form'' : it is a polynomial, $\ovl{f_{\cB ir, tr}}$, of some degree greater or equal than two in the combinations (``degenerate actions'') $\t_i=\frac{\widetilde{\rm p}_i^2+\widetilde{\rm q}_i^2}{2}$, $i=1$, $\cdots$, $2n$, plus a remainder with a higher order. {Roughly, Arnold projected to solve the limiting degeneracy by conjugating the planetary system to a new system, whose unperturbed part was just the truncated, normalized Hamiltonian $${\rm h}_{\rm Kep}+\m\,\ovl{f_{\cB ir, tr}}$$ so as to recover the standard set up of KAM theory.} With these ideas in mind, he proved the following impressive result and next applied it to the planar three-body problem. It states that stable trajectories occupy a positive measure set of the phase space, and are more and more dense {closer} to the elliptic equilibrium. Hence, the smaller eccentricities and inclinations are, the larger the number of stable motions is.
 \vskip.1in
 \noi {\bf `The Fundamental Theorem'' (V. I. Arnold, \cite{arnold63})} {\it If the Hessian matrix of ${\rm h}$ and the matrix of the coefficients of the second-order term in $\t_i$ in $\ovl{f_{\cB ir}}$ (``torsion'', or ``second-order Birkhoff invariants'') do not vanish identically, and if $\m$ is suitably small with respect to $\varepsilon$, the system affords a positive measure set $\cK_{\m, \varepsilon}$ of quasi-periodic motions in phase space such that its density goes to $1$ as $\varepsilon\to 0$.}

\vskip.1in \noi Arnold perfectly knew that, in order to apply the Fundamental Theorem to the problem in space, one should previously treat an unpleasant fact:  {one} of the first order Birkhoff invariants vanishes identically. {He was aware that the {reason for} this first-order degeneracy was to be sought into the existence of two non-commuting integrals}, the two horizontal components of the total angular momentum of the systems. If, apparently, a vanishing eigenvalue strongly violates the construction of the normalized system (a deeper analysis of the symmetries of the perturbing function \cite{maligeRL02, chierchiaPi11c}, however, shows that the identically vanishing eigenvalue is not a real obstruction), a major problem definitely prevents the application of the Fundamental Theorem: an infinite number of coefficients of {\sl any order} of the (formal) Birkhoff series vanishes identically, {among which one entire row and a column in the torsion matrix}, which so is identically singular, and the reason is again the invariance by rotations. {The proof of this generalized degeneracy is in \cite{chierchiaPi11c}. We recall here that even Herman had raised a question about the degeneracy of torsion \cite[p. 24]{herman09}. }

\nl {We do not know weather Arnold was aware of the infinite degeneracy of the normalized system (he did not even mention the vanishing of torsion in his paper). He however suggested two different strategies for the three- and the many-body {case, of} which he provided very few and somewhat controversial details. As for the three-body problem (his ideas for the many-body case will be recalled a few below), he proposed} to reduce the integrals (hence, the number of degrees of freedom) of the system by switching to a system of canonical coordinates going back to the XIX century, worked out by Jacobi and Radau \cite{jacobi1842, radau1868}, which in literature go under the name of {\sl Jacobi reduction of the nodes}. The idea was later completely developed by P. Robutel \cite{robutel95}, who, in a deeply quantitative study, checked the non-degeneracy assumptions required by the Fundamental Theorem.

\nl Finding a system of canonical coordinates that do the job of Jacobi reduction of the nodes when the number of bodies is more than three has been a central difficulty for a long time \cite{arnold63, maligeRL02}. At this respect, {Arnold sadly commented: ``{In the case of more than three bodies there is no such elegant method [as Jacobi reduction of the nodes] of reducing the number of degrees of freedom.}}'' \cite[Ch. III, \S 5.5, p. 141]{arnold63}.

\nl Exactly twenty years later, F. Boigey and A. Deprit refuted this sentence \cite{boigey82, deprit83}. They indeed were able to extend Jacobi-Radau reduction to the four{-body and} general problem, respectively. It should be remarked, anyway, that, while the works by Jacobi, Radau and Boigey provide canonical coordinates on suitable sub-manifolds of the phase space, the one by Deprit is more general and clarifying, since {it} provides a set of canonical coordinates for the whole phase space and allows us to recover his predecessors by restriction.
 
\nl The utility of Boigey-Deprit's coordinates was not suddenly clear. {Neither} Boigey nor Deprit ever provided any motivation of their study, or foresaw applications. The only application that is known to the author up to 2008, concerning indeed Deprit's coordinates, stands in a paper by Ferrer and Os\'acar, in the 90s, to the three body problem \cite{ferrerO94}. But this case is not really exhaustive, since for three bodies Deprit's and Jacobi-Radau's coordinates coincide. A reason why Boigey-Deprit's coordinates have been forgotten so long might be that, for more than three bodies, they actually have a less natural aspect, compared to the classical case of Jacobi. A sort of ``hierarchical'' structure in the geometry of Deprit's coordinates discouraged the author himself, who, at the end of his paper, declared: ``Whether the new phase variables are practical in the general theory of perturbation is an open question. At least, for planetary theories, the answer is likely to be in the negative. But finding a natural system of coordinates for eliminating the nodes in a planetary cluster was not the intention of this note.'' \cite[p. 194]{deprit83}.

\nl In the meantime, in 2004, the first general proof of Arnold's stability statement appeared. It was by Jacques F\'ejoz, who completed investigations by the late Michael Herman \cite{fejoz04} -- but the different procedure that Herman had in mind did not {rely on} the necessity of handling, explicitly, good coordinates. Indeed, Herman conceived a proof based, besides on a ``twist-less'' KAM theory going back to H. Russmann \cite{Russmann:2001}, on indirect arguments of Lagrangian intersections in order to {bypass} the so-called ``secular resonances''. See \cite{chierchiaPi14} for more details.

\nl In 2008, Boigey-Deprit's coordinates were rediscovered by the author \cite{pinzari-th09}, in a slightly different, ``planetary'' form. {The rediscovery was motivated by the purpose of realizing Arnold's program (\ie, applying the Fundamental Theorem quoted above directly to the planetary Hamiltonian) in the general case, so as to obtain a detailed information about the tori frequencies, the measure of the invariant set and the symplectic structure of the phase space.} The utility of Boigey-Deprit's coordinates became suddenly clear: switching (in order to overcome certain singularities of the chart) to a regularized version, called ``RPS'' coordinates, (acronym standing for ``Regular, Planetary and Symplectic''), allowed {them} to derive the Birkhoff normal form of the planetary problem, to prove its non-degeneracy, and hence to complete the application of the Fundamental Theorem to the general problem. These results have been published in \cite{chierchiaPi10, chierchiaPi11a, chierchiaPi11b}.

\nl Qualitatively, RPS coordinates are very different from JRBD (Jacobi-Radau-Boigey-Deprit){; rather, they are} more similar to Poincar\'e coordinates. The mentioned parities and the elliptic equilibrium of the averaged system are still present in the RPS-averaged system. But, as an advantage with respect to Poincar\'e coordinates, the RPS perform\footnote{In the framework of the study of canonical coordinates for the planetary system, by ``partial reduction'', we mean a system of canonical coordinates where a couple of conjugated coordinates consists of integrals (\eg, functions of the three components of the total angular momentum). By ``full reduction'', we mean a partial reduction where also another integral appears among the coordinates. The terms ``partial reduction'', ``full reduction'' have been coined in \cite{maligeRL02}.} a ``partial reduction'' of the rotation symmetry -- {in contrast} with JRBD coordinates, which reduce ``fully''. This way, all the degeneracies of the Birkhoff series mentioned above are removed at once, and the non-degeneracy assumptions of the Fundamental Theorem may be checked. \\
{We like to recall now Arnold's strategy for the many-body case: more than forty years earlier, he foresaw to construct a system of coordinates analogous to RPS, via a Taylor series in Poincar\'e coordinates \cite[Ch III, \S 5, n. 5, p. 141]{arnold63}. \\
Indeed, both the reduction of the nodes and this latter reduction are available whatever the number of bodies {is.}}

\vskip.1in \nl The possibility of switching from Delaunay-Poincar\'e to the more fruitful JRBD, or even RPS coordinates, is an effect of the limiting degeneracy. This gives in fact the opportunity of remixing coordinates related to secular quantities, and, simultaneously, keeping the Keplerian term ${\rm h}_{\rm Kep}$ unvaried. 

\nl Following this idea, in this paper, we show that other systems of coordinates may be determined for the planetary problem which, as well as JRBD and RPS coordinates, are well adapted to overcome the degeneracy due to rotations, and, moreover, enjoy some different properties.

\nl We present a full reduction, which we call $\cP$-map, or perihelia reduction. It refines JRBD coordinates in two respects.

\nl Firstly, the $\cP$-map is well defined in the case of the planar problem, while JRBD coordinates are not. Everyone knows, in fact, that the starting point for the Radau-Jacobi reduction is the so-called ``line of the nodes'', the straight line determined by the intersection between the planes of the two orbits. When the orbits of the two planets belong to the same plane, this is not defined. A similar circumstance arises for Boigey-Deprit's coordinates, since their construction relies on certain straight lines in the space, which again {lose} their meaning in case of co-planarity. 

\nl The proof of Arnold's theorem given in \cite{pinzari-th09, chierchiaPi11b} is not affected by such singularity, since, as said, it relies on RPS coordinates, which, {at the expense of} one more degree of freedom, are well defined for co-planar motions -- in that case they reduce to the classical Poincar\'e coordinates.

\nl It has its consequences when one wants to compare results for the fully reduced systems, in space or in the plane. The singularity of the chart does not allow {one} to state that motions in the spatial problem with minimum number of independent frequencies starting with very small inclinations stay close to the corresponding planar motions. Notwithstanding further studies appearing in \cite{pinzari13}, where this problem is partially overcome (via the construction of regular coordinates for co-planar motions defined locally), it would be nice, in principle, to handle a global system of action-angle coordinates which completely reduces rotations and is shared simultaneously by the planar and the spatial problem.

\nl Secondly, the $\cP$-map is well adapted to reflection symmetries of the problem, while JRBD coordinates are not, as discussed in \cite{maligeRL02, pinzari14}.

\nl Reflection symmetries are parities of the Hamiltonian expressed in Cartesian coordinates. As known, this does not change under arbitrary changes of the signs of positions or momenta coordinates. They are not related to integrals. Therefore, it might be a nice fact, and in general useful for applications, to have a system of coordinates that, after integrals are reduced, parities associated to reflections are maintained. Quite often parities are associated to equilibria, and equilibria to stable motions; an example is provided a few lines below.

\nl We shall apply {the} $\cP$-map by proving a variant of Arnold's stability theorem. We shall {face up to} a question raised again by Arnold in his fantastic paper on the possibility of removing the constraint on eccentricities and inclinations. He indeed proved that, at least for the planar three-body problem, there is no need to assume their smallness. Rather, it is sufficient that the trajectories of the planets are away enough so as to avoid collisions. He obtained this stronger result by exploiting the convergence of the Birkhoff series associated to the averaged perturbation, a very particular and happy circumstance, due to the few degrees of freedom of the problem.

\nl From the mathematical point of view, the question is whether{strategies exist for finding stable motions other than} the one of exploring the neighborhood of the elliptic equilibrium.

\nl Concerning instead the physical relevance, asteroids or some trans-Neptunian objects have motions with relatively large eccentricities and inclinations and an almost continuous spectrum of frequencies.

\nl Besides the mentioned stronger result by Arnold, some other statements in the same direction have been obtained for the case of the spatial three-body problem and the planar problem with any number of bodies \cite{pinzari13}. Here, the measure of the invariant set has been estimated to be larger and larger as the planetary masses and the semi-axes ratios are small, no matter the smallness {of the} eccentricities and inclinations -- the proof {relies} on an argument of convergence of a {significant} approximation of the Birkhoff series. Other results in this direction have been announced by J. F\'ejoz, since late 2013 \cite{fejoz15}. 

\nl Even though the arguments of \cite{arnold63, pinzari13} do not apply to the general spatial problem, since no {significant} approximation of the Birkhoff series associated to the averaged perturbation is integrable, using the $\cP$-map, we shall prove the following.

\vskip.1in \noi {\bf Theorem A} {\it Fix numbers $0<\underline e_i<\ovl e_i<0.6627\ldots$, $i=1,\cdots,n$. There exists a number ${\rm N}$ depending only on $n$ and a number $\a_0$ depending on $\underline e_i$, $\ovl e_i${, and} $n$ such that, if $\a<\a_0$, $\m\le \a^{\rm N}$, in a domain of planetary motions where the semi-major axes $a_1<a_2<\cdots<a_n$ are spaced as follows 
\[ a_i^-\le a_i\le a_i^+\qquad {\rm with}\qquad a_{i}^\pm:= \frac{a_1^\pm}{\a^{\frac{1}{3}(2^{n+1}-2^{n-i+2}+1-i)}} \tag{$*$} \]
there exists a positive measure set $\cK_{\m, \a}$, the density of which in phase space can be bounded below as $${\rm dens}(\cK_{\m, \a})\ge 1-(\log\a^{-1})^{\rm p}\sqrt\a,$$ consisting of quasi-periodic motions with $3n-2$ frequencies where the planets' eccentricities $e_i$ verify $$\underline e_i\le e_i\le \ovl e_i.$$ }

\noi Before we switch to details, a few remarks. 

\nl Firstly, the claimed upper bound $0.6627\ldots$ is classical. It is related to the fact that, as well as in \cite{arnold63, pinzari13}, the proof uses the machinery of real-analytic functions. We refer the reader to \cite{tisserand1889, leviCivita1904} and references therein for general notices. A treatment of the argument, as needed in the present paper, is provided in Section \ref{Computing the domain of holomorphy}.

\nl Secondly, as it may be seen to the choice of $a_j^\pm$, the distances among the planets' semi-axes are not of the same order but grow super-exponentially going towards the sun. This resembles a sort of belt arrangement, observed in nature for asteroids. It is possible to prove an analogous result, with increasing distances in the opposite direction.

\nl Thirdly, the result in Theorem A (especially, the claimed growth of $a_i^\pm$) may be regarded as an alternative way of solving the problem of the limiting degeneracy -- without Birkhoff normal form.

\vskip.1in \noi {\bf Acknowledgments} I am indebted to Jacques F\'ejoz, who let me know the work by Harrington~\cite{harrington69}, without which {I would never have thought of} this application of the $\cP$-coordinates. Also, I am deeply grateful to A.~Celletti and R.~de la~Llave, for their interest and for encouraging me with precious advices. Thanks finally to L.~Biasco for his interest.
  \chapter{Kepler maps and the Perihelia reduction} \label{perihelia reduction}

\nl We introduce the {\sl Perihelia reduction}, or $\cP$-map, in the slightly general context of {\sl Kepler maps}.

\vskip.1in \noi Fix a reference frame ${\rm G}_0=(k^\ppu, k^\ppd, k^\ppt)$ in the Euclidean space $E^3$. We identify the three chosen directions $k^\ppu$, $k^\ppd$, $k^\ppt$ with the triples of coordinates with respect of the system of coordinates established by themselves: $$k^\ppu=\left(
\begin{array}
    {ccc} 1\\
    0\\
    0
\end{array}
\right)\qquad k^\ppd=\left(
\begin{array}
    {ccc} 0\\
    1\\
    0
\end{array}
\right)\qquad k^\ppt=\left(
\begin{array}
    {ccc} 0\\
    0\\
    1
\end{array}
\right). $$
\begin{definition}
    \rm An {\sl ellipse} (with a focus in the origin and non-vanishing eccentricity) is a quadruplet ${\eufm E}=(a, e, N, P)$, where $a\in\real_+$ is the {\sl semi-major axis}, $e\in (0,1)$ is the {\sl eccentricity}, $N\in \real^3\cap S^2$ is the {\sl normal direction} and $P\in N^\perp\cap S^2$ is the {\sl perihelion direction}.
\end{definition}

\begin{definition}
    [Kepler maps]\label{Kepler map}\rm Given $2n$ positive ``mass parameters'' ${\eufm m}_1$, $\cdots$, ${\eufm m}_n$, ${\eufm M}_1$, $\cdots$, ${\eufm M}_n$, a set ${\eufm X}\subset\real^{5n}$, we say that a map
    \[ \cK: \quad {{\rm K}}=({\rm X}_\cK, \ell)\in \cD:={\eufm X}\times \torus^n\to (y_\cK,x_\cK) \in \cC:=\cK(\cD)\subset (\real^{3})^n\times (\real^{3})^n \]
     where 
    \[
    \begin{array}
        {lll} \ell=(\ell_1,\cdots, \ell_n),\qquad &(y_\cK,x_\cK)=(y_\cK^\ppu, \cdots, y_\cK^\ppn, x_\cK^\ppu, \cdots, x_\cK^\ppn)\\
        \\
        \dst y_\cK^\ppj=y_\cK^\ppj({\rm X}_\cK, \ell_j) &x_\cK^\ppj=x_\cK^\ppj({\rm X}_\cK, \ell_j)\qquad j=1,\cdots n,
    \end{array}
    \]
    is a {\sl Kepler map} if there exists an injection \beqano \t_\cK: \quad {\rm X}_\cK\in {\eufm X}&\to& {\eufm E}_\cK=\big({\eufm E}_{1,\cK},\cdots,{\eufm E}_{n,\cK}\big)\quad \eeqano which assigns to any ${\rm X}_\cK\in{\eufm X}$ an n-plet $\big({\eufm E}_{1,\cK},\cdots,{\eufm E}_{n,\cK}\big)$ of (co-focal) ellipses $$ {\eufm E}_{j,\cK}=\big(a_{j,\cK}, e_{j,\cK}, N^\ppj_\cK, P^\ppj_\cK\big),\quad j=1,\ \cdots,\ n $$ and $\cK$ acts in the following way. Letting $Q_\cK^\ppj:=N_\cK^\ppj\times P_\cK^\ppj$, then \beq{xjyj}x^\ppj_\cK={\rm a}_{j,\cK} P_\cK^\ppj+{\rm b}_{j,\cK} Q_\cK^\ppj\qquad y_\cK^\ppj={\rm a}^\circ_{j,\cK} P_\cK^\ppj+{\rm b}_{j,\cK}^\circ Q_\cK^\ppj\eeq where, if $\zeta_{j,\cK}$, the {\sl eccentric anomaly}, is the solution of Kepler's Equation \beq{Kep Eq}\zeta_{j,\cK}-e_{j,\cK}\sin\zeta_{j,\cK}=\ell_j\eeq then 
    \begin{equation}\label{a e b}
    \begin{split}
&{\rm a}_{j,\cK}:=a_{j,\cK}\big( \cos\zeta_{j,\cK}-e_{j,\cK} \big)\qquad \qquad\qquad
  \dst {\rm b}_{j,\cK}:=a_{j,\cK}\sqrt{1-e_{j,\cK}^2}\sin\zeta_{j,\cK}\\
&{\rm a}^\circ_{j,\cK}:=-{\eufm m}_j\sqrt{\frac{{\eufm M}_j}{a_{j,\cK}}} \frac{\sin\zeta_{j,\cK}}{1-e_{j,\cK}\cos\zeta_{j,\cK}} 
\quad\qquad {\rm b}_{j,\cK}^\circ:={\eufm m}_j\sqrt{\frac{{\eufm M}_j(1-e_{j,\cK}^2)}{a_{j,\cK}}}\frac{\cos\zeta_{j,\cK}}{{1-e_{j,\cK}\cos\zeta_{j,\cK}}}.
 \end{split}
    \end{equation}\end{definition}

\begin{remark}
    \rm The definition implies that 
    \begin{enumerate}
    \item[{\rm (i)}] $\cK$ is a bijection of the sets $\cD$ and $\cC$; \item[{\rm (ii)}] the {\sl angular momenta} and the {\sl energies}\footnote{Here, $\|v\|:=\sqrt{v_1^2+v_2^2+v_3^2}$ denotes the usual Euclidean norm of $v=(v_1, v_2, v_3)\in \real^3$.}%
     \beq{1B}{\rm C}^\ppj_\cK:=x^\ppj_{\cK}\times y^\ppj_\cK,\qquad {\rm H}_{\cK}^\ppj:=\frac{\|y^\ppj_\cK\|^2}{2{\eufm m}_j}-\frac{{\eufm m}_j{\eufm M}_j}{\|x^\ppj_\cK\|} \eeq 
    do not depend on $\ell_j$ and are given by \beq{1B*}{\rm C}^\ppj_\cK={\eufm m}_j\sqrt{{\eufm M}_j a_{j,\cK}(1-e_{j,\cK}^2)} N^\ppj_\cK,\quad {\rm H}_{\cK}^\ppj=-\frac{{\eufm m}_j{\eufm M}_j}{2a_{j,\cK}};\eeq

\item[{\rm (iii)}] the couples $(y_\cK^\ppj, x_\cK^\ppj)$ verify the system of ODEs \beqa{Kepler motion equations} \arr{\dst{\eufm m}_j \sqrt{\frac{{\eufm M}_j}{a_{j,\cK}^3}}\partial_{\ell_j} x_\cK^\ppj= y_\cK^\ppj\\
    \\
    \dst \sqrt{\frac{{\eufm M}_j}{a_{j,\cK}^3}}\partial_{\ell_j} y_\cK^\ppj=-{\eufm m}_j{\eufm M}_j \frac{x_\cK^\ppj}{\|x_\cK^\ppj\|^3}. } \eeqa \item[{\rm (iv)}] Even though {\sl canonical} maps (with respect to the standard two-form) have a pre-eminent role in Hamiltonian Mechanics, Kepler maps are used also in different contexts in Astronomy, where {being canonical} is not required. For example, one can consider the Kepler map associated to the ``elliptic elements'' injection 
\[ \t_{\cE e\ell\ell}:\quad (a, e, P, i, \Omega)\to {\eufm E}_{\cE e\ell\ell} \]
where $a=(a_1, \cdots, a_n)$ are the {\sl semi-major axes}, $e=(e_1, \cdots, e_n)$ are the {\sl eccentricities}, \mbox{$P=(P^\ppu, \cdots, P^\ppn)$} are the {\sl perihelia}, $i=(i_1, \cdots, i_n)$ are the {\sl inclinations}, \mbox{$\O=(\O_1, \cdots, \O_n)$} are the {\sl nodes' {longitudes}}.

\nl The only known examples up to now of {\sl canonical Kepler maps} are the classical {\sl Delaunay map} $\cDel$ (its definition is recalled in the next Definition \ref{Def: Delaunay}) and the map $\cD ep$ \cite{pinzari-th09, chierchiaPi11a} related to Deprit's coordinates \cite{deprit83}, which is recalled in Appendix \ref{comparison}. Below, we introduce a new canonical Kepler map.
\end{enumerate}
\end{remark}

\begin{definition}
    [perihelia reduction, or $\cP$-map]\label{def: P map}\rm We denote as $\cP$, and call {\sl perihelia reduction}, or $\cP$-map, the Kepler map \beq{cal P map}\cP:\quad {\rm P}= ({\rm X}_\cP,\ell) \in\cD_\cP={\eufm X}_\cP\times \torus^n\to (y, x)\in\real^{3n}\times \real^{3n}\eeq associated to the bijection $$\t_{\cP}:\quad {\rm X}_\cP=(\Theta,\chi,\L,\vartheta, \k)\in {\eufm X}_\cP \to({\eufm E}_1,\cdots, {\eufm E}_n)\in{\rm E }_\cP=\t_{\cP}({\eufm X }_\cP)\subset E^{3n} $$ defined by means of Definition \ref{t-1} and Proposition \ref{inversion} below.
\end{definition}

\begin{definition}
    \label{t-1}\rm

For a given $({\eufm E}_1,\cdots, {\eufm E}_n)\subset E^{3}\times \cdots\times E^3$, with ${\eufm E}_j=(a_j, e_j, N^\ppj, P^\ppj)$, and masses ${\eufm m}_1$, $\cdots$, ${\eufm m}_n$, ${\eufm M}_1$, $\cdots$, ${\eufm M}_n$, define \beqa{Cj Sj} {\rm C}_\cE^\ppj:={\eufm m}_j\sqrt{{\eufm M}_ja_j(1-e_j^2)}N^\ppj\qquad {\rm S}_\cE^\ppj:=\sum_{i=j}^n{\rm C}_\cE^\ppi\qquad 1\le j\le n\ \eeqa be the {\sl angular momenta} associated to ${\eufm E}_j$ and the $j^{\rm th}$ {\sl partial angular momenta}, so that \beq{former and latter}{\rm S}_\cE^\ppu=\sum_{i=1}^n{\rm C}_\cE^\ppi\qquad {\rm S}_\cE^\ppn={\rm C}_\cE^\ppn\eeq are the {\sl total angular momentum} and the angular momentum of the last ellipse, respectively. Define the {\sl$\cP$-nodes} \begin{equation}\label{good nodes}\begin{split} \n_j:=\left\{
    \begin{array}
        {llll} \dst k^\ppt\times {\rm S}_\cE^\ppu\quad &j=1\\
        \\
        \dst P^{(j-1)}\times {\rm S}_\cE^{(j)} &j=2,\cdots, n
    \end{array}
    \right.\qquad {\rm n}_j:=\dst{\rm S}_\cE^\ppj\times P^\ppj \qquad j=1,\cdots, n. \end{split}\end{equation} Finally, define $$\cE_\cP:=\big\{ (({\eufm E}_1,\cdots, {\eufm E}_n)\subset E^{3}\times \cdots\times E^3):\quad 0<e_j<1,\quad \n_j\ne 0\quad {\rm n}_j\ne 0\quad \forall\ j=1,\cdots, n\big\},$$ and, on this set, the map $$ \t_\cP^{-1}:\quad ({\eufm E}_1,\cdots, {\eufm E}_n)\in \cE_\cP\to {\rm X}_\cP\in{\eufm X}_\cP= \t_\cP^{-1}(\cE_\cP) $$ where $${\rm X}_\cP=(\Theta, \chi, \L, \vartheta, \k)\in \real^n\times \real_+^n\times\real_+^n\times \torus^n\times \torus^n$$ with \beqano &&\Theta=(\Theta_0,\cdots,\Theta_{n-1}),\quad \vartheta=(\vartheta_0,\cdots,\vartheta_{n-1})\nonumber\\
    \nonumber\\
    &&\chi=(\chi_0,\cdots,\chi_{n-1}),\quad\ \k=(\k_0,\cdots,\k_{n-1})\nonumber\\
    \nonumber\\
    &&\L=(\L_1,\cdots,\L_n) \eeqano defined via the following formulae{:}
\begin{equation}\label{belle*}\begin{split}
   &  \Theta_{j-1}:=\left\{
        \begin{array}
            {lrrr} \dst Z:={\rm S}_\cE^\ppu\cdot k^\ppt \\
            \\
            \dst {\rm S}_\cE^{(j)}\cdot P^{(j-1)}
        \end{array}
        \right.
        \qquad\ \
         \vartheta_{j-1}:=\left\{
        \begin{array}
            {lrrr} \dst\zeta:=\a_{k^\ppt}(k^\ppu, \n_1)\qquad& j=1\\
            \\
            \dst \a_{P^{(j-1)}}({\rm n}_{j-1}, \n_{j})&2\le j\le n
        \end{array}
        \right.\\\\
        &\chi_{j-1}:=\left\{
        \begin{array}
            {lrrr}{\rm G}:=\|{\rm S}_\cE^\ppu\| \\
            \\
            \|{\rm S}_\cE^\ppj\|
        \end{array}
        \right. 
        \qquad\qquad \k_{j-1}:=\left\{
        \begin{array}
            {lrrr}{\eufm g}:=\a_{{\rm S}_\cE^\ppu}(\n_1, {\rm n}_1)\qquad& j=1\\
            \\
            \a_{{\rm S}_\cE^\ppj}(\n_j, {\rm n}_j)&2\le j\le n
        \end{array}
        \right.\\\\
            &\L_j:={\eufm M}_j\sqrt{{\eufm m}_j a_j}.
               \end{split}
    \end{equation}\end{definition}

\begin{proposition}
    \label{inversion} Let ${\eufm X}_\cP$ be the subset of $\real^n\times \real_+^n\times\real_+^n\times \torus^n\times \torus^n$ defined by the following inequalities \beqa{non zero or one eccentricities} &&\sqrt{\chi_{i-1}^2+\chi_{i}^2-2\Theta_{i}^2+2\sqrt{(\chi_{i}^2-\Theta_{i}^2)(\chi_{i-1}^2-\Theta_{i}^2)}\cos{\vartheta_{i}}}<\L_i\nonumber\\
    && (\chi_{i-1}-\chi_i, \vartheta_i)\ne (0,\p)\qquad 0<\chi_{n-1}<\L_n\quad i=1,\cdots, n-1 \eeqa and \beq{non parallelism}|\Theta_0|< \chi_0\qquad |\Theta_i|< \min(\chi_{i-1}, \chi_{i})\qquad i=1,\cdots, n-1.\eeq The map $\t_\cP^{-1}$ is a bijection of $\cE_\cP$ onto ${\eufm X}_\cP$. The formulae of the inverse map $$\t_\cP:\quad {\rm X}_\cP=(\Theta, \chi, \L,\vartheta, \k)\in\cD_\cP\to{\eufm E}_\cP= ({\eufm E}_{1,\cP},\cdots, {\eufm E}_{n,\cP})\in\cE_\cP\quad {\eufm E}_{j,\cP}=(a_{j,\cP}, e_{j,\cP}, N_\cP^\ppj, P_\cP^\ppj)$$ are as follows. Let $\iota_1$, $\cdots$, $\iota_n$, ${\rm i}_1$, $\cdots$, ${\rm i}_{n}\in (0,\p)$ be defined via \beq{good incli*}\cos\iota_{j}=\frac{\Theta_{j-1}}{\chi_{j-1}},\quad \cos{\rm i}_j:=\frac{\Theta_{j}}{\chi_{j-1}},\quad1\le j\le n\eeq $($with $\Theta_n:=0$, so that ${\rm i}_{n}=\frac{\p}{2})$ and $\cT_1$, $\cdots$, $\cT_n$, $\cS_1$, $\cdots$, $\cS_{n}\in {\rm SO}(3)$ via \beqa{CC} && {\rm T}_{j}:=\cR_3(\vartheta_{j})\cR_1(\iota_{j})\qquad {\rm S}_j:=\cR_3(\k_j)\cR_1({\rm i}_j),\quad 1\le j\le n\ \eeqa and let \beqa{P*ecc} {\rm C}_\cP^\ppj:={\rm T}_1{\rm S}_1\cdots {\rm T}_{j-1}{\rm S}_{j-1}{\rm T}_j\big(\chi_{j-1} k^\ppt-\chi_{j}{\rm S}_{j}{\rm T}_{j+1} k^\ppt\big) \eeqa with $\chi_n:=0$, so that 
\begin{equation}\label{CP}\begin{split}
\|{\rm C}_\cP^\ppj\|=\left\{
    \begin{array}{llll} \dst\sqrt{\chi_{j-1}^2+\chi_{j}^2-2\Theta_{j}^2+2\sqrt{(\chi_{j}^2-\Theta_{j}^2)(\chi_{j-1}^2-\Theta_{j}^2)}\cos{\vartheta_{j}}}\quad &j=1,\cdots, n-1\\\\
            \chi_{n-1}& j=n.
    \end{array}
    \right.
    \end{split}
\end{equation} Then ${\rm C}^\ppj_\cP={\rm C}^\ppj_\cE\circ \t_\cP$ and \beqa{P*ecc new}
    \begin{array}
        {lllll} \dst a_{j, \cP}=\frac{1}{{\eufm M}_j}(\frac{\L_j}{{\eufm m}_j})^2\qquad &\dst e_{j,\cP}=\sqrt{1-\frac{\|{\rm C}_\cP^\ppj\|^2}{\L_j^2}}\\
        \\ \dst N_\cP^\ppj=\frac{{\rm C}_\cP^\ppj}{\|{\rm C}_\cP^\ppj\|} &\dst P_\cP^\ppj={\rm T}_1{\rm S}_1\cdots {\rm T}_{j}{\rm S}_{j}k^\ppt.
    \end{array}
\eeqa
\end{proposition}

\begin{remark}
    \label{C perp to P}\rm  ~\\[-1em]
    \begin{enumerate}
    \item[{\rm (i)}] From ${\rm C}^\ppj_\cP={\rm C}^\ppj_\cE\circ\t_\cP$, \equ{1B}, \equ{1B*} and \equ{S and C}, there follows that ${\rm C}^\ppj_\cP=x^\ppj_\cP\times y^\ppj_\cP$. \item[{\rm (ii)}]$P_\cP^\ppj\perp N_\cP^\ppj$. Indeed, using the definitions, \beqano {\rm C}_{\cP}^\ppj\cdot P_{\cP}^\ppj&=&\chi_{j-1}k^\ppt\cdot\big({\rm S}_jk^\ppt\big)-\cT_{j+1} \chi_{j}k^\ppt\cdot\big(k^\ppt\big) \nonumber\\
    &=&\chi_{j-1}\cos\iota_j-\chi_{j}\cos{\rm i}_{j+1} \\&=& 0. \eeqano \item[{\rm (iii)}] ${\rm S}^\ppj_\cP:={\rm S}^\ppj_\cE\circ\t_\cP=\sum_{i=j}^n{\rm C}^\ppi_\cP=\chi_{j-1}{\rm T}_1{\rm S}_1\cdots {\rm T}_{j-1}{\rm S}_{j-1}{\rm T}_j k^\ppt$.
    \end{enumerate}
\end{remark}
\nl We shall prove that
\begin{theorem}
    \label{main} The $\cP$-map preserves the standard 2-form $$\sum_{j=1}^n dy^\ppj_\cP\wedge dx^\ppj_\cP=\sum_{i=1}^{n} \big(d \Theta_{i-1}\wedge d\vartheta_{i-1}+ d\chi_{i-1}\wedge d\k_{i-1}+ d\L_i\wedge d\ell_i\big).$$
\end{theorem}

\begin{remark}
    \rm Actually, we shall prove a finer result: the change ${\phi_{\cDel}^\cP}:= \cDel^{-1}\circ\cP$ which relates the $\cP$-coordinates to the classical Delaunay coordinates (see the Definition \ref{Def: Delaunay}) is homogeneous-canonical (compare Lemma \ref{homogeneous with respect to Del}).
\end{remark}
\vskip.in\textsc{Proof of  Proposition \ref{inversion}.} The formula for $a_{j,\cP}$ in \equ{P*ecc} is immediate from the definition of $\L_j$. Postponing to below that ${\rm C}^\ppj_\cP:={\rm C}^\ppj_\cE\circ\t_\cP$ has the expression in \equ{P*ecc} ({in} turn this implies \equ{CP}, the formula for $N^\ppj$ and the one for $e_{j,\cP}$ in \equ{P*ecc new}), we check that the image set $\t^{-1}_\cP(\cE_\cP)$ is included in the domain ${\eufm X}_\cP$ defined by inequalities \equ{non zero or one eccentricities}, \equ{non parallelism}. From the formula for $e_{j, \cP}$ in \equ{P*ecc new}, we have that conditions $0<e_{j,\cP}< 1$ for all $j=1, \cdots, n$ {correspond} to relations in \equ{non zero or one eccentricities}. Note that the first condition in the second line of \equ{non zero or one eccentricities} is equivalent to $e_{j,\cP}\ne 1$, as one sees rewriting \begin{equation}\label{CPjsq} \|{\rm C}_{\cP}^\ppj\|^2= \big(\sqrt{\chi_{j-1}^2-\Theta_{j}^2}-\sqrt{\chi_{j}^2-\Theta_{j}^2}\big)^2+2\sqrt{(\chi_{j}^2-\Theta_{j}^2)(\chi_{j-1}^2-\Theta_{j}^2)}(1+\cos{\vartheta_{j}}). \end{equation}

\nl Next, recalling the definitions of $\Theta_0$, $\chi_0$ in \equ{belle*}, and noticing the relations $$\Theta_j={\rm S}_\cE^{(j+1)}\cdot P^{(j)}=({\rm S}_\cE^{(j)}-{\rm C}_\cE^{(j)})\cdot P^{(j)}={\rm S}_\cE^{(j)}\cdot P^{(j)}\qquad j=1,\cdots, n-1,$$ we immediately see that conditions $\n_i\ne 0\ne {\rm n}_i$ imply \equ{non parallelism}. We have so checked what we wanted.

\nl Now it remains to check the formula for ${\rm C}_\cP^\ppj$ in \equ{P*ecc} and the one for $P_\cP^\ppj$ in \equ{P*ecc new}, for any ${\rm X}_\cP\in{\eufm X}_\cP$. To this end, we consider the following chain of vectors
 \begin{equation}\label{chain}
\begin{split}
\begin{array}
    {cccccccccccccccccc} \dst k^\ppt&\to&{\rm S}_\cE^\ppu&\to& P^\ppu&\to&\cdots&\to& {\rm S}_\cE^\ppj&\to& P^\ppj&\to&
    \cdots&\to& P^\ppn\\
    \\
    \dst&&\Downarrow&&\Downarrow&&\vdots&&\Downarrow&&\Downarrow&&
    \vdots&&\Downarrow\\
    \\
    \dst&&\n_1&&{\rm n}_1&&\vdots&&\n_j&&{\rm n}_j&&\vdots
    &&{\rm n}_n\\
\end{array}
\end{split}
\end{equation}
where $\n_1, {\rm n}_1, \cdots, \n_n, {\rm n}_n$ are the {\sl $\cP$-nodes} in \equ{good nodes}, given by the skew-product of the two consecutive vectors in the chain.

\nl We associate to this chain of vectors the following chain of frames
 \begin{equation}\label{P chain}
 \begin{split}
\begin{array}
    {cccccccccccccccccc} \dst {\rm G}_0&\to&{\rm F}_1&\to& {\rm G}_1&\to&\cdots&\to& {\rm F}_j&\to& {\rm G}_j&\to& {\rm F}_{j+1}&\to&\cdots&\to &{\rm G}_n
\end{array}
    \end{split}
\end{equation} where ${\rm G}_0=(k^\ppu, k^\ppd, k^\ppt)$ is the initial prefixed frame and the frames, while ${\rm F}_i$, ${\rm G}_i$ are frames defined via \beqa{FG}&&{\rm F}_j=(\n_j, \ \cdot , {\rm S}^\ppj)\quad{\rm G}_j=({\rm n}_j,\ \cdot,P^\ppj)\qquad j=1,\cdots, n. \eeqa By construction, each frame in the chain has its first axis coinciding with the intersection of  horizontal plane with the horizontal plane of the previous frame (hence, in particular, $\n_j\perp {\rm S}^\ppj$ and $ {\rm n}_j\perp P^\ppj$). Denote as ${\rm T}_j$ the rotation matrix which describes the change of coordinates from ${\rm G}_{j-1}$ to ${\rm F}_j$ and as ${\rm S}_j$ {the} one from ${\rm F}_j$ to ${\rm G}_j$. The matrices $\cT_j$, $\cS_j$ have just the expressions claimed in \equ{good incli*}, \equ{CC}. This follows from the definitions of $(\Theta, \chi,\vartheta, \k)$ in \equ{belle*}. Then we have the following sequence of transformations \beqano
\begin{array}
    {cccccccccccccccccc} &{\rm T}_1&&{\rm S}_1& &&\cdots&&&{\rm S}_j&&
     &\cdots&{\rm S}_n& \\
    \\
    \dst {\rm G}_0&\to&{\rm F}_1&\to& {\rm G}_1&\to&\cdots&\to& {\rm F}_j&\to& {\rm G}_j&\to&
    \cdots&\to &{\rm G}_n
\end{array}
\eeqano connecting ${\rm G}_0$ to any other frame in the chain. From this, and the definitions of the frames \equ{FG}, the formulae for $P_\cP^\ppj$ in \equ{P*ecc new} and $${\rm S}_\cP^\ppj=\chi_{j-1}{\rm T}_1{\rm S}_1\cdots {\rm T}_{j-1}{\rm S}_{j-1}{\rm T}_j k^\ppt$$ follow at once. Hence, also the formulae for ${\rm C}_\cP^\ppj$, which is given by ${\rm C}_\cP^\ppj={\rm S}_\cP^\ppj-{\rm S}_\cP^{(j+1)}$, with ${\rm S}_\cP^{(n+1)}\equiv0$. \qed

\nl For the proof of Theorem \ref{main}, we shall use three auxiliary maps, that we shall denote as $\widetilde\cP$, $\widetilde{\cDel}$ and ${\cDel}$. The map $\widetilde\cP$ is very closely related to $\cP$; $\widetilde{\cDel}$ and ${\cDel}$ are well known: in the literature they are often referred to as {(two variants of)} {\sl Delaunay} maps.

\vskip.1in\noi\paragraph{\bf The map $\widetilde\cP$}

\nl Define the set $$\cC_{\widetilde\cP}:=\Big\{ (y, x)\in\real^{3n}\times \real^{3n}: \quad x^\ppj\ne 0,\quad \widetilde{\rm n}_j:\ne 0,\quad \widetilde\n_j\ne 0\quad \forall\ j=1,\cdots, n \Big\},$$ where, for $(y, x)\in\real^{3n}\times\real^{3n}$, with $y=(y^\ppu,\cdots, y^\ppn)$, $x=(x^\ppu,\cdots, x^\ppn)$, $x^\ppj\ne 0$, {we let} \beqano \widetilde\n_j:=\left\{
\begin{array}
    {llll} \dst k^\ppt\times {\rm S}_\cC^\ppu\quad &j=1\\
    \\
    \dst \frac{x^{(j-1)}}{\|x^{(j-1)}\|}\times {\rm S}^{(j)}_\cC &j=2,\cdots, n
\end{array}
\right.\qquad \widetilde{\rm n}_j:=\dst{\rm S}_\cC^\ppj\times \frac{x^\ppj}{\|x^\ppj\|} \eeqano with $j=1$, $\cdots$, $n$ and \beqa{CScart} &&{\rm C}_\cC^\ppj:=x^\ppj\times y^\ppj,\quad {\rm S}_\cC^\ppj:=\sum_{i=j}^n{\rm C}^\ppi. \eeqa

\nl Define a map $$\widetilde\cP^{-1}:\quad (y,x)\in\cC_{\widetilde\cP}\to (\widetilde\Theta,\widetilde\chi,\widetilde{\rm R},\widetilde\vartheta,\widetilde\k,\widetilde{\rm r})\in \real^n\times \real_+^n\times\real^n\times \torus^n\times \torus^n\times\real_+^n$$ with \beqano
\begin{array}
    {lll} \dst\widetilde\Theta=(\widetilde\Theta_0,\cdots,\widetilde\Theta_{n-1})\qquad &\widetilde\vartheta=(\widetilde\vartheta_0,\cdots,\widetilde\vartheta_{n-1})\\
    \\ \dst\widetilde\chi=(\widetilde\chi_0,\cdots,\widetilde\chi_{n-1})\qquad &\widetilde\k=(\widetilde\k_0,\cdots,\widetilde\k_{n-1})\\
    \\ \widetilde{\rm R}=(\widetilde{\rm R}_1,\cdots, \widetilde{\rm R}_n) &\widetilde{\rm r}=(\widetilde{\rm r}_1,\cdots, \widetilde{\rm r}_n)
\end{array}
\eeqano via the following formulae: \beqano
\begin{array}
    {lllllll} \dst\widetilde{\rm R}_j=\frac{y^\ppj\cdot x^\ppj}{\|x^\ppj\|}\qquad \qquad &\widetilde{\rm r}_j=\|x^\ppj\|\qquad &j=1,\cdots, n \\
    \\
    \dst\widetilde\chi_{j-1}=\|{\rm S}_\cC^\ppj\| & \widetilde\k_{j-1}=\a_{{\rm S}_\cC^\ppj}(\widetilde\n_j, \widetilde{\rm n}_j)&j=1,\cdots, n \\
    \\
    \dst \widetilde\Theta_{j-1}=\left\{
    \begin{array}
        {lrrr} \dst {\rm S}_\cC^\ppu\cdot k^\ppt \\
        \\
        \dst{\rm S}_\cC^{(j)}\cdot\frac{ x^{(j-1)}}{\|x^{(j-1)}\|}
    \end{array}
    \right.& \widetilde\vartheta_{j-1}=\left\{
    \begin{array}
        {lrrr} \dst\a_{k^\ppt}(k^\ppu, \widetilde\n_1)\qquad\\
        \\
        \dst \a_{\frac{x^{(j-1)}}{\|x^{(j-1)}\|}}(\widetilde{\rm n}_{j-1}, \widetilde\n_{j})
    \end{array}
    \right.&
    \begin{array}
        {lll} j=1\\
        \\
        j=2,\cdots, n.
    \end{array}
    
\end{array}
\eeqano
\begin{lemma}
    \label{tilde P injective} Let $\cD_{\widetilde\cP}$ be the set of $(\widetilde\Theta,\widetilde\chi,\widetilde{\rm R},\widetilde\vartheta,\widetilde\k,\widetilde{\rm r})\in \real^n\times \real_+^n\times\real^n\times \torus^n\times \torus^n\times\real_+^n$ such that $(\widetilde\Theta,\widetilde\chi,\widetilde\vartheta,\widetilde\k)$ satisfies \equ{non parallelism}, and let $\widetilde{\rm T}_j$, $\widetilde{\rm S}_j$ and ${\rm C}_{\widetilde\cP}^\ppj$ {be} the functions of $(\widetilde\Theta,\widetilde\chi,\widetilde\vartheta,\widetilde\k)$ defined in \equ{good incli*}-\equ{P*ecc}, with $(\widetilde\Theta,\widetilde\chi,\widetilde\vartheta,\widetilde\k)$ replacing $(\Theta,\chi,\vartheta,\k)$.

\nl The map $\widetilde\cP^{-1}$ is a bijection from $\cC_{\widetilde\cP}$ onto the set $\cD_{\widetilde\cP}$. Its inverse map \beqano \widetilde\cP:\quad (\widetilde\Theta,\widetilde\chi,\widetilde{\rm R},\widetilde\vartheta,\widetilde\k,\widetilde{\rm r})\in \cD_{\widetilde\cP}\to (y_{\widetilde\cP}, x_{\widetilde\cP})\in \real^n\times \real^n \eeqano has the following analytical expression: \beqa{map} \arr{\dst x^\ppj_{\widetilde\cP}:= \widetilde{\rm r}_j \widetilde{\rm T}_1\widetilde{\rm S}_1\cdots \widetilde{\rm T}_{j}\widetilde{\rm S}_{j}k^\ppt \\
    \\
    \dst y^\ppj_{\widetilde\cP}:=\frac{\widetilde{\rm R}_j}{\widetilde{\rm r}_j}x_{\widetilde\cP}^\ppj+\frac{1}{\widetilde{\rm r}_j^2}{\rm C}_{\widetilde\cP}^\ppj\times x_{\widetilde\cP}^\ppj\qquad 1\le j\le n } \eeqa

\nl Moreover, the following relation holds: \beq{S and C} {\rm C}_{\widetilde\cP}^\ppj={\rm C}_\cC^\ppj\circ \widetilde\cP=x^\ppj_\cP\times y^\ppj_\cP. \eeq
\end{lemma}
\begin{proof} With similar arguments as the ones of the proof of Proposition \ref{inversion}, but replacing, in the diagram \equ{chain}, ${\rm S}^\ppj_\cE$ with ${\rm S}^\ppj_{\cC}$, $P_\cP^\ppj$ with $\frac{x^\ppj}{\|x^\ppj\|}$ and the nodes $\n_k$, ${\rm n}_k$ with $\widetilde\n_k$, $\widetilde{\rm n}_k$, one finds the formula for $x^\ppj_{\widetilde\cP}$ in \equ{map}, the formula for $${\rm S}_{\widetilde\cP}^\ppj:={\rm S}^\ppj_\cC\circ\widetilde\cP=\widetilde\chi_{j-1}\widetilde{\rm T}_1\widetilde{\rm S}_1\cdots \widetilde{\rm T}_{j-1}\widetilde{\rm S}_{j-1}\widetilde{\rm T}_j k^\ppt$$ and hence the formula for $${\rm C}^\ppj_\cC\circ\widetilde\cP={\rm S}_{\widetilde\cP}^\ppj-{\rm S}_{\widetilde\cP}^{(j+1)}={\rm C}_{\widetilde\cP}^\ppj$$ being just the formula for ${\rm C}^\ppj_\cP$ in \equ{P*ecc}, with $(\Theta,\chi,\vartheta, \k)$ replaced by $(\widetilde\Theta,\widetilde\chi,\widetilde\vartheta, \widetilde\k)$. With the same argument as in Remark \ref{C perp to P}(ii), we see that $x^\ppj_{\widetilde\cP}\perp{\rm C}^\ppj_{\widetilde\cP}$. Finally, the formula for $y^\ppj_{\widetilde\cP}$ is found taking for $y^\ppj_{\widetilde\cP}$ the unique vector verifying $$y^\ppj_\cP\cdot\frac{x^\ppj_\cP}{\|x^\ppj_\cP\|}={\rm R}_j\qquad x^\ppj_\cP\times y^\ppj_\cP={\rm C}^\ppj_\cP.$$
\end{proof}
\begin{lemma}
    \label{lem: good change} $\widetilde\cP$ preserves the standard Liouville 1-form: \beq{1form}\sum_{j=1}^ny_{\widetilde\cP}^\ppj\cdot dx_{\widetilde\cP}^\ppj=\sum_{j=1}^{n}\big(\widetilde\Theta_{j-1}d\widetilde\vartheta_{j-1}+\widetilde\chi_{j-1}d\widetilde\k_{j-1}+\widetilde{\rm R}_j d\widetilde{\rm r}_j\big).\eeq
\end{lemma}
\nl The proof of Lemma \ref{lem: good change} uses the {following easy lemma:}
\begin{lemma}
    [\cite{chierchiaPi11a}]\label{lemma on reduction} Let $$x=\cR_3(\theta)\cR_1(i)\bar x,\quad y=\cR_3(\theta)\cR_1(i)\bar y,\quad {\rm C}:=x\times y,\quad \bar{\rm C}:=\bar x\times \bar y,$$ with $x,\bar x, y, \bar y\in\real^3$.  Then, $$y\cdot dx={\rm C}\cdot k^{(3)}d\theta+\bar{\rm C}\cdot k^{(1)} di+\bar y \cdot d\bar x.$$
\end{lemma}
\vskip.in\textsc{Proof of Lemma \ref{lem: good change}.} We may write $$x_{\widetilde\cP}^\ppj=\widetilde{\rm T}_1\widetilde{\rm S}_1\cdots\widetilde {\rm T}_j\widetilde{\rm S}_j\tilde x^\ppj,\quad y_{\widetilde\cP}^\ppj=\widetilde{\rm T}_1\widetilde{\rm S}_1\cdots\widetilde {\rm T}_j\widetilde{\rm S}_j\tilde y^\ppj,\quad {\rm C}_{\widetilde\cP}^\ppj=\widetilde{\rm T}_1\widetilde{\rm S}_1\cdots\widetilde {\rm T}_j\widetilde{\rm S}_j\tilde {\rm C}^\ppj$$ where \beqa{Rr} \tilde x^\ppj&:=&\widetilde{\rm r}_j k^\ppt\qquad j=1,\cdots, n-1\nonumber\\
\tilde y^\ppj&:=&\widetilde{\rm R}_jk^\ppt+\frac{1}{\widetilde{\rm r}_j}\tilde {\rm C}^\ppj\times k^\ppt\nonumber\\
\tilde {\rm C}^\ppj&:=&\widetilde\chi_{j-1}\widetilde{\cS}_j^{-1}k^\ppt-\widetilde\chi_j{\widetilde\cT}_{j+1}k^\ppt=\tilde x^\ppj\times \tilde y^\ppj \eeqa with $\widetilde\chi_n:=0$, $\widetilde\cS_n:=\id$. We also let, for $1\le k\le j\le n$ and $1\le i\le n-1$, \beqano &&\hat{\rm C}^\ppj_k=\widetilde\cS_k(\widetilde\cT_{k+1}\widetilde\cS_{k+1}\cdots\widetilde\cT_j\widetilde\cS_j)\tilde{\rm C}^\ppj,\quad \check{\rm C}^\ppj_k=\widetilde\cT_{k}\widetilde\cS_{k}\cdots\widetilde\cT_j\widetilde\cS_j\tilde{\rm C}^\ppj,\quad \check{\rm C}^\ppj_{j+1}:=\tilde{\rm C}^\ppj\nonumber\\
&&\hat{\rm S}^\ppj_k:=\sum_{m=j}^n\hat{\rm C}^\ppm_k,\quad \check{\rm S}^\ppj_k:=\sum_{m=j}^n\check{\rm C}^\ppm_k,\quad \check{\rm S}^\ppi_{i+1}:=\tilde{\rm C}^\ppi+\check{\rm S}^{(i+1)}_{i+1} \eeqano where the product $\widetilde\cT_{k+1}\widetilde\cS_{k+1}\cdots\widetilde\cT_j\widetilde\cS_j$ is to be replaced with the identity when $k=j$. We have the following identities (implied by ${\rm S}^\ppj=\sum_{k=j}^n{\rm C}^\ppk$):

\begin{equation}\label{SC} \begin{split}\check{\rm S}^\ppj_j=\sum_{k=j}^n\check{\rm C}^\ppk_j=\widetilde\chi_{j-1}\widetilde\cT_jk^\ppt,\quad \hat{\rm S}^\ppj_j=\sum_{k=j}^n\hat{\rm C}^\ppk_j=\widetilde\chi_{j-1}k^\ppt,\quad \check{\rm S}^\ppi_{i+1}=\widetilde\chi_{j-1}{\widetilde\cS}_i^{-1}k^\ppt. 
\end{split}\end{equation} Applying Lemma \ref{lemma on reduction} repeatedly and using (as it follows from \equ{Rr}) $$\tilde y^\ppj\cdot d\tilde x^\ppj=\widetilde{\rm R}_j d\widetilde{\rm r}_j,$$ we have, for $1\le j\le n$, \[ y_{\widetilde\cP}^\ppj\cdot x_{\widetilde\cP}^\ppj ~=~ \sum_{k=1}^j\big( \check{\rm C}^\ppj_{k}\cdot k^\ppt d\widetilde\vartheta_{k-1}+\hat{\rm C}^\ppj_k\cdot k^\ppu d\widetilde\iota_k+\hat{\rm C}^\ppj_k\cdot k^\ppt d\widetilde\k_{k-1}+\check{\rm C}^\ppj_{k+1}\cdot k^\ppu d\widetilde{\rm i}_k \big)+\widetilde{\rm R}_jd\widetilde{\rm r}_j \] where, as in the proof of Lemma \ref{tilde P injective}, $\widetilde\iota_j$, $\widetilde{\rm i}_j$ denote the functions $\iota_j$, ${\rm i}_j$ in \equ{good incli*}, with $\Theta_i$, $\chi_i$ replaced by $\widetilde{\Theta}_i$, $\widetilde\chi_i$. Note that we have used $d\,\widetilde{\rm i}_n\equiv 0$, since, by definition, $\widetilde{\rm i}_n=\frac{\p}{2}$. Taking the sum over $j=1$, $\cdots$, $n$, \[ \sum_{j=1}^ny_{\widetilde\cP}^\ppj \!\cdot  dx_{\widetilde\cP}^\ppj ~=~ \sum_{j=1}^n\check{\rm S}_{j}^{(j)} \!\cdot  k^\ppt d\widetilde\vartheta_{j-1}+\hat{\rm S}_{j}^{(j)} \!\cdot  k^\ppu d\widetilde\iota_{j}+\hat{\rm S}_{j}^\ppj \!\cdot  k^\ppt d\widetilde\k_{j-1}+\check{\rm S}^{(j)}_{j+1} \!\cdot  k^\ppu d\widetilde{\rm i}_{j} + \sum_{j=1}^n\widetilde{\rm R}_j d\widetilde{\rm r}_j. \] In view of \equ{SC} and of the definitions in \equ{good incli*}-\equ{CC}, we then find \equ{1form}. \qed

 \vskip.1in\noi\paragraph{\bf The map $\widetilde{\cDel}$}

The map \beqano \widetilde{\cDel}:\quad (\widetilde{\rm H}, \widetilde\G,\widetilde{\rm R} , \widetilde{\rm h}, \widetilde{\rm g}, \widetilde{\rm r})\in \cD_{\widetilde{\cDel}}\to (y_{\widetilde{\cDel}}, x_{\widetilde{\cDel}})\in \real^{3n}\times\real^{3n}\eeqano is defined on the set \beqano \cD_{\widetilde{\cDel}}&:=&\Big\{ (\widetilde{\rm H}, \widetilde\G,\widetilde{\rm R} , \widetilde{\rm h}, \widetilde{\rm g}, \widetilde{\rm r})=(\widetilde{\rm H}_1,\cdots, \widetilde{\rm H}_n, \widetilde\G_1,\cdots, \widetilde\G_n,\widetilde{\rm R}_1,\cdots, \widetilde{\rm R}_n , \widetilde{\rm h}_1, \cdots,\widetilde{\rm h}_n,\nonumber\\
&& \widetilde{\rm g}_1,\cdots, \widetilde{\rm g}_n, \widetilde{\rm r}_1,\cdots, \widetilde{\rm r}_n)\in \real^{3n}\times \torus^{2n}\times \real_+^n: \quad \widetilde{\rm r}_j>0,\quad \widetilde\G_j>0,\quad \frac{|\widetilde{\rm H}_j|}{\widetilde\G_j}<1\nonumber\\
&& \forall\ j=1,\cdots, n \Big\}\eeqano via the following formulae: \beqano &&x^\ppj_{\widetilde{\cDel}}:={\rm R}_3(\widetilde{\rm h_j})\cR_1(\widetilde i_j)\ovl x^\ppj_{\widetilde{\cDel}},\qquad y^\ppj_{\widetilde{\cDel}}:={\rm R}_3(\widetilde{\rm h_j})\cR_1(\widetilde i_j)\ovl y^\ppj_{\widetilde{\cDel}} \eeqano where \beqano \widetilde i_j&:=&\cos^{-1}\frac{\widetilde{\rm H}_j}{\widetilde\G_j}\in (0,\p)\nonumber\\
\ovl x^\ppj_{\widetilde{\cDel}}&:=& \widetilde{\rm r}_j\cos\widetilde{\rm g}_j k^\ppu+\widetilde{\rm r}_j\sin\widetilde{\rm g}_jk^\ppd\nonumber\\
\ovl y^\ppj_{\widetilde{\cDel}}&:=&\big( \widetilde{\rm R}_j\cos\widetilde{\rm g}_j-\frac{\widetilde{\G}_j}{\widetilde{\rm r}_j}\sin\widetilde{\rm g}_j\big)k^\ppu+\big(\widetilde{\rm R}_j\sin\widetilde{\rm g}_j+\frac{\widetilde\G_j}{\widetilde{\rm r}_j}\cos\widetilde{\rm g}_j\big) k^\ppd.\eeqano
\begin{lemma}
    [Delaunay]\label{Del} $\widetilde{\cDel}$ is a bijection from the domain $\cD_{\widetilde{\cDel}}$ onto the set \beqano {\rm C}_{\widetilde{\cDel}}&:=&\Big\{(y, x)=(y^\ppu, \cdots, y^\ppn, x^\ppu, \cdots, x^\ppn)\in \real^{3n}\times \real^{3n}:\nonumber\\
    &&\widetilde{\eufm n}_j:=k^\ppt\times {\rm C}_\cC^\ppj\ne 0,\quad x^\ppj\ne 0\qquad\forall\ j=1,\cdots, n\Big\} \eeqano where ${\rm C}_\cC^\ppj$ is as in \equ{CScart}. The formulae for the inverse map $$\widetilde{\cDel}^{-1}:\quad (y, x)\in\cC_{\widetilde{\cDel}}\to (\widetilde{\rm H}, \widetilde\G, \widetilde{\rm R}, \widetilde{\rm h}, \widetilde{\rm g}, \widetilde{\rm r})\in\cD_{\widetilde{\cDel}}$$ are \begin{equation}\label{auxiliary Del}\left\{
    \begin{array}
        {l} \dst\widetilde{\rm H}_j={\rm C}_\cC^\ppj\cdot k^\ppt\\
        \\
        \dst\widetilde{\rm h}_j:=\a_{k^\ppt}(k^\ppu,\widetilde{\eufm n}_j)
    \end{array}
    \right. \qquad\left\{
    \begin{array}
        {l} \dst\widetilde\G_j:=\|{\rm C}_\cC^\ppj\|\\
        \\
        \dst\widetilde{\rm g}_j:=\a_{{\rm C}_\cC^\ppj}(\widetilde{\eufm n}_j,x^\ppj)
    \end{array}
    \right. \qquad \arr{\dst\widetilde{\rm R}_j=\frac{y^\ppj\cdot x^\ppj}{\|x^\ppj\|}\\
    \\
    \dst \widetilde{\rm r}_j=\|x^\ppj\| } \end{equation} Finally, $\widetilde{\cDel}$ preserves the standard Liouville 1-form \[ \sum_{i=1}^ny_{\widetilde{\cDel}}^\ppi\cdot dx_{\widetilde{\cDel}}^\ppi=\sum_{i=1}^n\big(\widetilde{\rm H}_id\widetilde{\rm h}_i+\widetilde\G_id\widetilde{\rm g}_i+\widetilde{\rm R}_i d\widetilde{\rm r}_i\big). \]
\end{lemma}
We omit the proof of Lemma \ref{Del}, which may be found in classical textbooks.

 \vskip.1in\noi\paragraph{\bf The map ${\cDel}$}
\begin{definition}
    [Delaunay map]\label{Def: Delaunay}\rm Let \beqano {\eufm X}_{{\cDel}}&:=&\Big\{ {\rm X}_{\cDel}:=({\rm H}, \G,\L , {\rm h}, {\rm g})=({\rm H}_1,\cdots, {\rm H}_n, \G_1,\cdots, \G_n,\L_1,\cdots, \L_n , {\rm h}_1, \cdots,{\rm h}_n,\nonumber\\
    && {\rm g}_1,\cdots, {\rm g}_n)\in \real^{3n}\times \torus^{2n}: \quad \quad \G_j>0,\quad \frac{|{\rm H}_j|}{\G_j}<1,\quad \L_j>0\nonumber\\
    && \forall\ j=1,\cdots, n \Big\}\eeqano and let $\cE_{\cDel}$ be the set of n-plets $({\eufm E}_1,\cdots, {\eufm E}_n)$ where ${\eufm E}_j=(a_j, e_j, N^\ppj, P^\ppj)$ satisfies $$0<e_j<1,\qquad {\eufm n}_j:=k^\ppt\times N^\ppj\ne 0,\quad \forall\ j=1,\cdots, n.$$

\nl Fix positive numbers ${\eufm M}_1$, $\cdots$, ${\eufm M}_n$, ${\eufm m}_1$, $\cdots$, ${\eufm m}_n$. {Define} $$\t_{\cDel}:\ {\rm X}_{\cDel}:=({\rm H},\G,\L, {\rm h}, {\rm g})\in {\eufm X}_{\cDel}\to {\eufm E}_{\cDel}=({\eufm E}_{1,\cDel},\cdots, {\eufm E}_{n,\cDel}) $$ with ${\eufm E}_{j,\cDel}=(a_{j,\cDel}, e_{j,\cDel}, N_{\cDel}^\ppj, P_{\cDel}^\ppj)$ and \beqano
    \begin{array}
        {llll} &\dst a_{j,\cDel}=\frac{1}{{\eufm M}_j}(\frac{\L_j}{{\eufm m}_j})^2,\qquad &e_{j,\cDel}=\sqrt{1-(\frac{\G_j}{\L_j})^2}\\
        \\
        &\dst N_{\cDel}^\ppj={\rm R}_3({\rm h_j})\cR_1(i_j)k^\ppt&P_{\cDel}^\ppj={\rm R}_3({\rm h_j})\cR_1(i_j) \cR_3({\rm g}_j)k^\ppu
    \end{array}
    \eeqano where $i_j:=\cos^{-1}\frac{{\rm H}_j}{\G_j}$. 

\nl We call {\sl Delaunay map} the map \beqa{Del map} {\cDel}: \quad {\rm Del}=({\rm H},\G,\L, {\rm h}, {\rm g}, \ell)\in \cD_{{\cDel}}\to (y_{{\cDel}}, x_{{\cDel}})\in \real^{3n}\times\real^{3n}\eeqa which is defined on the domain \beqano \cD_{{\cDel}}&:=&{\eufm X}_{{\cDel}}\times \torus^n \eeqano as the Kepler map associated to $\t_{\cDel}$ via the following lemma (the proof of which may be found in classical textbooks).
\end{definition}

\begin{lemma}
    [Delaunay]\label{Delaunay symplectic} $\t_{\cDel}$ is a bijection of $ {\eufm X}_{\cDel}$ onto $\cE_{\cDel}$. Its inverse map $$\t_{\cDel}^{-1}:\quad {\eufm E}_{\cDel}=({\eufm E}_{1,\cDel},\cdots, {\eufm E}_{n,\cDel})\in\cE_{\cDel} \to {\rm X}_{\cDel}\in {\eufm X}_{\cDel} $$ is defined by equations \begin{equation}\label{Delauninv}\begin{split}
     \left\{
    \begin{array}
        {l} \dst{\rm H}_j={\rm C}_\cE^\ppj\cdot k^\ppt\\
        \\
        \dst{\rm h}_j:=\a_{k^\ppt}(k^\ppu,{\eufm n}_j)
    \end{array}
    \right. \quad\left\{
    \begin{array}
        {l} \dst\G_j=\|{\rm C}_\cE^\ppj\|\\
        \\
        \dst{\rm g}_j:=\a_{{\rm C}_\cE^\ppj}({\eufm n}_j,P^\ppj)
    \end{array}
    \right. \quad \L_j={\eufm m}_j\sqrt{{\eufm M}_j a_j},
    \end{split}\end{equation} where ${\rm C}_\cE^\ppj$ is as in \equ{former and latter}. Furthermore, ${\rm D}e\ell$ preserves the standard 2-form $$\sum_{j=1}^ndy_{{\cDel}}^\ppj\wedge dx_{{\cDel}}^\ppj =\sum_{j=1}^n\big(d{\rm H}_j\wedge d{\rm h}_j+d\G_j\wedge d{\rm g}_j+d\L_j\wedge d\ell_j\big).$$
\end{lemma}
\nl Now we are ready to complete the proof of Theorem~\ref{main}. \vskip.1in\noi\textsc{Proof of  Theorem \ref{main}.} Let $$ \cD_\cP^*:=\Big\{ {\rm P}=(\Theta, \chi, \L,\vartheta, \k,\ell)\in \cD_\cP: {\rm P}({\rm P})\in \cC_{\cDel} \Big\}. $$ It is enough to prove Theorem \ref{main} on $\cD_\cP^*$, since indeed the $\cP$-map is regular on $\cD_\cP=\ovl{\cD_\cP^*}$. On $\cD_\cP^*$, we consider the map 
\begin{equation} 
\begin{split} &{\phi_{\cDel}^\cP}:= \cDel^{-1}\circ\cP: \nonumber\\
& {\rm P}=(\Theta,\chi, \L,\vartheta, \k,\ell)\in \cD_\cP^*\to{\rm Del}=({\rm H}, \G, \L,{\rm h}, {\rm g}, \ell)\in \cD^*_{\cDel}:={\phi_{\cDel}^\cP}(\cD_\cP^*)\subset \cD_{\cDel}.
\end{split} \end{equation} 
$\phi_{\cDel}^\cP$ gives the Delaunay coordinates at left hand side in \equ{Del map} in terms of the P-coordinates at left hand side of \equ{cal P map} in the subset $\cD_\cP^*$ of $\cD_\cP$ the $\cP$-image of which lies in the $\cDel$-image of $\cD_{\cDel}$. Clearly, ${\phi_{\cDel}^\cP}$ leaves the $(\L,\ell)$ unvaried. More precisely, ${\phi_{\cDel}^\cP}$ decouples into two disjoint maps: the identity on the $(\L,\ell)$, and a $4n$-dimensional map $$\widehat{\phi_{\cDel}^\cP}:\quad (\Theta,\chi,\vartheta, \k)\in \widehat{\cD_\cP^*}\to ({\rm H}, \G, {\rm h}, {\rm g})\in \widehat{\cD^*_{\cDel}}={\phi_{\cDel}^\cP}(\widehat{\cD_\cP^*})\subset \widehat{\cD_{\cDel}} $$ on the remaining coordinates, which turns out to be a bijection of the sets $\widehat{\cD_\cP^*}$ and $\widehat{\cD^*_{\cDel}}$. Here, the map $\widehat{\phi_{\cDel}^\cP}$ and the sets $\widehat{\cD_\cP^*}$ and $ \widehat{\cD_{\cDel}}$ do not depend on $(\L,\ell)$. Indeed, the explicit expressions of $\widehat{\phi_{\cDel}^\cP}$, $\widehat{\cD_\cP^*}$ in terms of $ {\rm P}=(\Theta,\chi, \L,\vartheta, \k,\ell)$; or {of} $ \widehat{\cD_{\cDel}}$ in terms of ${\rm Del}=({\rm H}, \G, \L,{\rm h}, {\rm g}, \ell)$ involve only the ${\rm C}^\ppj_\cP$, $P^\ppj_\cP$; the ${\rm C}^\ppj_{\cDel}$, $P^\ppj_{\cDel}$, that do not depend on $(\L,\ell)$: \equ{Delauninv} (where one has to replace $\cC$ with $\cP$), \equ{CC} and \equ{P*ecc new}.

\nl In view of the previous consideration and of Lemma \ref{Delaunay symplectic}, Theorem \ref{main} is implied by
\begin{lemma}
    \label{homogeneous with respect to Del} The map $\widehat{\phi_{\cDel}^\cP}$ preserves that standard 1-form: $$\sum_{j=1}^n\big({\rm H}_jd{\rm h}_j+\G_jd{\rm g}_j\big)=\sum_{j=1}^n\big(\Theta_{j-1}d\vartheta_{j-1}+\chi_{j-1}d\k_{j-1}\big).$$
\end{lemma}
\begin{proof} We look at the {analogous} map $$\widehat{\phi_{\widetilde{\cDel}}^{\widetilde{\cP}}}:\quad (\widetilde\Theta,\widetilde\chi,\widetilde\vartheta,\widetilde \k)\in \widehat{\cD_{\widetilde\cP}^*}\to (\widetilde{\rm H}, \widetilde\G, \widetilde{\rm h}, \widetilde{\rm g})\in \widehat{\cD^*_{\widetilde{\cDel}}}={\phi_{\widetilde{\cDel}}^{\widetilde\cP}}(\widehat{\cD_{\widetilde\cP}^*})\subset \widehat{\cD_{\widetilde{\cDel}}}. $$ The analytical expression of this map is identical to the one of $\widehat{\phi_{\cDel}^\cP}$. This follows from the fact that $\widehat{\phi_{\widetilde{\cDel}}^{\widetilde{\cP}}}$ depends on the coordinates $(\widetilde\Theta,\widetilde\chi,\widetilde\vartheta,\widetilde \k)$ only via ${\rm C}^\ppj_{\widetilde\cP}$ and $\frac{x^\ppj_{\widetilde\cP}}{\|x^\ppj_{\widetilde\cP}\|}$ exactly as $\widehat{\phi_{\cDel}^\cP}$ depends on $(\Theta,\chi,\vartheta, \k)$ only via ${\rm C}^\ppj_{\cP}$ and $\P^\ppj_{\cP}$, that ${\rm C}^\ppj_{\widetilde\cP}$ and $\frac{x^\ppj_{\widetilde\cP}}{\|x^\ppj_{\widetilde\cP}\|}$ have exactly the same expressions of ${\rm C}^\ppj_{\cP}$ and $P^\ppj$, apart for replacing $(\Theta,\chi,\vartheta, \k)$ with $(\widetilde\Theta,\widetilde\chi,\widetilde\vartheta,\widetilde \k)$.  Compare \equ{auxiliary Del} (where one has to replace ${\rm C}^\ppj_\cC$ with ${\rm C}^\ppj_{\widetilde\cP}$), \equ{Delauninv} (where one has to replace ${\rm C}^\ppj_\cE$ with ${\rm C}^\ppj_\cP$), \equ{CC}, \equ{P*ecc new}, \equ{map} and \equ{S and C}. But Lemmata \ref{lem: good change} and \ref{Del} imply that $\widehat{\phi_{\widetilde{\cDel}}^{\widetilde{\cP}}}$ preserves that standard 1-form: $$\sum_{j=1}^n\big(\widetilde{\rm H}_jd\widetilde{\rm h}_j+\widetilde\G_jd\widetilde{\rm g}_j\big)=\sum_{j=1}^n\big(\widetilde\Theta_{j-1}d\widetilde\vartheta_{j-1}+\widetilde\chi_{j-1}d\widetilde\k_{j-1}\big).$$ Then $\widehat{\phi_{\cDel}^\cP}$ does.
 \end{proof}

 \section[The $\cP$-map vs rotations and reflections]{The $\cP$-map vs rotations and reflections}\label{P-map vs rotations and reflections} Now we discuss how the $\cP$-map behaves in presence of symmetries in the Hamiltonian due to rotations or reflections. 

\nl Let ${\rm H}={\rm H}(y, x)$ be the Hamiltonian governing the motion of $n$ particles, where such particles are expressed in the canonical coordinates $(y^\ppu, x^\ppu)$, $\cdots$, $(y^\ppn, x^\ppn)$. Assume that ${\rm H}$ is left unvaried by rotations and reflections. Namely, if $$\phi_{\cR, \cS}:\quad (y^\ppj, x^\ppj)\to ({\rm R}y^\ppj, \cS x^\ppj),\qquad j=1,\cdots, n$$ where $\cR$, $\cS$ are real a $3\times 3$ matrices, then rotation invariance is $${\rm H}\circ\phi_{\cR, \cR}={\rm H}\qquad \forall\ \cR:\ \cR \cR^{\rm t}=\id$$ while reflection invariance is \begin{equation}\begin{split}
&{\rm H}\circ\phi_{\cS_\s,\cS_\t}={\rm H}\quad{\rm for \ some} \quad \cS_\s=\left(
\begin{array}
    {lll} \s_1&0&0\\
    0&\s_2&0\\
    0&0&\s_3
\end{array}
\right)\quad \cS_\t=\left(
\begin{array}
    {lll} \t_1&0&0\\
    0&\t_2&0\\
    0&0&\t_3
\end{array}
\right)\nonumber\\
& \s_i,\ \t_i=\pm 1. 
\end{split}
\end{equation} Rotation invariance is associated to the conservation, through the motion, of the total angular momentum ${\rm S}^\ppu_\cC$ is \equ{CScart}. Reflection invariance is not associated to integrals.

\nl The Hamiltonian ${\rm H}_{\rm hel}$ in \equ{Helio} is rotation and reflection invariant, and reflection invariance holds with any choice of $\s$, $\t$.

\nl Let $${\rm H}_\cP:={\rm H}\circ\cP.$$ The fact that ${\rm S}^\ppu_\cC$ is preserved along the motions of ${\rm H}$ implies that the coordinates $$\Theta_0=Z,\quad \vartheta_0=\zeta,\quad \k_0={\eufm g}$$ do not appear in ${\rm H}_\cP$. Indeed, $Z$ and $ \zeta$ are integrals, while ${\eufm g}$ is conjugated to ${\rm G}=\|{\rm S}_\cP^\ppu\|$, which is an integral for ${\rm H}_\cP$. Thus, the number of degrees of freedom is naturally reduced by two units, once one regards ${\rm G}$ as a prefixed external parameter. Namely, for any fixed $\chi_0={\rm G}$, ${\rm H}_\cK$ may be regarded as a function of the $2(3n-1)$ dimensional coordinates \[\ovl{\rm P}:=(\ovl\Theta,\chi,\L,\ovl\vartheta,\k,\ell)\] {\sl which does not depend on $\k_0$}. Here, \beqano \bar\Theta=(\Theta_1,\cdots,\Theta_{n-1}),\quad \bar\vartheta=(\vartheta_1,\cdots,\vartheta_{n-1}). \eeqano An analogue property is also shared with the action-angle coordinates $(\Psi$, $\G$, $\L$, $\psi$, $\g$, $\ell)$ described in \cite{pinzari-th09, chierchiaPi11a}, and related to a set of coordinates discovered by A. Deprit \cite{deprit83} in the 80s (compare also \cite{zhao14} or the Appendix \ref{comparison}).

\nl The main novelty introduced by the $\cP$-coordinates (that does not hold for the coordinates of \cite{chierchiaPi11a}) is how $\cP$ behaves {relative to} reflections.

\nl We denote as $$\cR_2^-:=\phi_{\cS_{\s^\ppd},\cS_{\s^\ppd}}\qquad \s^\ppd=(1,-1,1)$$ the reflection of the second coordinate both for the $y^\ppj$'s and the $x^\ppj$'s and we let \[\cS^-(\Theta,\chi,\L,\vartheta,\k,\ell):=(-\Theta,\chi,\L, -\vartheta, \k,\ell).\]
\begin{proposition}
    \label{integrability} \beq{commutation}{\rm R}^-_2\circ\cP=\cP\circ \cS^-.\eeq Therefore, if ${\rm H}={\rm H}(y, x)$ satisfies $${\rm H}\circ\cR_2^-={\rm H}$$ then ${\rm H}_\cP:={\rm H}\circ\cP$ satisfies $${\rm H}_\cP\circ \cS^-={\rm H}_\cP.$$ Hence, any of the points $$\Theta_0=\cdots=\Theta_{n-1}=0,\quad (\vartheta_0, \cdots,\vartheta_{n-1})=(k_0,\cdots, k_{n-1})\p \quad \mod 2\p\integer^{n}$$ is an equilibrium point for ${\rm H}_\cP$, for any $(\chi,\L,\k,\ell)$.
\end{proposition}
\begin{proof} Defining $\cR^\ppj:={\rm T}_j{\rm S}_j$, $s^\ppj:={\rm T}_jk^\ppt$, we write the vectors $P^\ppj_\cP$ and ${\rm S}^\ppj_\cP$ (compare Eq.~\equ{P*ecc new} and Remark \ref{C perp to P}(iii)) as $$P_\cP^\ppj=\cR^\ppu\cdots\cR^\ppj k^\ppt,\quad {\rm S}_\cP^\ppj=\chi_{j-1}\cR^\ppu\cdots\cR^\ppj s^\ppj.$$ The explicit expressions of $\cR^\ppj$ and $s^\ppj$ are

 \beqano \cR^\ppj_{11}&=&\cos\k_{j-1}\cos\vartheta_{j-1}-\sin\k_{j-1}\cos\iota_j\sin\vartheta_{j-1}\nonumber\\
\cR^\ppj_{21}&=&\cos\k_{j-1}\sin\vartheta_{j-1}+\sin\k_{j-1}\cos\iota_j\cos\vartheta_{j-1}\nonumber\\
\cR^\ppj_{31}&=&\sin\k_{j-1}\sin\iota_j\nonumber\\
\cR^\ppj_{12}&=&-\cos{\rm i}_j\sin\k_{j-1}\cos\vartheta_{j-1}+\sin\vartheta_{j-1}(-\cos{\rm i}_j\cos\iota_j\cos\k_{j-1}+\sin\iota_j\sin{\rm i_j})\nonumber\\
\cR^\ppj_{22}&=&-\cos{\rm i}_j\sin\k_{j-1}\sin\vartheta_{j-1}-\cos\vartheta_{j-1}(-\cos{\rm i}_j\cos\iota_j\cos\k_{j-1}+\sin\iota_j\sin{\rm i_j})\nonumber\\
\cR^\ppj_{32}&=&\cos{\rm i}_j\cos\k_{j-1}\sin\iota_j+\sin{\rm i}_j\cos{\iota}_j\nonumber\\
\cR^\ppj_{13}&=&\sin{\rm i}_j\sin\k_{j-1}\cos\vartheta_{j-1}+\sin\vartheta_{j-1}(\sin{\rm i}_j\cos\iota_j\cos\k_{j-1}+\sin\iota_j\cos{\rm i_j})\nonumber\\
\cR^\ppj_{23}&=&\sin{\rm i}_j\sin\k_{j-1}\sin\vartheta_{j-1}-\cos\vartheta_{j-1}(\sin{\rm i}_j\cos\iota_j\cos\k_{j-1}+\sin\iota_j\cos{\rm i_j})\nonumber\\
\cR^\ppj_{33}&=&-\sin{\rm i}_j\cos\k_{j-1}\sin\iota_j+\cos{\rm i}_j\cos{\iota}_j\nonumber\\
s^\ppj_1&=&\sin\iota_j\sin\vartheta_{j-1}\nonumber\\
s^\ppj_2&=&-\sin\iota_j\cos\vartheta_{j-1}\nonumber\\
s^\ppj_3&=&\cos\iota_j. \eeqano Then $\cS^-$ lets $P_{\cP}^\ppj$ and ${\rm S}_{\cP}^\ppj$ respectively, into $$(P_{\cP}^\ppj)^-:=\cR_2^-P_{\cP}^\ppj \quad{\rm and}\quad({\rm S}_{\cP}^\ppj)^-:=-\cR_2^-{\rm S}_{\cP}^\ppj.$$

\nl Therefore, ${\rm C}^\ppj_{\cP}={\rm S}^\ppj_{\cP}-{\rm S}^{(j+1)}_{\cP}$ (with ${\rm S}^{(n+1)}_{\cP}:=0$) and $Q^\ppj_{\cP}=\frac{{\rm C}_{\cP}^\ppj}{\|{\rm C}_{\cP}^\ppj\|}\times P^\ppj_{\cP}$ are transformed, respectively, into $$({\rm C}^\ppj_{\cP})^-:=-\cR_2^-{\rm C}_{\cP}^\ppj,\qquad (Q^\ppj_{\cP})^-:=\cR_2^-Q_{\cP}^\ppj.$$ On the other hand, $a_{j,\cP}$ and $e_{j,\cP}$ are left unvaried by $\cS^-$. In view of Definition \ref{Kepler map} and Definition \ref{def: P map}, the thesis \equ{commutation} follows. \end{proof}

\chapter{The $\cP$-map and the planetary problem}

\nl After the reduction of the invariance by translations, a Hamiltonian governing the motions of $n$ planets with masses $\m m_1$, $\cdots$, $\m m_n$ interacting among themselves and with a star with mass $m_0$ can be taken to be the ``heliocentric'' one 
\beq{Helio}{\rm H}_{\rm hel}:=\sum_{1\le i\le n}\bigg(\frac{\|y^\ppi\|^2}{2 {\eufm m}_i}-\frac{ {\eufm m}_i{\eufm M}_i}{\|x^\ppi\|}\bigg)+\m \sum_{1\le i<j\le n}\bigg(\frac{y^\ppi\cdot y^\ppj}{m_0}-\frac{m_im_j}{\|x^\ppi-x^\ppj\|}\bigg)\eeq where $(y,x) = (y^\ppu, \cdots, y^\ppn, x^\ppu, \cdots, x^\ppn)$ are ``Cartesian coordinates'' taking values on the ``collision-less'' phase space $\real^{3n}\times \real^{3n}\setminus\D$, where $$\D=\Big\{x=(x^\ppu, \cdots, x^\ppn)\in \real^3\times \cdots\times\real^3:\quad 0\ne x^\ppi\ne x^\ppj\quad \forall\ 1\le i<j\le n\Big\}$$ endowed with the standard 2- form $$\O:=dy\wedge dx:=\sum_{i=1}^n \sum_{j=1}^3dy_j^\ppi\wedge dx^\ppi_j$$ and with \beq{reduced masses}{\eufm M}_i=m_0+\m m_i\qquad {\eufm m}_i=\frac{m_0m_i}{m_0+\m m_i}\eeq being the so-called ``reduced masses''.

\nl In the following Section \ref{A remark on a general property of Kepler maps} we describe a general property of Kepler maps, in relation to their application to the Hamiltonian ${\rm H}_{\rm hel}$. Then (in Section \ref{The case of the P-map}) we shall specialize to the case of the $\cP$-map.

 \section{A general property of Kepler maps}\label{A remark on a general property of Kepler maps}

 For a general Kepler map $\cK$, we denote $${\rm H}_\cK({{\rm K}}):={\rm H}_{\rm hel}\circ\cK=-\sum_{j=1}^n\frac{{\eufm m}_j{\eufm M}_j}{2 a_{j,\cK}({\rm X}_\cK)}+\m f_{\cK}({{\rm K}}),$$ where  $$f_{\cK}({{\rm K}}):=\sum_{1\le i<j\le n}\bigg(\frac{y_\cK^\ppi\cdot y_\cK^\ppj}{m_0}-\frac{m_im_j}{\|x_\cK^\ppi-x_\cK^\ppj\|}\bigg)$$ and $y^\ppj_\cK$, $x^\ppj_\cK$ are as in Definition \ref{Kepler map}. 

\nl We denote \beq{average K}\ovl{ f_\cK}({\rm X}_\cK):=\frac{1}{(2\p)^n}\int_{\torus^n}f_{\cK}({\rm X}_\cK, \ell)d\ell, \eeq so that \beqano
\begin{array}
    {lll} \dst f_{\cK}=\sum_{1\le i<j\le n}f_{\cK}^{ij},\qquad\qquad &\dst \ovl{ f_\cK}=\sum_{1\le i<j\le n}\ovl {f_{\cK}^{ij}}\\
    \\
    \dst f_{\cK}^{ij}:=\frac{y_\cK^\ppi\cdot y_\cK^\ppj}{m_0}-\frac{m_im_j}{\|x_\cK^\ppi-x_\cK^\ppj\|},&\dst \ovl {f_{\cK}^{ij}}:=\frac{1}{(2\p)^n}\int_{\torus^n}{ f_\cK^{ij}}d\ell_1\cdots d\ell_n.
\end{array}
\eeqano

 \noi For a general Kepler map, one always has, as a consequence of \equ{Kepler motion equations},  
\begin{align}
\label{yi} 
 -\frac{1}{2\p}\int_\torus T^\ppj_\cK d\ell_j
&= \frac{1}{2\p}\int_\torus \frac{V^\ppj_\cK}{2}d\ell_j=T^\ppj_\cK+V^\ppj_\cK=-\frac{{\eufm m}{\eufm M}}{2a_{j,\cK}} \notag\\
\frac{1}{2\p}\int_\torus y_\cK^\ppj d\ell_j
&= 0
\qquad\qquad
\frac{1}{2\p}\int_\torus \frac{x_\cK^\ppj}{\|x_\cK^\ppj\|^3}d\ell_j
=0, 
\end{align}
 where we have denoted as $$T^\ppj_\cK:=\frac{\|y_\cK^\ppj\|^2}{2{\eufm m}_j}\qquad V^\ppj_\cK:=-\frac{{\eufm m}_j{\eufm M}_j}{\|x_\cK^\ppj\|}$$ the kinetic, potential part of ${\rm H}^\ppj_\cK$ in \equ{1B}, respectively.

\nl Consider the average $\ovl{ f_\cK}({\rm X}_\cK)$ in \equ{average K}. Due to the fact that $y^\ppj_\cK$ has zero-average, one has that only the Newtonian part contributes to $\ovl{ f_\cK}({\rm X}_\cK)$: $$\ovl {f_\cK}=-\sum_{1\le i<j\le n}\frac{m_im_j}{(2\p)^2}\int_{\torus^2}\frac{d\ell_id\ell_j}{\|x_\cK^\ppi-x_\cK^\ppj\|}. $$ We now consider any of the contributions to this sum

$$\ovl {f_{\cK}^{ij}}=-\frac{m_im_j}{(2\p)^2}\int_{\torus^2}\frac{d\ell_id\ell_j}{\|x_\cK^\ppi-x_\cK^\ppj\|}\qquad 1\le i<j\le n $$ and expand {any such} terms \[ \ovl {f_{\cK}^{ij}}=\ovl {f_{\cK}^{ij}}^\ppo+ \ovl {f_{\cK}^{ij}}^\ppu+ \ovl {f_{\cK}^{ij}}^\ppd+\cdots\] where $$ \ovl {f_{\cK}^{ij}}^\pph:=-\frac{m_im_j}{(2\p)^2}\int_{\torus^2}\frac{1}{h!}\frac{d^h}{d\varepsilon^h}\frac{1}{\|\varepsilon x_\cK^\ppi-x_\cK^\ppj\|}\Big|_{\varepsilon=0}d\ell_id\ell_j$$ is proportional to {$\frac{1}{a_j}(\frac{a_i}{a_j})^h$}. Then the formulae in \equ{yi} imply that the two first terms of this expansion are given by $$\ovl {f_{\cK}^{ij}}^\ppo= -\frac{m_im_j}{a_{j,\cK}},\qquad \ovl {f_{\cK}^{ij}}^\ppu= 0.$$

\nl Namely, whatever is the Kepler map that is used, the first term that depends on the secular coordinates ${\rm X}_\cK$ is the double average of the second order term \[ \ovl {f_{\cK}^{ij}}^\ppd({\rm X}_\cK)=-\frac{m_im_j}{(2\p)^2}\int_{\torus^2}\frac{3(x_\cK^\ppi\cdot x_\cK^\ppj)^2-\|x_\cK^\ppi\|^2\|x_\cK^\ppj\|^2}{\|x_\cK^\ppj\|^5}d\ell_id\ell_j.\]

\nl Now we specialize to the case of the $\cP$-map. 

\section[The case of the $\cP$-map]{The case of the $\cP$-map}\label{The case of the P-map} We denote as \beq{HP}{\rm H}_\cP({\rm X}_\cP,\ell)={\rm h}_{\eufm{fast}}^0(\L)+\m f_\cP({\rm X}_\cP, \ell)\qquad {\rm X}_\cP:=(\Theta,\chi, \L,\vartheta, \k)\eeq where \beq{hk0}{\rm h}_{\eufm{fast}}^0(\L):=-\sum_{j=1}^n\frac{{\eufm m}_j^3{\eufm M}_j^2}{2\L_j^2},\eeq {is} the Hamiltonian \equ{Helio} expressed in $\cP$-coordinates.

\nl Using the definitions, it not difficult to see that
\begin{lemma}
    \label{fij av in P coordinates} $\ovl{f_\cP^{ij}}$ , $f_\cP^{ij}$ depend, respectively, only on the coordinates
    \begin{align*}
    	\ovl{{\rm X}_\cP^{ij}}
	&:= \big(\Theta_{i},\cdots,\Theta_{j\wedge(n-1)},~ \chi_{i-1},\cdots,\chi_{j\wedge (n-1)},~ \L_i,~ \L_j,~ \vartheta_{i},\cdots,\vartheta_{j\wedge (n-1)},~ \k_{i},\cdots,\k_{j-1}\big) \\
	\ovl{{\rm P}^{ij}} 
	&:= (\ovl{{\rm X}_\cP^{ij}}, \ell_i,\ell_j) 
	\end{align*}
	with $a\wedge b$ denoting the minimum of $a$ and $b$.
\end{lemma}
Accordingly to the previous lemma, the ``nearest-neighbor'' terms $\ovl {f_{\cP}^{i, i+1}}$, with $i=1$, $\cdots$, $n-1$, depend only on 
\begin{equation}\label{XX} 
\begin{split}\ovl{{\rm X}_\cP^{i, i+1}}=\left\{
\begin{array}
    {lll} \big(\Theta_{i} ,\ \Theta_{i+1} ,\ \chi_{i-1} ,\ \chi_i ,\ \chi_{i+1} ,\ \L_i ,\ \L_{i+1} ,\ \vartheta_{i} ,\ \vartheta_{i+1} ,\ \k_{i}\big) &\dst n\ge 3\ \&\ \atop{\dst 1\le i\le n-2}\\
    \\
    \big(\Theta_{n-1} ,\ \chi_{n-2} ,\ \chi_{n-1} ,\ \L_{n-1} ,\ \L_{n},\ \vartheta_{n-1} ,\ \k_{n-1}\big) &i=n-1.
\end{array}
\right.
\end{split} \end{equation}

\nl However, for the functions $\ovl {f_{\cP}^{i, i+1}}^\ppd$, we have a special rule. Indeed, for any Kepler map $\cK$, the ``exterior'' angular momentum $\|{\rm C}^{(i+1)}_\cK\|$ is an integral for $\ovl {f_{\cK}^{i, i+1}}^\ppd$. This readily implies that any $\ovl {f_{\cK}^{i, i+1}}^\ppd$ is {\sl integrable}, for having four degrees of freedom and four independent, commuting integrals ($\|{\rm C}^{(i)}_\cK+{\rm C}^{(i+1)}_\cK\|$, $({\rm C}^{(i)}_\cK+{\rm C}^{(i+1)}_\cK)\cdot k^\ppt$, $\|{\rm C}^{(i+1)}_\cK\|$ and $\ovl {f_{\cK}^{i, i+1}}^\ppd$ itself). This fact has been firstly noticed, in the three-body case ($i=1$, $n=2$), by R. Harrington \cite{harrington69} who, using the Jacobi reduction of the nodes $\cJ ac$, where the coordinates are named $${\rm G}_i,\qquad {\rm g}_i,\qquad \L_i,\qquad \ell_i,\qquad i=1,\ 2$$ (with ${\rm G}_i=\|{\rm C}^\ppi\|$, ${\rm g}_i$ related to the perihelia directions, and ${\rm G}:=\|{\rm C}\|$, ${\rm C}={\rm C}^\ppu+{\rm C}^\ppd$ appearing as an external parameter){,} noticed that $\ovl{f^{12}_{\cJ ac}}^\ppd$ depends only on $({\rm G}, {\rm G}_1, {\rm G}_2, \g_1, \L_1, \L_2)$. 

\nl Let us now inspect how the integrability of $\ovl {f_{\cP}^{i, i+1}}^\ppd$ is exhibited in terms of the $\cP$-map. Since $\|{\rm C}_\cP^\ppn\|=\chi_{n-1}$, one has that $\ovl {f_{\cP}^{n-1, n}}^\ppd$ does not {depend on} $\k_{n-1}$, and hence, by \equ{XX}, depends only on $$\ovl{\ovl{{\rm X}_\cP^{n-1,n}}}:=\big(\Theta_{n-1} ,\ \chi_{n-2} ,\ \chi_{n-1} ,\ \L_{n-1} ,\ \L_{n},\ \vartheta_{n-1}\big).$$ This fact, for $n\ge 3$, is no longer true for $i=1$, $\cdots$, $n-2$, because in that case $\chi_i\neq\|{\rm C}_\cP^{(i+1)}\|$ (indeed, $\chi_i=\|{\rm S}_\cP^{(i+1)}\|$). However, since, for $(\Theta_{i+1}, \vartheta_{i+1})=(0,\p)$, $\|{\rm C}^{(i+1)}_\cP\|$ reduces to $$\|{\rm C}^{(i+1)}_\cP\|\Big|_{(\Theta_{i+1}, \vartheta_{i+1})=(0,\p)}=\chi_i-\chi_{i+1}\qquad i=1,\ \cdots,\ n-2,$$ one has that the functions $$\ovl{\ovl {f_{\cP}^{i, i+1}}}^{(2)}:=\ovl {f_{\cP}^{i, i+1}}^\ppd\Big|_{(\Theta_{i+1}, \vartheta_{i+1})=(0,\p)},\qquad i=1 ,\ \cdots,\ n-2$$ do not depend on $\k_i$ and hence, by \equ{XX} depend only on $$\ovl{\ovl{{\rm X}_\cP^{i,i+1}}}:=\big(\Theta_{i} ,\ \chi_{i-1} ,\ \chi_i ,\ \chi_{i+1} ,\ \L_i ,\ \L_{i+1} ,\ \vartheta_{i} \big),\qquad i=1,\ \cdots,\ n-2.$$

\nl In the following lemma we provide their explicit expressions.
\begin{lemma}
    \label{integrabilityyyy} The function $\ovl {f_{\cP}^{n-1, n}}^\ppd$ and , for $n\ge 3$ and $1\le i\le n-2$, the functions $\ovl{\ovl {f_{\cP}^{i, i+1}}}^{(2)}$ have the following expressions 

 \begin{equation}\label{fn-1n} \begin{split}&\ovl {f_{\cP}^{n-1, n}}^\ppd=m_{n-1}m_n \frac{a_{n-1}^2}{4a_n^3}\frac{\L_n^3}{\chi_{n-1}^5}\Big[ \frac{5}{2}(3\Theta_{n-1}^2-\chi_{n-1}^2)\\
 &-\frac{3}{2}\frac{4\Theta_{n-1}^2-\chi_{n-1}^2}{\L_{n-1}^2}\Big(\chi_{n-2}^2+\chi_{n-1}^2-2\Theta_{n-1}^2+2\sqrt{(\chi_{n-1}^2-\Theta_{n-1}^2)(\chi_{n-2}^2-\Theta_{n-1}^2)}\cos{\vartheta_{n-1}}\Big)\\
   &+\frac{3}{2}\frac{(\chi_{n-1}^2-\Theta_{n-1}^2)(\chi_{n-2}^2-\Theta_{n-1}^2)}{\L_{n-1}^2}\sin^2{\vartheta_{n-1}} \Big] \end{split}\end{equation} 
    and 
    \beqa{ovl f} \ovl{\ovl {f_{\cP}^{i, i+1}}}^{(2)}&=&m_{i} m_{i+1} \frac{a_{i}^2}{4a_{i+1}^3}\frac{\L_{i+1}^3}{\chi_i^2(\chi_i-\chi_{i+1})^3}\Big[ \frac{5}{2}(3\Theta_i^2-\chi_i^2)\nonumber\\
    &-&\frac{3}{2}\frac{4\Theta_i^2-\chi_i^2}{\L_{i}^2}\Big(\chi_{i-1}^2+\chi_i^2-2\Theta_i^2+2\sqrt{(\chi_i^2-\Theta_i^2)(\chi_{i-1}^2-\Theta_i^2)}\cos{\vartheta_i}\Big)\nonumber\\
    &+&\frac{3}{2}\frac{(\chi_i^2-\Theta_i^2)(\chi_{i-1}^2-\Theta_i^2)}{\L_{i}^2}\sin^2{\vartheta_i}\Big]. \eeqa
\end{lemma}
\vskip.1in \noi

\nl Lemma \ref{integrabilityyyy} is proved in Appendix \ref{Proof of Lemma {integrabilityyyy}}. Here, we limit to the following.
\begin{remark} \rm  ~\\[-1em]
\begin{enumerate} 
    \rm \item[{\rm (i)}] The formula in \equ{ovl f} holds also for complex values of the coordinates, provided that \\ $\arg(\chi_i-\chi_{i+1})\in (-\frac{\p}{2},\ \frac{\p}{2}]$ mod $2\p$.
    \item[{\rm (ii)}] The importance of the formulae in \equ{fn-1n} and \equ{ovl f}, which is the main feature of the $\cP$-map, is that, exploiting the equilibrium for $(\Theta_i,\vartheta_i)=(0,\p)$, the integration of $\ovl {f_{\cP}^{n-1, n}}^\ppd$ and of $\ovl{\ovl {f_{\cP}^{i, i+1}}}^{(2)}$ can be performed {\sl explicitly}, switching to a suitable associated {\sl convergent} Birkhoff series, as Lemma \ref{integration***} below states. Direct integrations of $\ovl {f_{\cK}^{n-1, n}}^\ppd$, for example, starting with Hamiltonian computed in \cite{harrington69}, appear technically much more involved and, up to now, are not known.
\end{enumerate}
\end{remark}

\begin{lemma}
    \label{integration***} It is possible to find complex domains $\ovl{{\eufm B}_i}$ with non-empty real part and a canonical, real-analytic change of coordinates \beqa{tr*} \ovl{\phi_{\eufm{int}}^{i}}:\quad (p_i, q_i, {\rm y}^*_i, {\rm x}^*_i)\in\ovl{{\eufm B}_i} \to (\Theta_i, \vartheta_i, {\rm y}_i, {\rm x}_i) \eeqa where \beqa{yx} {\rm y}^*_i&:=&\arr{ \dst (\chi^*_{n-2}, \chi^*_{n-1}, \L^*_{n-1}, \L^*_{n})\qquad\qquad i=n-1\\
    \\
    \dst(\chi^*_{i-1}, \chi^*_i, \chi^*_{i+1}, \L^*_i, \L^*_{i+1})\qquad\ \ \ i=1,\cdots, n-2\ \&\ n\ge 3 }\nonumber\\
    \nonumber\\
    {\rm x}^*_i&:=&\arr{ \dst (\k^*_{n-2}, \k^*_{n-1}, \ell^*_{n-1}, \ell^*_{n})\qquad \qquad\ i=n-1\\
    \\
    \dst(\k^*_{i-1}, \k^*_i, \k^*_{i+1}, \ell^*_i, \ell^*_{i+1})\qquad\quad \ \ i=1,\cdots, n-2\ \&\ n\ge 3 }\nonumber\\
    \nonumber\\
    {\rm y}_i&:=&\arr{ \dst (\chi_{n-2}, \chi_{n-1}, \L_{n-1}, \L_{n})\qquad\quad \ \ i=n-1\\
    \\
    \dst(\chi_{i-1}, \chi_i, \chi_{i+1}, \L_i, \L_{i+1})\qquad\ \ \ i=1,\cdots, n-2\ \&\ n\ge 3 }\nonumber\\
    \nonumber\\
    {\rm x}_i&:=&\arr{ \dst (\k_{n-2}, \k_{n-1}, \ell_{n-1}, \ell_{n})\qquad \qquad\ i=n-1\\
    \\
    \dst(\k_{i-1}, \k_i, \k_{i+1}, \ell_i, \ell_{i+1})\qquad\quad \ \ \ i=1,\cdots, n-2\ \&\ n\ge 3 } \eeqa such that \beq{nude integration}\ovl{{\rm h}_{\eufm{sec}}^{i}}:=\arr{ \dst {\ovl {f_{\cP}^{n-1, n}}}^{(2)}\circ\ovl{\phi_{\eufm{int}}^{n-1}}\qquad i=n-1\\
    \\
    \dst\ovl{\ovl {f_{\cP}^{i, i+1}}}^{(2)}\circ\ovl{\phi_{\eufm{int}}^{i}}\qquad\ \ \ i=1,\cdots, n-2\ \ \&\ \ n\ge 3 } \eeq depends only on  \[{\rm Y}^*_i:=\arr{ \dst (\tfrac{p_{n-1}^2+q_{n-1}^2}{2}, \L^*_{n-1}, \L^*_{n}, \chi^*_{n-2}, \chi^*_{n-1})\qquad \ \ i=n-1\\
    \\
    \dst(\tfrac{p_i^2+q_i^2}{2}, \L^*_i, \L^*_{i+1}, \chi^*_{i-1}, \chi^*_i, \chi^*_{i+1})\qquad\qquad i=1,\cdots, n-2\ \&\ n\ge 3. } \] The transformation $\ovl{\phi^i_{\eufm{int}}}$ may be chosen so as to verify \beqa{symmetries}&&{\rm y}_i^*={\rm y}_i,\quad (\Theta_i, \ \vartheta_i,\ {\rm x}_i-{\rm x}_i^*)={\rm F}_i(p_i, q_i, {\rm y}_i^*)\nonumber\\
    \nonumber\\
    && \ovl{\phi^i_{\eufm{int}}}(-p_i, -q_i, {\rm y}_i^*, {\rm x}_i^*)=(-\Theta_i, -\vartheta_i, {\rm y}_i, {\rm x}_i)\eeqa if \[ \ovl{\phi^i_{\eufm{int}}}(p_i, q_i, {\rm y}_i^*, {\rm x}_i^*)=(\Theta_i, \vartheta_i, {\rm y}_i, {\rm x}_i).\]
\end{lemma}
Lemma \ref{integration***} is proved in the following Section \ref{Local coordinates}. \chapter[Global Kolmogorov tori]{Global Kolmogorov tori in the planetary problem} In this section we show how the $\cP$-map can be used to prove Theorem A. We defer to the next Section \ref{Proofs} more technical parts.

\section{A domain of holomorphy}\label{A sub-domain of holomorphy} A typical practice, in order to use perturbation theory techniques, is to extend Hamiltonians governing dynamical systems to the complex field, and then to study their holomorphy properties.

\nl In this section we aim to discuss a domain of holomorphy for the perturbing function $f_\cP$ in \equ{HP}, regarded as a function of complex coordinates. We shall choose it of the following form

 $${\eufm D}_{\cP}:={\rm T}_{\Theta^+,\vartheta^+}\times\big({\rm X}_\theta\times \ovl\torus^n_{{s}}\big)\times\big({\rm A}_\theta\times \ovl\torus^{n}_{{s}}\big),$$ where, for given positive numbers \beqano \Theta_{j}^+,\quad \vartheta_{j}^+,\quad{\rm G}_i^\pm,\quad \L_i^\pm,\quad \theta_i,\quad s \eeqano with $i=1$, $\cdots$, $n$, $j=1$, $\cdots$, $n-1$, \begin{equation}\label{domain***} \begin{split}{\rm T}_{\Theta^+,\vartheta^+}:=&\Big\{(\ovl\Theta, \ovl\vartheta)=(\Theta_1,\cdots, \Theta_{n-1}, \vartheta_1,\cdots,\vartheta_{n-1})\in \complex^{n-1}\times \torus_\complex^{n-1}:\\
 &|\vartheta_j-\p|\le {\vartheta^+_j},\quad |\Theta_{j}|\le {\Theta_j^+},\quad \forall\ j=1,\cdots, n-1\Big\}\\
\\
{\rm X}_\theta:=&\Big\{ \chi=(\chi_0,\cdots, \chi_{n-1})\in \complex^n:\ {\rm G}_j^-\le |\chi_{j-1}-\chi_j|\le{\rm G}_j^+,\quad |\Im (\chi_{j-1}-\chi_j)|\le \theta_j\\
 &\forall\ j=1,\cdots, n\Big\}\\\\
{\rm A}_\theta:=&\Big\{\L=(\L_1,\cdots, \L_n)\in \complex^n:\quad \L_j^-\le | \L_j|\le \L_j^+,\quad |\Im \L_j|\le \theta_j\\
 &\forall\ j=1,\cdots, n\Big\}\\
\ovl\torus_{{s}}:=&\torus+{\rm i}[-{s}, {s}]\end{split}\end{equation} with $\chi_n:=0$.

\nl The domain ${\eufm D}_\cP$ will be determined as the intersection of the ``collision-less'' set, where, as functions of complex variables, the mutual distances of the planets \[ d_{j,\cP}:=\|x^\ppj_\cP-x^{(j+1)}_\cP\|\] are far away from zero, {with} the holomorphy domain of $\cP$, where, again as as functions of complex variables, the absolute values $|e_{j,\cP}|$ of eccentricities in \equ{P*ecc new} are bounded away from $0$ and $1$, those of the inclinations $|\iota_j|$, $|{\rm i}_j|$ in \equ{good incli*} are away from $0$ and, finally, {\sl Kepler equation} \equ{Kep Eq} provides a holomorphic solution.

\nl The latter issue is not a peculiarity of this problem, since it naturally arises in the context of the two-body problem's equations. In the early XX century, T. Levi Civita \cite{leviCivita1904} studied the holomorphy of the solution of Kepler's Equation with respect to the eccentricity. The holomorphy with respect to the mean anomaly has been investigated, using similar arguments as in \cite{leviCivita1904}, in \cite{cellettiPi05}. Here, we address the problem of determining the holomorphy with respect to both the arguments.
\begin{proposition}
    \label{Kepler equation*} Let $\widehat e=\;${$0.6627\ldots$} be the solution of \beq{hat e}0\le \r\le 1\quad \&\quad\frac{\r\, e^{\sqrt{1+\r^2}}}{1+\sqrt{1+\r^2}}=1.\eeq Then for any $0<\ovl e<\widehat e$, one can find a positive number $\bar\ell$ depending on $\ovl e$ such that, for any ${\rm e}={\rm e}_1+{\rm i}{\rm e}_2\in \complex$, with $|{\rm e}|\le \ovl e$, the complex {Kepler equation} $$\zeta-{\rm e}\sin\zeta=\ell$$ has a unique solution $\zeta(\ell, {\rm e})$ which turns out to be real-analytic for $\ell\in \ovl{\torus}_{\bar\ell}$.
\end{proposition}
\nl The following result completes the study of the holomorphy of $f_\cP$.
\begin{proposition}
    \label{prop: domain} Let $\widehat e$ be as in Proposition \ref{Kepler equation*}. For any given $\underline e_i$, $\overline e_i$, with \[0<\underline e_i<\overline e_i<\widehat e\qquad i=1,\cdots, n\] it is possible to find positive numbers $$\cA_j,\quad \cB_j,\quad \ovl{\rm C}_i>\underline{\rm C}_i,\quad \bar d_j,\quad {s}\in (0,1),\quad \s\in (0,1)$$ such that, if the following inequalities are satisfied 
\begin{equation}\label{initial bounds}
\begin{split}
	&\underline{\rm C}_i\L^+_i<{\rm G}^-_i< {\rm G}^+_i<\ovl{\rm C}_i\L^-_i;\\
    &\max\bigg\{~ \frac{\theta_i}{\L_i^-},~~ \frac{\theta_{i}}{{\rm G}_{i}^-},~~ \sum_{i=1}^{n-1}\bigg|\sin^{-1}\!\Big(\frac{{\rm G}_i^+}{{\rm G}_{i+1}^-}\bigg)\bigg|,~~ \frac{ \Theta_{j}^+}{{\rm G}_n^-},~~ \sum_{i=1}^{n-1}\frac{{\rm G}_i^+}{{\rm G}_n^-},~~ \vartheta_j^+,~~ |\Im\k_j|,~~ |\Im \ell_i| ~\bigg\}\le s \\
    &\vartheta_{j}^+\le\min\bigg\{\frac{\cA}{{\rm G}_n^+}\sqrt{({\rm G}_j^-)^2-(\underline{\rm C}_j\L^+_j)^2},~ \frac{\cB}{{\rm G}_n^+}\sqrt{(\ovl{\rm C}_j\L^-_j)^2-({\rm G}_j^+)^2}\bigg\}, 
    \end{split}
    \end{equation}
then the eccentricities $e_{i,\cP}$, inclinations $\iota_i$, ${\rm i}_i$ and the mutual distances $d_{i,\cP}$ verify \beq{eccentricities} \underline e_i\le |e_{i,\cP}|\le\ovl e_i,\quad \max_{i,j}\Big\{|\cos\iota_i|,\ |\cos{\rm i}_j|\Big\}\le\s,\quad |d_{i,\cP}|\ge \bar d \eeq
\end{proposition}
Proposition \ref{Kepler equation*} and Proposition \ref{prop: domain} are proved in Appendix \ref{Computing the domain of holomorphy} and \ref{Computing the domain of holomorphy***}, respectively. We shall use them in the form below. We remark that the super-exponential decay of the semi-major axes ratio will be used only in Section \ref{Secular normalizations} below.

\nl
\begin{corollary}
    [choice of parameters] \label{cor: parameters} Fix $\underline e_i<\overline e_i$, $c\in(0,1)$, and let $\underline\cC_i<\underline\cC_i^*<\overline\cC_i^*<\overline\cC_i$, \mbox{$\cD_i:=\min\{\cA\sqrt{(\underline\cC_i^*)^2-(\underline\cC_i)^2}$}, $\cB\sqrt{(\ovl\cC_i)^2-(\ovl\cC^*_i)^2}\}$, $\cD:=\min_{1\le j\le n-1}\frac{\cD_j}{\ovl\cC^*_n}\frac{{\eufm m}_j\sqrt{{\eufm M}_j}}{{\eufm m}_n\sqrt{{\eufm M}_n}}$, $\a<\frac{s}{\cD}$. Define, for $i=1,\cdots,n$ and $j=1,\cdots,n-1$,
\begin{equation}\label{Choice of parameters}\begin{split}
	\L_i^\pm:=&{\eufm m}_i\sqrt{{\eufm M}_ia_i^\pm},\quad {\rm G}_i^+:=\ovl\cC_i^* \L_i^-,\qquad {\rm G}_i^-:=\underline\cC_i^* \L_i^+,\quad\Theta_j^+:=s{\rm G}_n^-,\qquad \vartheta_j^+:=\cD_i\frac{{\L_i^-}}{{\rm G}_n^+}\\
    \theta_i:=&s\sqrt{\L_i^-}
    \end{split}
\end{equation}
where $a_{i}^\pm$ is as in $(*)$. Then $f_\cP$ is real-analytic in the domain ${\eufm D}_\cP$.
\end{corollary}
\section{A normal form for the planetary problem}
\begin{definition}
    [\cite{arnold63}]\rm Given $m$, $\n_1$, $\cdots$, $\n_m\in \natural$, $\n:=\n_1+\cdots+\n_m$, let $${\eufm L}_0\supset{\eufm L}_1\supset{\eufm L}_2\supset \cdots \supset {\eufm L}_m=\{0\}$$ be a decreasing sequence of sub-lattices of $\integer^\n$ defined by \beq{Li} {\eufm L}_0:=\integer^\n,\quad {\eufm L}_i:=\big\{k=(k_1,\cdots, k_{m})\in \integer^\n,\quad k_j\in\integer^{\n_j} : \quad k_1=\cdots=k_{i}=0\big\}\eeq with $i=1$, $\cdots$, $n$. Next, given $\g$, $\g_1$, $\cdots$, $\g_m$, $\t\in \real_+$, we define the set ${\rm D}^{\n}_{\g_1\cdots\g_m;\t}$ of the $(\g_1\cdots\g_m;\t)$-{\sl diophantine numbers} via the following formulae: \beqano &&{\rm D}^{\n, K, i}_{\g;\t}:=\Big\{\o\in \real^\n:\ |\o\cdot k|\ge \frac{\g}{|k|^\t}\quad \forall k\in {\eufm L}_{i-1}\setminus{\eufm L}_i,\ |k|_1\le K\Big\}\nonumber\\
    &&{\rm D}^{\n, K}_{\g_1\cdots\g_m;\t}:=\bigcap_{i=1}^m {\rm D}^{\n, K, i}_{\g_i;\t}\qquad {\rm D}^{\n}_{\g_1\cdots\g_m;\t}:=\bigcap_{K\in \natural} {\rm D}^{\n, K}_{\g_1\cdots\g_m;\t} .\eeqano In other words $\o=(\o_1,\cdots,\o_m)\in{\rm D}^{\n}_{\g_1\cdots\g_m;\t}$ if, for any $k=(k_1,\cdots,k_m)\in \integer^\n\setminus\{0\}$, with $k_j\in \integer^{\n_j}$, 
    \beq{dioph2sc} |\o\cdot k|=\bigg|\sum_{j=1}^m\o_j\cdot k_j\bigg|\geq \left\{
    \begin{array}
        {l} \dst\frac{\fg}{|\tk|^{\t}}\quad {\rm if}\quad \fk\neq0;\\
        \ \\
        \dst \frac{\g_2}{|k|^{\t}}\quad {\rm if}\quad \fk= 0,\quad \sk\neq 0;\\
        \ \\
        \;\;\;\vdots\\
        \ \\
        \dst \frac{\g_m}{|k_m|^{\t}}\quad {\rm if}\quad k_1=\cdots=k_{m-1}= 0,\ \cdots,\ k_{m}\neq 0.
    \end{array}
    \right. \eeq
\end{definition}

\begin{remark}
    \rm The choice $m=1$, $\g_1:=\g$ gives the usual Diophantine set $\cD^\n_{\g,\t}$. The $m=2$-case, ${\mathcal D}^{\n}_{\g_1,\g_2,\t}$, with $\g_1={\rm O}(1)$ and $\g_2={\rm O}(\m)$, where $\m$ is the strength of the planetary masses has been considered in \cite{arnold63} for the proof of the Fundamental Theorem, mentioned in the introduction.
\end{remark}
\nl The following result in proven in the next Section \ref{Proofs}. It is unavoidably detailed.
\begin{proposition}
    \label{exponential average} Let ${\eufm m}_j$, ${\eufm M}_j$ be as in \equ{reduced masses} and ${\rm m}_j:=\sum_{i=1}^{j-1} m_i$, with $j=2,\cdots, n$. There exists a number ${\eufm c}$, depending only on $n$, $m_0,\cdots,m_n$, $a_n^\pm$, $\underline e_j$, $\ovl e_j$, and a number $0<\ovl{\eufm{c}}<1$, depending only on $n$ such that, for any fixed positive numbers $\ovl\g<1<\bar K$, $\a>0$ verifying \beqa{bar K} &&\bar K\le \frac{{\eufm c} }{\a^{3/2}}\eeqa and \beqa{small secular}\frac{1}{{\eufm c}}\max\Big\{ \m(\frac{a_n^+}{a_1^-})^{5}\frac{\bar K^{2\bar\t+2}}{\bar\g^2},\ \frac{\bar K^{2(\bar\t+1)}\a}{\bar\g^2}\Big\}<1\eeqa there exist natural numbers $\n_1,\cdots,\n_{2n-1}$, with $\sum_j\n_j=3n-2$, open sets $B_j^*\subset B^{2}_{\varepsilon_j}, {\rm X}^*\subset {\rm X}$, positive real numbers \mbox{$\g_1> \cdots >\g_{2n-1} \varepsilon_1, \cdots, \varepsilon_{n-1}, \ovl r_1, \cdots, \ovl r_{n-1}, \widetilde r_1, \cdots, \widetilde r_{n}$}, a domain $${\eufm D}_{\rm n}:=B_{\sqrt{2\ovl r}}\times {\rm X}_{\ovl r}\times\cA_{\widetilde r} \times\torus^n_{\ovl{\eufm{c}}s}\times \torus^n_{\ovl{\eufm{c}}s}$$ a sub-domain of the form $${\eufm D}^*_{\rm n}:=B^*_{\sqrt{2\ovl r}}\times {\rm X}^*_{\ovl r}\times\cA_{\widetilde r} \times\torus^n_{\ovl{\eufm{c}}s}\times \torus^n_{\ovl{\eufm{c}}s}$$ verifying \beq{good set***}\meas{\eufm D}^*_{\rm n}\ge\big(1-\frac{\bar\g}{\ovl{\eufm c}}\big)\meas{\eufm D}_{\rm n}\eeq a real-analytic transformation $$\phi_{\rm n}:\quad (p,q,\chi,\L,\k,\ell)\in {\eufm D}^*_{\rm n}\to {\eufm D}_\cP$$ which conjugates ${\rm H}_\cP$ to $${\rm H}_{\rm n}(p,q, \chi,\L,\k,\ell) :={\rm H}_\cP\circ\phi_{\rm n}={\rm h}_{\eufm{fast}, \eufm{sec}}(p,q,\chi,\L)+\m\,{f}_\eufm{exp}(p,q, \chi,\L,\k,\ell) $$ where ${f}_\eufm{exp}(p,q, \chi,\L,\k,\ell)$ is independent of $\k_0$, and the following holds.

\vskip.1in\noi\paragraph{\bf 1.} The function ${\rm h}_{\eufm{fast}, \eufm{sec}}(p,q,\chi,\L)$ is a sum $${\rm h}_{\eufm{fast}, \eufm{sec}}(p,q,\chi,\L)={\rm h}_{\eufm{fast}}(\L)+\m\, {\rm h}_{\eufm{sec}}(p,q,\chi,\L)$$ where, if \beqano \hat{\rm y}_i:=\bigg(\frac{p_i^2+q_i^2}{2},\ \cdots,\ \frac{p_{n-1}^2+q_{n-1}^2}{2},\ \chi_{i-1},\ \cdots,\ \chi_{n-1},\ \L_i,\ \cdots,\ \L_n\bigg) \eeqano then ${\rm h}_{\eufm{fast}}$ and ${\rm h}_{\eufm{sec}}$ are given by $${\rm h}_{\eufm{fast}}(\L)=-\sum_{j=1}^n\frac{{\eufm m}_j^3{\eufm M}_j^2}{2\L_j^2}-\m\sum_{j=2}^n\frac{{\eufm M}_j{\eufm m}_j^2m_j{\rm m}_j}{\L_j^2} ,\qquad {\rm h}_{\eufm{sec}}(p,q,\chi,\L)=\sum_{i=1}^{n-1}{\rm h}_{\eufm{sec}}^i(\hat{\rm y}_i)$$ where the functions ${\rm h}_{\eufm{sec}}^i$ have an analytic extension on ${\eufm D}_{\rm n}$ and verify \[{\eufm c}\frac{(a_{n-j}^+)^2}{(a_{n-j+1}^-)^3}\le | {\rm h}_{\eufm{sec}}^j(\hat{\rm y}_j) |\le \frac{1}{\eufm c}\frac{(a_{n-j}^+)^2}{(a_{n-j+1}^-)^3}.\]

\vskip.1in\noi\paragraph{\bf 2.} The function ${f}_\eufm{exp}$ satisfies $$|{f}_\eufm{exp}|\le \frac{1}{\eufm c}\frac{e^{-{\eufm c}\bar K}}{a_1^-}.$$ \vskip.1in\noi\paragraph{\bf 3.} If $\zeta$ is $\hat{\rm y}_1$ deprived of $\chi_0$, the frequency-map $$\zeta\to\o_{\eufm{fast}, \eufm{sec}}(\zeta):= \partial_{\zeta}{\rm h}_{\eufm{fast}, \eufm{sec}}(\zeta)$$ is a diffeomorphism of $\P_\zeta(B^*_{\sqrt{2\ovl r}}\times {\rm X}^*_{\ovl r}\times\cA^*_{\widetilde r})$ and, moreover, it satisfies \equ{dioph2sc}, with $m=2n-1$, $\t=\bar\t>2$, and
 \begin{equation}\label{nugamma}  \begin{split}
 \n_j:=&\left\{
    \begin{array}
        {llll} \dst\ 1& j=1,\cdots, n\\
        \\
        \dst \ 2\qquad &j=3,\ n=2\\
        \\
        \dst \ 3 & j=n+1,\ n\ge 3 \\
        \\
        \dst \ 2& n+2\le j\le 2n-2,\ n\ge 4 \\
        \\
        \dst\ 1& j=2n-1,\ n\ge 3
    \end{array}
    \right.
\\\\
    \omega_j:=&\left\{
    \begin{array}
        {llll} \dst\partial_{\L_j}{\rm h}_{\eufm{fast}, \eufm{sec}}& j=1,\cdots, n\\
        \\
        \dst \partial_{(\frac{p_{1}^2+q_{1}^2}{2},\chi_1)}\,{\rm h}_{\eufm{fast}, \eufm{sec}}\qquad &j=3,\ n=2\\
        \\
        \dst \partial_{(\frac{p_{n-1}^2+q_{n-1}^2}{2}, \chi_{n-2},\chi_{n-1})}\,{\rm h}_{\eufm{fast}, \eufm{sec}}& j=n+1,\ n\ge 3 \\
        \\
        \dst \partial_{(\frac{p_{2n-j}^2+q_{2n-j}^2}{2}, \chi_{2n-j-1})}\,{\rm h}_{\eufm{fast}, \eufm{sec}}& n+2\le j\le 2n-2,\ n\ge 4 \\
        \\
        \dst \partial_{\frac{p_{1}^2+q_{1}^2}{2}}\,{\rm h}_{\eufm{fast}, \eufm{sec}}& j=2n-1,\ n\ge 3
    \end{array}
    \right.\\\\
    \g_j:=& \left\{
    \begin{array}
        {llll} \dst\frac{1}{a_j^-}\frac{\ovl\g}{\theta_j} \qquad &1\le j\le n\\
        \\
        \dst\frac{\m(a_{j-n}^+)^2}{(a_{j+1-n}^-)^3}\frac{\ovl\g}{\theta_{j-n}} &n+1\le j\le 2n-1
    \end{array}
    \right.\end{split}\end{equation}

\vskip.1in\noi\paragraph{\bf 4.} The mentioned constants are \beqano \varepsilon_j:={\eufm c}\,\sqrt{\theta_j},\quad \ovl r_j:=\frac{\theta_j\ovl\g}{\bar K^{\bar\t+1}} ,\quad \widetilde r_i:={\eufm c}\,\theta_j \eeqano with $\bar\t>2$.
\end{proposition}
\section{A ``multi-scale'' KAM Theorem and proof of Theorem A} In this section we state a ``multi-scale'' KAM Theorem and next we show how this theorem applies to the Hamiltonian ${\rm H}_{\rm n}$ so as to obtain the proof of Theorem A.
\begin{theorem}
    [Multi-scale KAM Theorem]\label{two scales KAM} 
    Let $m,\ell,\n_1,\cdots,\n_m\in \natural$, $\n:=\n_1+\cdots+\n_m\ge \ell$, $\t_*>\n$, $\g_1\ge \cdots\ge\g_m>0$, $0<4s\leq \bar{s}<1$, $\r_1, \cdots, \r_\ell, r_1, \cdots, r_{\n-\ell}, \varepsilon_1, \cdots, \varepsilon_\ell>0$, $B_1, \cdots, B_\ell\subset \real^2$, $D_j:=\{\frac{x^2+y^2}{2}\in \real: (x,y)\in B_j\}\subset \real$, $B:=B_1\times\cdots\times B_\ell\subset \real^{2\ell}$, $D:=D_1\times\cdots\times D_\ell\subset \real^\ell$, $C\subset \real^{\n-\ell}$, $A:=D_\r\times C_r$.
     Let \beqano {\rm H} (p,q, I,\psi)={\rm h}(p,q,I)+{f}(p,q, I,\psi) \eeqano be real-analytic on $B_{\sqrt{2\r}}\times C_r\times \torus_{\bar{\eufm s}+s}^{\n-\ell}$, where ${\rm h}(p,q,I)$ depends on $(p,q)$ only via \[ J(p,q):=\Big(\frac{p_1^2+q_1^2}{2},\ \cdots,\ \frac{p_\ell^2+q_\ell^2}{2}\Big).\] Assume that $\o_0:=\partial_{(J(p,q),I)} {\rm h}$ is a diffeomorphism of $A$ with non singular Hessian matrix $U:=\partial^2_{(J(p,q),I)}{\rm h}$ and let $U_k$ denote the $ (\n_k+\cdots+\n_m)\times \n$ submatrix of $U$, \ie, the matrix with entries $(U_k)_{ij}=U_{ij}$, for $\n_{1}+\cdots+\n_{k-1}+1\leq i\leq \n$, $1\leq j\leq \n$, where $2\le k\le m$. Let \beqano && {\rm M}\geq\sup_{A}\|U\|,\quad {\rm M}_k\geq\sup_{A}\|U_k\|,\quad \bar {\rm M} \geq\sup_{A}\|U^{-1}\|,\quad \pertnorm\geq\|{f}\|_{\r,\bar{\eufm s}+s}\nonumber\\
    &&\bar {\rm M}_k\geq \sup_{A}\|T_k\|\quad {\rm if}\quad \dst U^{-1}=\left(
    \begin{array}
        {lrr} T_1\\
        \vdots\\
        T_m
    \end{array}
    \right)\qquad 1\le k\le m.\eeqano Define \beqano && \dst K:=\frac{6}{s}\ \log_+{\left(\frac{\pertnorm {\rm M}_1^2\,L}{\gamma_1^2}\right)^{-1}}\quad {\rm where}\quad \log_+ a :=\max\{1,\log{a}\}\\
    && \dst \hat\r_k:=\frac{\g_k}{3{\rm M}_kK^{\t_*+1}},\quad \hat\r:=\min\left\{\hat\r_1,\ \cdots,\ \hat\r_m,\ \r_1,\ \cdots,\ \r_\ell,\ r_1,\ \cdots ,\ r_{\n-\ell}\right\}\\
    \\
    && \dst L:=\max \Big\{\bar {\rm M} , \ {\rm M}_1^{-1},\ \cdots,\ {\rm M}_m^{-1}\Big\} \\
    && \hat E:=\frac{E L}{\hat\r^2}. \eeqano Then one can find two numbers $\hat c_\n>c_\n$ depending only on $\n$ such that, if the perturbation ${f}$ {is} so small that the following ``KAM condition'' holds \[ \hat c_\n\KAM<1, \] for any $\o\in\O_*:=\o_0({D})\cap\cD_{\g_1,\cdots,\g_m,\t_*}$, one can find a unique real-analytic embedding \beqano \phi_\o:\quad \vartheta=(\hat\vartheta,\bar\vartheta)\in\torus^{\n } &\to&(\hat v(\vartheta;\o),\hat\vartheta+\hat u(\vartheta;\o), \cR_{\bar\vartheta+\bar u(\vartheta;\o)}w_1,\ \cdots,\ \cR_{\bar\vartheta+\bar u(\vartheta;\o)}w_\ell)\nonumber\\
    &&\in \Re C_r\times \torus^{\n-\ell}\times \Re B^{2\ell}_{\sqrt{2r}} \eeqano where $r:= c_\n \KAM \hat\r$ such that ${\rm T}_\o:=\phi_{\o}(\torus^\n)$ is a real-analytic $\n$-dimensional ${\rm H}$-invariant torus, on which the ${\rm H}$-flow is analytically conjugated to $\vartheta\to \vartheta+\o\,t$. Furthermore, the map $(\vartheta;\o)\to\phi_\o(\vartheta)$ is Lipschitz and one-to-one and the invariant set $\dst {{\rm K}}:=\bigcup_{\o\in\O_*}{\rm T}_\o$ satisfies the following measure estimate \[
	\meas\Big(\!\Re({D}_r)\times\torus^\td\setminus{{\rm K}}\Big)
        \leq  c_\n\Big(\!\meas({D}\setminus{D}_{\g_1,\cdots, \g_m,\t_*}\times\torus^\td)+\meas(\Re({D}_r)\setminus{D})\times\torus^\td\Big),
    \] where ${D}_{\g_1,\cdots, \g_m,\t_*}$ denotes the $\o_0$-pre-image of $\cD_{\g_1,\cdots, \g_m,\t_*}$ in ${D}$. Finally, on $\torus^\n\times \O_*$, the following uniform estimates hold 
    \begin{align*}
    | v_k(\cdot;\o)-I_k^0(\o)|
    &\leq c_\n \Big(\frac{\bar {\rm M}_k}{\bar {\rm M}}+\frac{{\rm M}_k}{{\rm M}_1}\Big)\KAM\,\hat\r \\
    |u(\cdot;\o)| &\leq c_\n\KAM\,s
    \end{align*}
    where $v_k$ denotes the projection of $v=(\hat v, \bar v)\in \real^{\n_1}\times\cdots\times\real^{\n_m}$ over $\real^{\n_k}$, $\dst\bar v_k:=\frac{|w_k|^2}{2}$ and $I^0(\o) = (I^0_1(\o),\cdots, I^0_\n(\o)) \in D$ is the $\o_0$-pre-image of $\o\in\O_*$.
\end{theorem}
\nl Theorem \ref{two scales KAM} generalizes \cite[Proposition 3]{chierchiaPi10} in two respects. The {first} generalization concerns {the consideration of $m\ge 2$ scales} (in \cite{chierchiaPi10} only the case $m=2$ was treated). The {second} consists of taking ${\rm H}$ depending also on the rectangular variables $(p,q)\in B^{2\ell}$. Such generalizations can be easily obtained, and hence will be not discussed here.

\vskip.in\textsc{Proof of  Theorem A.} Let $$ \bar\g:={\ovl{\eufm c}}\sqrt\a(\log\a^{-1})^{\bar\t+1},\quad \bar K=\frac{1}{\widetilde{\eufm c}}\log\frac{1}{\a}$$ where $\ovl{\eufm c}$ is as in \equ{good set***} and $\widetilde{\eufm c}$ will be fixed later. We {aim to} apply Theorem \ref{two scales KAM} to the Hamiltonian ${\rm H}_{\rm n}$ of Proposition \ref{exponential average}, with these choices of $\bar\g$ and $\bar K$. To this end, we take \beqano &&{\rm M}_j=\left\{
\begin{array}
    {llll} \dst\frac{1}{{\eufm c}_1a_j^-\theta_j^2}\qquad\ &1\le j\le n\\
    \\
    \dst\frac{\m(a_j^+)^2}{{\eufm c}_1(a_{j+1}^-)^3\theta_j^2} & n+1\le j\le 2n-1
\end{array}
\right. \qquad L=\bar{\rm M}=\frac{1}{{\eufm c}_2}\,\frac{\theta_1^2(a_2^+)^3}{\m(a_1^-)^2} \nonumber\\
&&E=\frac{1}{{\eufm c}_3}\frac{\m}{a_1^-}e^{-{\eufm c}\bar K} \qquad\qquad\qquad\qquad\qquad \qquad\qquad\quad K=\frac{1}{{\eufm c}_4}\log_+\Big(\frac{1}{\ovl\g^2}\frac{(a_2)^3}{(a_1^-)^3}e^{-{\eufm c}\bar K}\Big)^{-1}\nonumber\\
&&\hat \r_j=\arr{\dst{{\eufm c}_5}\frac{\ovl\g\theta_j}{K^{\t_*+1}}\quad 1\le j\le n\\
\\
\dst{\eufm c}_5\frac{\ovl\g\theta_{j-n}}{{}K^{\t_*+1}}\quad n+1\le j\le 2n-1 }\qquad \qquad\hat\r:=\frac{\theta_1\ovl\g}{\hat K^{\t_*+1} }\quad \t_*>3n-2\nonumber\\
&&\hat E=\frac{1}{{\eufm c}_6}\frac{1}{\ovl\g^2}\frac{(a_2)^3}{(a_1^-)^3}e^{-{\eufm c}\bar K}\hat K^{2(\t_*+1)}\nonumber\\
\eeqano where $\hat K:=\max\{K,\bar K\}$. The number $\frac{1}{\ovl\g^2}\frac{(a_2)^3}{(a_1^-)^3}$ can be bounded by $\frac{1}{\a^N}$ for a sufficiently large $N$ depending only on $n$. Hence, if $\widetilde{\eufm c}<\frac{\eufm c}{N}$ and $\a<{\eufm c}_6$, we have $\hat E<1$ and the theorem is proved. \qed

\chapter{Proofs}\label{Proofs} In this section we provide the proof of Proposition \ref{exponential average}. This is divided in two steps: normalization of fast angles and of secular coordinates.

\section{Normalization of fast angles} Let $\ovl{f_{\cP}^{ij}}$, $\ovl{f_{\cP}^{ij}}^\ppk$ as in Lemma \ref{fij av in P coordinates}, and let \beq{split barfij}\ovl{{f}_{\cP}^{ij}}^{(\ge 2)} :=~ \ovl{{f}_{\cP}^{ij}}-\ovl{{f}_{\cP}^{ij}}^\ppo.\eeq
\begin{proposition}
    \label{Keplerian averages} There exist two small numbers $\widehat{\eufm c}$, ${\eufm c}_1$, where $\widehat{\eufm c}$ depends only on $n$, while ${\eufm c}_1$ depends only on $n$, $m_1,\cdots,m_n$, such that, if the inequality in \equ{bar K} and \beqa{fast averaging smallness assumption} \frac{1}{\eufm c}\m \bar K\bigg(\frac{a_n^+}{a_1^-}\bigg)^{\frac{3}{2}}<1 \eeqa hold, one can find a real-analytic and symplectic transformation $$\phi_{\eufm{ fast}}:\quad (\ovl\Theta,\ovl\vartheta,\chi,\L,\k,\ell)\in {\eufm D}_{\eufm{ fast}}:={\rm T}_{\widehat{\eufm c}\Theta^+,\widehat{\eufm c}\vartheta^+}\times {\rm X}_{\widehat{\eufm c}\theta}\times\cA_{\widehat{\eufm c}\theta} \times\torus^n_{\widehat{\eufm c}s}\times \torus^n_{\widehat{\eufm c}s}\to {\eufm D}_\cP$$ which conjugates ${\rm H}_\cP$ to \beq{fij fast averaged}{\rm H}_{\eufm{ fast}, \eufm{exp}}(\ovl\Theta,\chi,\L,\ovl\vartheta,\k,\ell):={\rm H}_\cP\circ\phi_{\eufm{fast}}={\rm h}_{\eufm{fast}}(\L)+\m\,{ {f}_{\eufm{fast}}}(\ovl\Theta,\chi,\L,\ovl\vartheta,\k)+\m\,{f}_{\eufm{fast}, \eufm{exp}}(\ovl\Theta,\chi,\L,\ovl\vartheta,\k,\ell)\eeq where ${\rm h}_{\eufm{fast}}$ is as in Proposition \ref{exponential average}, and \beqa{fast split} { {f}_{\eufm{fast}}}&:=&\sum_{i=1}^{n-1}{ {f}^{i}_{\eufm{fast}}},\qquad {f}_{\eufm{fast}, \eufm{exp}}:=\sum_{i=1}^{n-1}{f}_{\eufm{fast}, \eufm{exp}}^{i}. \eeqa Here, \vskip.1in\noi\paragraph{\bf 1.} The ``fast frequency-map'' $$\o_{\eufm{fast}}:=\partial {\rm h}_{\eufm{fast}}$$ is a diffeomorphism of $\cA$ with non-vanishing Jacobian matrix on $\cA_{\hat{\eufm c}\theta}$ and, moreover, $$ \o_\eufm{fast}\in \cD_{\g_{\eufm{fast}}, \t}^{\bar K, \n_{\eufm{fast}}}\qquad \forall\ \L\in \cA,$$ with $$\g_{\eufm{fast}}:=(\g_1,\ \cdots,\ \g_n)\qquad \n_{\eufm{fast}}:=(\n_1,\ \cdots,\ \n_n)$$ and $\n_i$, $\g_i$ as in \equ{nugamma};

 \vskip.1in\noi\paragraph{\bf 2.} the functions ${ {f}^{i}_{\eufm{fast}}}$, ${f}_{\eufm{fast}, \eufm{exp}}^{i}$ do not depend on $\k_0$; the ${ {f}^{i}_{\eufm{fast}}}$'s are given by \begin{equation}\label{averaged bounds*} \begin{split}&{ {f}^{i}_{\eufm{fast}}}={ {f}^{i}_{\eufm{fast}}}({\rm t}_i, {\rm y}_i, {\rm x}_i)=\ovl{{f}_{\cP}^{i}}^{(\ge 2)} ({\rm t}_i, {\rm y}_i, {\rm x}_i)+\widetilde{ {f}^{i}_{\eufm{fast}}}({\rm t}_i, {\rm y}_i, {\rm x}_i) ,\qquad i=1,\ \cdots ,\ n-1,\end{split}
 \end{equation} with 
\begin{align*}
	\ovl{{f}_{\cP}^{i}}^{(\ge 2)} 
	&:= \sum_{j=i+1}^n\ovl{{f}_{\cP}^{ij}}^{(\ge 2)} \\
    {\rm t}_i 
    &:=\big(\Theta_{i}, \ \cdots,\ \Theta_{n-1}, \ \vartheta_{i},\ \cdots,\ \vartheta_{n-1}\big), \\
    {\rm y}_i
    &:=(\chi_{i-1}, \ \cdots, \chi_{n-1},\ \L_i,\ \cdots,\ \L_{n}\big)\ \\
    {\rm x}_i 
    &:=(\k_{i}, \cdots, \k_{n-1}).
   \end{align*}
     In particular, $\widetilde{ {f}^{i}_{\eufm{fast}}}$ do not depend on $\ell_1, \cdots, \ell_n$;

\vskip.1in\noi\paragraph{\bf 3.} finally, $\widetilde{ {f}^{i}_{\eufm{fast}}}$, ${{f}_{\rm exp, \eufm{ fast}}^{i}}$ satisfy the following bounds \beqa{averaged bounds} &&\|\widetilde{ {f}^{i}_{\eufm{fast}}}\|_{ {\eufm D}_{\eufm{fast}}}\le\frac{1}{{\eufm c}_1}\m \bar K\bigg(\frac{a_n^+}{a_1^-}\bigg)^{\frac{3}{2}}\frac{1}{a_{i+1}^-},\qquad \|{f}_{\eufm{fast}, \eufm{exp}}^{i}\|_{ {\eufm D}_{\eufm{fast}}}\le \frac{1}{{\eufm c}_1}\frac{e^{-\widehat{\eufm c}\bar Ks}}{a_{i+1}^-}. \eeqa
\end{proposition}
\nl Let ${\rm L}_0$, $\cdots$, ${\rm L}_n$ be defined as ${\eufm L}_i$ in \equ{Li}, with $\n=m=n$ and $\n_1=\cdots=\n_n=1.$
\begin{lemma}
    \label{non resonance lemma new} If $\bar K$ verifies the inequality in \equ{bar K}, then one can find a number ${\eufm c}_3$, depending only on $m_0$, $\cdots$, $m_n$, such that \beqano |\o_{\rm k,\eufm{fast}}(\L)\cdot k|\ge\frac{{\eufm c}_3}{(a_j^+)^{3/2}}\quad \forall\ k\in {\rm L}_{j-1}\setminus {\rm L}_j,\quad |k|\le \bar K,\quad \forall\ \L\in {\rm A}_{\theta},\quad \forall\ j=1,\cdots, n.\eeqano
\end{lemma}
\begin{proof} For $\L\in {\rm A}_{\theta}$, $\o_{\rm k,\eufm{fast}, j}:=\frac{{\eufm M}_j^2 {\eufm m}_j^3}{\L_j^3}$ verifies $ \frac{\sqrt{{\eufm M}}_j}{(a_j^+)^{3/2}}\le |\o_{\rm k,\eufm{fast}, j}|\le \frac{\sqrt{\eufm M}_j}{(a_j^-)^{3/2}}$. In the case $j=n$, we find $|\o_{\rm k,\eufm{fast}}\cdot k|=|\o_{\rm k,\eufm{fast}, n} k_n|\ge \frac{\sqrt{{\eufm M}}_j}{(a_n^+)^{3/2}}$, since $k_n\ne 0$. Let then $j\ne n$. For $k\in {\rm L}_{j-1}\setminus {\rm L}_j$, $k_j\ne 0$, so, inequality \equ{bar K}, with ${\eufm c}_2\le \frac{\min_j\sqrt{{\eufm M}_j}}{\max_j\sqrt{{\eufm M}_j}}$, and \equ{Choice of parameters} imply $$\bar K \le\frac{\min_j\sqrt{{\eufm M}_j}}{\max_j\sqrt{{\eufm M}_j}}\min_{1\le j\le n-1}\Big(\frac{a_{j+1}^-}{a_j^+}\Big)^{3/2}$$ and hence \beqano |\o_{\rm k,\eufm{fast}}\cdot k|&=&|\sum_{i=j}^n \o_{\rm k,\eufm{fast}, i} k_i|\ge \inf_{\cA_\theta}|\o_{\rm k,\eufm{fast}, j}|-\bar K \max_{j< i\le n}\sup_{\cA_\theta}|\o_{\rm k,\eufm{fast},i}|\nonumber\\
&\ge& \frac{\sqrt{{\eufm M}}_j}{(a_j^+)^{3/2}}-\bar K\frac{\max_{i>j}\sqrt{{\eufm M}_i}}{(a_{j+1}^-)^{3/2}}\ge\frac{\sqrt{{\eufm M}}_j}{2(a_j^+)^{3/2}}eq \eeqano
\end{proof}

\vskip.in\textsc{Proof of  Proposition \ref{Keplerian averages}.} The proof proceeds by recursion, in $n$ steps. We describe the $h^{\rm th}$ step of this recursion, with $h=1$, $\cdots$, $n$. We start with {a Hamiltonian} of the form \beq{start***}{\rm H}_{h-1}={\rm h}_{\eufm{fast}}^0+\m\,{f}_{h-1}\eeq where ${\rm h}_{\eufm{fast}}^0$ is as in \equ{hk0}, and a domain $${\eufm D}_{h-1}={\rm T}_{\Theta^{+(h-1)}, \vartheta^{+(h-1)}}\times{\rm X}_{\theta^{(h-1)}}\times \cA_{\theta^{(h-1)}}\times \torus^n_{s^{(h-1)}}\times \torus^n_{s^{(h-1)}}. $$ When $h=1$, we take ${\rm H}_{0}:={\rm H}_\cP$, $\Theta_+^\ppo:=\Theta^+$, $\vartheta_+^\ppo:=\vartheta^+$, $\theta^\ppo:=\theta$, $s^\ppo:=s$, $ {f}_{0}:={f}_\cP$ and we decompose $${f}_{0}:=\widehat{{f}_{0}}:=\sum_{i=1}^{n-1}\widehat{ f_0^i}\qquad {\rm with}\qquad f_0^i:=\sum_{j=i+1}^n f_{\cP}^{ij}.$$ We observe that $\widehat{ f_0^i}$ depends on the coordinates \beqano &&\Theta_{i},\quad \cdots,\quad \Theta_{n-1},\quad \chi_{i-1},\quad \cdots, \chi_{n-1},\quad\L_i,\quad\cdots,\quad \L_n\nonumber\\
&&\vartheta_i,\quad \cdots,\quad \vartheta_{n-1},\quad \k_i,\quad \cdots,\quad \k_{n-1},\quad \ell_i ,\quad\cdots,\quad \ell_n.\eeqano For $n\ge 3$ and $2\le h\le n-1$, we assume, inductively, that ${f}_{h-1}$ is a sum \beqa{thesis***}{f}_{h-1}=\widehat{f_{h-1}}+f_{\eufm{exp}, h-1}=\sum_{1\le i\le n}\widehat{ {f}_{h-1}^{i}}+\sum_{1\le i\le n}{f}_{\eufm{exp}, h-1}^{i},\eeqa where, in turn, \beqano\widehat{ {f}_{h-1}^{i}}=\ovl{{f}_{h-1}^{i}}+\widetilde{{f}_{h-1}^{i}}\eeqano with $\ovl{{f}_{h-1}^{i}}$, $\widetilde{{f}_{h-1}^{i}}$ depending only on the coordinates \beqano &&\Theta_{i},\quad \cdots,\quad \Theta_{n-1},\quad \chi_{i-1},\quad \cdots, \chi_{n-1},\quad\L_i,\quad\cdots,\quad \L_n\nonumber\\
&&\vartheta_i,\quad \cdots,\quad \vartheta_{n-1},\quad \k_i,\quad \cdots,\quad \k_{n-1},\quad \ell_{i\vee h} ,\quad\cdots,\quad \ell_n\eeqano and $\ovl{{f}_{h-1}^{i}}$, $\widetilde{{f}_{h-1}^{i}}$, $f_{\eufm{exp}, h-1}$ verifying the following bounds and identities \beqa{thesis} \ovl{{f}_{h-1}^{i}}&=&\P_{{\rm L}_{h-1}}T_{\bar K}\widehat{{f}_{h-2}^{i}}\nonumber\\
\|\widetilde{{f}_{h-1}^{i}}\|_{{\eufm D}_{h-1}}&\le&{\eufm C}_{1, h-1}{\m \bar K}\big(\frac{a_n^+}{a_1^-}\big)^{\frac{3}{2}}\|\widehat{ {f}_{h-2}^{i}}\|_{{\eufm D}_{h-2}}\nonumber\\
\| {f}_{\eufm{exp}, h-1}^{i}\|_{{\eufm D}_{h-1}}&\le&{\eufm C}_{2, h-1}e^{-Ks^{(h)}}\|\widehat{ {f}_{h-2}^{i}}\|_{{\eufm D}_{h-2}}. \eeqa Here $\P_{{\rm L}_h}$ denotes the projection over the module ${\rm L}_h$. In any case, $h=1$, or $2\le h\le n-1$, we focus on the Hamiltonian \beq{start}\widehat{{\rm H}_{h-1}}={\rm h}_{\eufm{fast}}^0+\m\,\widehat{{f}_{h-1}}={\rm h}_{\eufm{fast}}^0+\m\,\sum_{i=1}^{n-1}\widehat{ {f}_{h-1}^{i}}.\eeq Our purpose is to apply Proposition \ref{iterative lemma} to this Hamiltonian, in the case that the abstract system \equ{abstract system} does not depend on the coordinates $(p,q)$. To this end, we take the coordinates $$I:=\L,\quad \f:=\ell,\quad \eta:=(\ovl\Theta,\chi),\quad \xi:=(\ovl\vartheta,\k),$$ the functions $f_i$ in \equ{f***} to be the $\widehat{{f}_{h-1}^{n-i}}$, and \beqano &&{N=n-1,\qquad \n=n,\qquad m_{i}:=2i}\nonumber\\
&&{(I_1,\cdots, I_{\n}):=(\L_n,\cdots,\L_{1})}\nonumber\\
&&(\f_1,\cdots, \f_{\n_{i}}):=(\ell_n,\cdots,\ell_{\max\{n-i,h\}})\nonumber\\
&& (\eta_1,\cdots, \eta_{m_{i}}):=(\Theta_{n-1}, \cdots, \Theta_{n-i}, \chi_{n-1}, \cdots, \chi_{n-i-1})\nonumber\\
&& (\xi_1,\cdots, \xi_{m_{i}}):=(\vartheta_{n-1}, \cdots, \vartheta_{n-i}, \k_{n-1}, \cdots, \k_{n-i})\nonumber\\
&&u_i:=(\L_n,\cdots, \L_1, \Theta_{n-1}, \cdots, \Theta_{n-i}, \chi_{n-1}, \cdots, \chi_{n-i-1},\vartheta_{n-1}, \cdots, \vartheta_{n-i}, \k_{n-1}, \cdots, \k_{n-i}). \eeqano The non-resonance assumption \equ{non res} for $\o=\o_{\rm k, fast}=\partial_\L{\rm h}_{\rm k, fast}$, with $${\eufm Z}_{i}={\rm L}_{h-1},\qquad {\eufm Z}=\cup_{i}{\eufm Z}_{i}={\rm L}_{h-1},\qquad {\eufm L}={\rm L}_h\qquad K=\bar K$$ is ensured by Lemma \ref{non resonance lemma new}, with $${\eufm a}=\frac{{\eufm c}_3}{(a_h^+)^{3/2}},\quad A=\cA,\quad r=\theta_1^{(h-1)}.$$ Now we have to check condition \equ{new smallness cond}. In the case $2\le h\le n-1$ the inductive assumptions \equ{thesis} and assumption \equ{fast averaging smallness assumption} imply
\begin{align}
\|\widehat{ {f}_{h-1}^{i}}\|_{{\eufm D}_{h-1}} &\le \|\ovl{{f}_{h-1}^{i}}\|_{{\eufm D}_{h-1}}+\|\widetilde{{f}_{h-1}^{i}}\|_{{\eufm D}_{h-1}}\le \Big(1+{\eufm C}_1{\m \bar K}\big(\tfrac{a_n^+}{a_1^-}\big)^{\frac{3}{2}}\Big)\|\widehat{ {f}_{h-2}^{i}}\|_{{\eufm D}_{h-2}}\notag\\
&\le\cdots\le (1+{\eufm C}_{1, h-1}{\eufm c}_1)^{h-1}\|\widehat{ {f}_{0}^{i}}\|_{{\eufm D}_{0}}\le\frac{{\eufm C}_{4, h-1}}{a_i^-}=:E_i. \label{thesis******}
\end{align}
An {analogous} bound holds also for $h=1$. The numbers $c_i$ and $d_i$ in \equ{c and d} may be evaluated as $$c_{i}=e(1+2ie)/2\qquad d_{i}=\min\{ \theta_1^{(h-1)}s^{(h-1)}, \Theta_i^{+(h-1)}\vartheta_i^{+(h-1)} \}={\eufm c}_2\theta_1^{(h-1)}.$$ From these bounds it is immediate to see that inequality \equ{new smallness cond} is implied by \equ{fast averaging smallness assumption}, provided ${\eufm c}_1<2^{-7}\frac{6}{7}(\frac{8}{9})^{n-2}{\eufm c}_2/({\eufm C}_4c_n)$. Then Proposition \ref{iterative lemma} applies. Its thesis implies that $\widehat{{\rm H}_{h-1}}$ in \equ{start} can be conjugated to a suitable ${\rm H}^*_{h}={\rm h}_{\rm k, \eufm{ fast}}+\m f^*_{h}$, where $f^*_h$ verifies equalities and inequalities in \equ{thesis***}-\equ{thesis} with $h$ replaced by $h+1$ and ${\eufm C}_{1, h-1}$, ${\eufm C}_{2, h-1}$ replaced by suitable ${\eufm C}^*_{1, h}$, ${\eufm C}^*_{2, h}$. Then, applying the same transformation to ${{\rm H}_{h-1}}$ in \equ{start***}, we shall conjugate ${\rm H}_{h-1}$ to ${\rm H}_{h}={\rm h}_{\rm k, \eufm{ fast}}+\m f_{h}$, where $f_h$ satisfies the same equalities and inequalities {as} $f^*_h$, with suitable ${\eufm C}_{1, h}\ge {\eufm C}^*_{1, h}$, ${\eufm C}_{2, h}\ge {\eufm C}^* _{2, h}$.

\nl After we have performed $n$ steps, we let ${\eufm D}_{\eufm{fast}}:={\eufm D}_n$, $ {\rm H}_{\eufm{fast}, \eufm{exp}}:={\rm H}_{n}$, ${f_{\eufm{fast}}^i}:=\widehat{f_n^i}$, $\widetilde{f_{\eufm{fast}}^i}:={f_{\eufm{fast}}^i}-\ovl{f_\cP^i}$, $f_{\eufm{fast}, \eufm{exp}}^i:=f_{\rm exp, {n}}^i$, $\widehat f_{\eufm{fast}}:=\sum_{i=1}^{n-1}\widehat{f^i}$, $\widetilde f_{\eufm{fast}}:=\sum_{i=1}^{n-1}\widetilde{f^i}$, $f_{\eufm{fast}, \eufm{exp}}:=\sum_{i=1}^{n-1}f_{\eufm{fast}, \eufm{exp}}^i$ , with $\ovl{f_\cP^i}:=\sum_{j=i+1}^n \ovl{f_\cP^{ij}}$. Therefore, $$ {\rm H}_{\eufm{fast}}={\rm h}_{\eufm{fast}}^\ppo+\m\big(\widehat f_{\eufm{fast}}+f_{\eufm{exp}, \eufm{fast}}\big)={\rm h}_{\eufm{fast}}^\ppo+\m\big(\sum_{1\le i<j\le n}\ovl{f_\cP^{ij}}+ \widetilde f_{\eufm{fast}}+f_{\rm exp,\eufm{fast}}\big)$$ reduces to \equ{fij fast averaged} and the formulae given below, using \equ{split barfij}.

\nl It remains to check the bound on the left in \equ{averaged bounds} (the one on the right follows by construction). This follows by telescopic arguments. Indeed,

\beqano \|\widetilde{f_{\eufm{fast}}^i}\|_{{\eufm D}_n}&=&\|{f_{\eufm{fast}}^i}-\ovl{f_\cP^i}\|_{{\eufm D}_n}=\|\widehat{f_n^i}-\ovl{f_\cP^i}\|_{{\eufm D}_n}=\|\P_{ {\rm L}_n}\widehat{f_n^i}-\P_{ {\rm L}_n}{f_\cP^i}\|_{{\eufm D}_n} \nonumber\\
&\le&\sum_{h=1}^n \|\P_{ {\rm L}_n}\widehat{f^i_h}-\P_{ {\rm L}_n}T_{\bar K} \widehat{f_{h-1}^i}\|_{{\eufm D}_n}\nonumber\\
&=&\sum_{h=1}^n \|\P_{ {\rm L}_n}\widehat{f^i_h}-\P_{ {\rm L}_n}\P_{ {\rm L}_h}T_{\bar K} \widehat{f_{h-1}^i}\|_{{\eufm D}_n}\nonumber\\
&\le&\sum_{h=1}^n \|\widehat{f^i_h}-\P_{ {\rm L}_h}T_{\bar K} \widehat{f_{h-1}^i}\|_{{\eufm D}_n}\nonumber\\
&\le&\sum_{h=1}^n \|\widehat{f^i_h}-\P_{ {\rm L}_h}T_{\bar K} \widehat{f_{h-1}^i}\|_{{\eufm D}_h}\nonumber\\
&\le&{\m \bar K}\big(\frac{a_n^+}{a_1^-}\big)^{\frac{3}{2}}\frac{\sum_{h=1}^n{\eufm C}_{1, h}{\eufm C}_{4, h-1}}{a_{i+1}^-} .\eeqano Here, we have used the second bound in \equ{thesis}, \equ{thesis******}, that $\widehat{f_n^i}$ does not depend on $\ell_1$, $\cdots$, $\ell_n$, and, finally, $\P_{{\rm L}_n}=\P_{{\rm L}_n} T_{\bar K}=\P_{{\rm L}_n} \P_{{\rm L}_h} $, for all $1\le h\le n$. 
\qed

\section{Secular normalizations}\label{Secular normalizations} Consider the following truncation $${{\rm H}_{\eufm{ fast}}}(\ovl\Theta,\chi,\L,\ovl\vartheta,\k):={\rm h}_{\eufm{fast}}(\L)+\m\,{ {f}_{\eufm{fast}}}(\ovl\Theta,\chi,\L,\ovl\vartheta,\k)$$ of the Hamiltonian ${\rm H}_{\eufm{fast}, \eufm{exp}}$ in \equ{fij fast averaged}. The purpose of this section is to describe an iterative scheme which, after $(n-1)$ steps, conjugates ${{\rm H}_{\eufm{ fast}}}$ to a close-to be integrable system, with an arbitrarily small remainder. 

\nl Let us firstly establish the following notation.
\begin{enumerate}
    \item[\tiny\textbullet] Given a Taylor-Fourier expansion of the form \beqano g(p, q, \k)=\sum_{(a,b)\in \natural^{2m_1}\atop k\in \integer^{m_2}}g_{a,b,k} \bigg(\frac{p-{\rm i} q}{\sqrt2}\bigg)^a\bigg(\frac{p+{\rm i} q}{\sqrt2{\rm i}}\bigg)^be^{{\rm i}k\cdot\k}\qquad (p, q, \k)\in B^{2m_1}(0)\times \torus^{m_2}. \eeqano we denote as $$\P_{p,q,\k} g:=\sum_{a \in \natural^{m_1}}g_{0,a,a}\bigg(\frac{p^2+q^2}{2{\rm i}}\bigg)^a. $$
\end{enumerate}
\nl
\begin{proposition}
    \label{claim} There exists number $\ovl{\eufm c}_h$, depending only on $n$, $m_0$, $\cdots$, $m_n$, $a^\pm_n$ such that, for any $h=1$, $\cdots$, $n-1$ and any $\bar K$, $\bar\g>0$ such that \equ{small secular} hold with ${\eufm c}$ replaced by $\ovl{\eufm c}_h$, one finds open sets $$B_j^*\subset B^{2}_{\varepsilon_j},\qquad {\rm G}^*_{j}\subset{\rm G}_j:=\Big[{\rm G}_j^+,\ {\rm G}_j^+\Big],\qquad j=n-h,\ \cdots,\ n-1$$ verifying \beq{good set}\meas\big(B^*_j\times{\rm G}_j^*\big)\ge\Big(1-\frac{\bar\g}{\ovl{\eufm c}_h}\Big)\meas(B^{2}_{\varepsilon_j}\times{\rm G}_j\big)\eeq {such that,} defining

\begin{equation}\label{domain Dsec} \begin{split}&{\rm T}_{\ovl{\eufm c}_h\theta}^h :=\left\{
    \begin{array}
        {llll}\dst \Big\{ (\Theta_1,\cdots, \Theta_{n-h-1}, \vartheta_1,\cdots,\vartheta_{n-h-1})\in \complex^{n-1}\times \torus_\complex^{n-1}:&\\
        \qquad\qquad\ \ |\vartheta_j-\p|\le \ovl{\eufm c}_h\frac{\theta_j}{{\rm G}_n^-},\quad |\Theta_{j}|\le \ovl{\eufm c}_h {\rm G}_n^+\ \\
        \qquad\qquad\ \ \forall\ j=1,\cdots, n-h-1\Big\}\ &n\ge 3,\ 1\le h<n-2\\
        \\ \emptyset&{\rm otherwise}
    \end{array}
    \right.\\\\
    &B^{*h}_{\ovl{\eufm c}_h\ovl r} :=(B^*_{n-h})_{\ovl{\eufm c}_h\sqrt{\ovl r_{n-h}}}\times\cdots\times (B^*_{n-1})_{\ovl{\eufm c}_h\sqrt{\ovl r_{n-1}}}\\\\
    &{\rm X}^{*h}_{ \ovl{\eufm c}_h\theta, \ovl{\eufm c}_h\bar r}:=\Big\{ \chi=(\chi_{0},\cdots,\ \chi_{n-1}):\quad \chi_{i-1}-\chi_{i}\in ({\rm G}^*_i)_{\ovl{\eufm c}_h\theta_i},\ \chi_{j-1}-\chi_{j}\in ({\rm G}_j)_{ \ovl{\eufm c}_h\bar r_j }\\
&\qquad\qquad\qquad \forall\ i=1,\cdots, n-h-1,\ j=n-h,\cdots, n,\quad \chi_n:=0\Big\}\\\\
&{\eufm D}^{h}_{\eufm{sec}}:= {\rm T}^h_{\ovl{\eufm c}\theta} \times B^{*h}_{\ovl{\eufm c}\varepsilon}\times{\rm X}^{*h}_{\ovl{\eufm c}\theta, \ovl{\eufm c}\ovl r} \times{\rm A}_{\ovl{\eufm c}_h\widetilde r}\times \torus^n_{\ovl{\eufm c}_h s}\times \torus^n_{\ovl{\eufm c}_h s} \end{split}\end{equation} 
    a real-analytic transformation $$\Phi_{\eufm{sec}, h} :\quad {\eufm D}^{h}_{\eufm{sec}}\to {\eufm D}_{\eufm{fast}},$$ may be found, which conjugates ${f_\eufm{fast}}$ to a new function $$f_{\eufm{sec}, h}:={ {f}_{\eufm{fast}}}\circ\Phi_{\eufm{sec}, h} $$ enjoying the following properties.

\nl \vskip.1in\noi\paragraph{\bf 1.} Denoting by $({\rm t}^\pph, {\rm z}^\pph, {\rm y}^\pph, {\rm x}^\pph)$, where \beqa{new coordinates*} &&{\rm t}^\pph=(\Theta^\pph, \vartheta^\pph)=(\Theta^\pph_1,\ \cdots,\ \Theta^\pph_{n-h-1},\ \vartheta^\pph_1,\ \cdots,\ \vartheta^\pph_{n-h-1})\nonumber\\
    &&{\rm z}^\pph=(p^\pph, q^\pph)=(p^\pph_{n-h},\ \cdots,\ p^\pph_{n-1},\ q^\pph_{n-h},\ \cdots,\ q^\pph_{n-1})\nonumber\\
    &&{\rm y}^\pph=(\chi^\pph, \L^\pph)=(\chi^\pph_0,\ \cdots,\ \chi^\pph_{n-1},\ \L^\pph_1,\ \cdots,\ \L^\pph_n)\nonumber\\
    &&{\rm x}^\pph=(\k^\pph, \ell^\pph)=(\k^\pph_0,\ \cdots, \ \k^\pph_{n-1},\ \ell^\pph_1,\ \cdots,\ \ell^\pph_n), \eeqa coordinates on ${\eufm D}^{h}_{\eufm{sec}}$ then ${\Phi_{\eufm{sec}, h}}$ is co-variant with the symmetry: \beqano &&{\Phi_{\eufm{sec}, h}}(-{\rm t}^\pph, -{\rm z}^\pph, {\rm y}^\pph, {\rm x}^\pph)=(-{\rm t}^\ppo, {\rm y}^\ppo, {\rm x}^\ppo) \qquad {\rm if}\nonumber\\
    \nonumber\\
    &&{\Phi_{\eufm{sec}, h}}({\rm t}^\pph, {\rm z}^\pph, {\rm y}^\pph, {\rm x}^\pph)=({\rm t}^\ppo, {\rm y}^\ppo, {\rm x}^\ppo) \eeqano and hence, $f_{\eufm{sec}, h}$ is even around $${\rm t}^\pph=(0, k\p),\qquad {\rm z}^\pph=0\qquad k\in\{0,\ 1\}^{n-h-1}$$

 \vskip.1in\noi\paragraph{\bf 2.} Defining \beqa{good transformation} &&{\rm t}_i^\pph:=\arr{\big(\Theta^\pph_i,\cdots,\Theta^\pph_{n-h-1}, \vartheta^\pph_i,\cdots,\vartheta^\pph_{n-h-1}\big)\qquad i\le n-h-1\\
    \\
    \emptyset\qquad\qquad{\rm otherwise}} \nonumber\\
    &&\hat{\rm y}_i^\pph=\arr{\Big(\frac{(p^\pph_i)^2+(q^\pph_i)^2}{2},\ \cdots ,\ \frac{(p^\pph_{n-1})^2+(q^\pph_{n-1})^2}{2},\ \chi^\pph_{i-1},\ \cdots,\ \chi_{n-1},\\
    \qquad \L^\pph_{i},\ \cdots,\ \L^\pph_n\Big)\qquad i\ge n-h\\
    \\
    \Big(\frac{(p^\pph_{n-h})^2+(q^\pph_{n-h})^2}{2},\ \cdots ,\ \frac{(p^\pph_{n-1})^2+(q^\pph_{n-1})^2}{2},\ \chi^\pph_{i-1},\ \cdots,\ \chi_{n-1},\\
    \qquad \L^\pph_{i},\ \cdots,\ \L^\pph_n\Big)\qquad {\rm otherwise} } \nonumber\\
    &&\hat{\rm x}_i^\pph=\arr{\big(\k^\pph_{i},\ \cdots,\ \k^\pph_{n-h-2}\big)\quad n\ge4\ \&\ 1\le h\le n-3\ \&\ 1\le i\le n-h-2\\
    \\ \emptyset\qquad\qquad{\rm otherwise}} \nonumber\\
    \eeqa and $\hat{\rm y}:=\hat{\rm y}_1$, $\hat{\rm x}:=\hat{\rm x}_1$, $f_{\eufm{sec}, h}$ has the form \beqa{step h} f_{\eufm{sec}, h}({\rm t}^\pph, {\rm z}^\pph,{\rm y}^\pph,{\rm x}^\pph)&=&{\rm h}_{\eufm{sec},h}(\hat{\rm y}_{n-h}^\pph)+f_{\eufm{norm}, h}({\rm t}^\pph, \hat{\rm y}^\pph, \hat{\rm x}^\pph )\nonumber\\
    &+&f_{\eufm{exp}, \eufm{sec}, h}({\rm t}^\pph, {\rm z}^\pph,{\rm y}^\pph,{\rm x}^\pph) \eeqa with \beqa{3split} &&{\rm h}_{\eufm{sec}}(\hat{\rm y}_{n-h}^\pph)=\sum_{i=n-h}^{n-1}{\rm h}_{\eufm{sec}}^i(\hat{\rm y}_i^\pph)\nonumber\\
    \nonumber\\
    &&f_{\eufm{norm}, h}({\rm t}^\pph, \hat{\rm y}^\pph, \hat{\rm x}^\pph )=\sum_{i=1}^{n-h-1}f_{\eufm{norm}, h}^i({\rm t}_i^\pph, \hat{\rm y}_i^\pph, \hat{\rm x}_i^\pph ) \eeqa where

\vskip.1in\noi\paragraph{\bf 3.} the functions ${\rm h}_{\eufm{sec}}^i$ $f_{\eufm{norm}, h}^i$ may be decomposed as \beqa{split step h} &&{\rm h}_{\eufm{sec}}^i(\hat{\rm y}^\pph_i)=\ovl{{\rm h}_{\eufm{sec}}^i}(\hat{\rm y}^\pph_i)+\widetilde{{\rm h}_{\eufm{sec}}^i}i(\hat{\rm y}^\pph_i)\nonumber\\
    \nonumber\\
    && f_{\eufm{norm}, h}^i({\rm t}_i^\pph, \hat{\rm y}^\pph_i, \hat{\rm x}^\pph_i)=\ovl{f_{\eufm{norm}, h}^i}({\rm t}_i^\pph, \hat{\rm y}^\pph_i, \hat{\rm x}^\pph_i)+\widetilde{f_{\eufm{norm}, h}^i}({\rm t}_i^\pph, \hat{\rm y}^\pph_i, \hat{\rm x}^\pph_i) \eeqa where \beq{zero approximation}\ovl{f_{\eufm{norm}, h}^i}=\sum_{j=i+1}^n\P_h\big(\ovl{f_\cP^{ij}}^{(\ge 2)}\circ\ovl{\phi_{\eufm{int}}^{n-1}}\circ\cdots\circ\ovl{\phi_{\eufm{int}}^{n-h}}\big)\eeq and $\ovl{{\rm h}_{\eufm{ sec}}^i}$, $\ovl{\phi_{\eufm{int}}^{i}}$ as in Lemma \ref{integration***}. The functions $\widetilde{{\rm h}_{\eufm{sec}, h}}$, $\widetilde{f_{\eufm{norm}, h}}$, $f_{\eufm{exp}, \eufm{sec}, h}$ in \equ{step h} may be bounded as \beqa{h-1 bound} && |\widetilde{{\rm h}_{\eufm{sec}, h}^i}|\le \frac{1}{{\eufm c}_h}\max\Big\{\m\bar K\big(\frac{a_n}{a_1}\big)^{3/2}\frac{1}{a_{i+1}^-} , \quad {\frac{\bar K^{\bar\t+1}\sqrt\a}{\bar\g}} \frac{(a_i^+)^2}{(a_{i+1}^+)^3},\quad \frac{\varepsilon_{i+1}^2}{\theta_{i+1}} \frac{(a_i^+)^2}{(a_{i+1}^-)^3}\Big\}\nonumber\\
    && |\widetilde{f_{\eufm{norm}, h}^i}| \le \frac{1}{{\eufm c}_h}\max\Big\{\m\bar K\big(\frac{a_n}{a_1}\big)^{3/2}\frac{1}{a_{i+1}^-} , \quad {\frac{\bar K^{\bar\t+1}\sqrt\a}{\bar\g}}\frac{(a_i^+)^2}{(a_{i+1}^-)^3}\Big\} \nonumber\\
    &&|f_{\eufm{exp}, \eufm{sec}, h}|\le \frac{1}{{\eufm c}_h}\frac{(a_{n-1}^+)^2}{(a_n^-)^3}e^{-{\eufm c}_h\bar K}\eeqa

 \vskip.1in\noi\paragraph{\bf 4.} Defining $$\zeta^\pph:=\Big(\frac{(p^\pph_{n-h})^2+(q^\pph_{n-h})^2}{2},\ \cdots ,\ \frac{(p^\pph_{n-1})^2+(q^\pph_{n-1})^2}{2},\ \chi^\pph_{i-1},\ \cdots,\ \chi_{n-1}\Big)$$ so that $$\hat{\rm y}_{n-h}^\pph=(\zeta^\pph, \L^\pph_{n-h},\ \cdots,\ \L^\pph_n) $$ for any $\L^\pph_{n-h}$, $\cdots$, $\L^\pph_n$, the map $$\zeta^\pph\to \o_{\eufm{sec}, h}:=\partial_{\zeta^\pph}{\rm h}_{\eufm{sec}, h}(\zeta^\pph, \L^\pph)$$ is a diffeomorphism of $D_{{\ovl r}}\times {\rm X}_{\ovl r}$, with non-vanishing Jacobian matrix. The set $D^*_{{\ovl r}}\times {\rm X}^*_{\ovl r}$ consists of the subset of $D_{{\ovl r}}\times {\rm X}_{\ovl r}$ such that $\o_{\eufm{fast}, \eufm{sec}}\in \cD^{\bar K,\n_{\eufm{sec}}}_{\g_{\eufm{ sec}};\ovl\t}$, where, if $\n_j$, $\g_j$ are as in \equ{nugamma}, $$\n_{\eufm{sec}}:=(\n_{n+1},\ \cdots,\ \n_{2n-1})\qquad \g_{\eufm{sec}}:=(\g_{n+1},\cdots,\ \g_{2n-1}).$$
\end{proposition}
\vskip.1in \noi We shall give the complete details of the proof of Proposition \ref{claim} along the following sections \ref{Local coordinates}-\ref{first normalization}. In this section we just provide main ideas. \vskip.1in\noi\textsc{
 Scheme of Proof.} The proof is by recursion. The $h^{\rm th}$ step of this recursion starts with \beqano f_{\eufm{sec}, h-1}={\rm h}_{\eufm{sec},h-1}+f_{\eufm{norm}, h-1}+f_{\eufm{exp}, \eufm{sec}, h-1}, \eeqano where, for $h=1$ \beq{first step}{\rm h}_{\eufm{sec},0}\equiv 0\\
,\qquad f_{\eufm{exp}, \eufm{sec}, 0}\equiv 0 ,\qquad f_{\eufm{sec}, 0}:=f_{\eufm{norm}, 0}:={f_{\eufm{ fast}}}, \eeq while, for $n\ge 3$ and $h=2$, $\cdots$, $n-1$, we assume, inductively, that ${\rm h}_{\eufm{sec},h-1}$, $f_{\eufm{sec}, h-1}$ and $f_{\eufm{exp}, \eufm{sec}, h-1}$ satisfy the theses of Proposition \ref{claim}, with $h$ replaced by $(h-1)$. 

\nl The transformation $\phi_{\eufm{sec}}^{n-h}$ conjugating $f_{\eufm{sec}, h-1}$ to $f_{\eufm{sec}, h}$ will be constructed as a product {$\phi_{\eufm{sec}}^{n-h}=\phi_{\eufm{int}}^{n-h}\circ\phi_{\eufm{norm}}^{n-h}$} of an ``integrating'' and a ``normalizing'' transformation.

\nl Due to the bound on $f_{\eufm{exp}, \eufm{sec}, h-1}$, it is enough to focus on the truncation \begin{equation}\label{beginning}\begin{split} \widehat{f_{\eufm{sec}, h-1}}:={\rm h}_{\eufm{sec},h-1}+f_{\eufm{norm}, h-1}={\rm h}_{\eufm{sec},h-1}+\sum_{i=1}^{n-h}f_{\eufm{norm}, h-1}^i({\rm t}_i^{(h-1)}, \hat{\rm y}_i^{(h-1)}, \hat{\rm x}_i^{(h-1)})\end{split}\end{equation} of $f_{\eufm{ sec}, h-1}$. We split $$f_{\eufm{norm}, h-1}=f_{\eufm{norm}, h-1}^{n-h}({\rm t}_{n-h}^{(h-1)}, \hat{\rm y}_{n-h}^{(h-1)}, \hat{\rm x}_{n-h}^{(h-1)})+\sum_{i=1}^{n-h-1}f_{\eufm{norm}, h-1}^i({\rm t}_i^{(h-1)}, \hat{\rm y}_i^{(h-1)}, \hat{\rm x}_i^{(h-1)})$$ and we distinguish two cases.

\vskip.1in \noi {\sl Case $n\ge 3$, $h=2$, $\cdots$, $n-1$.} By the inductive assumption (see \equ{good transformation} with $h$ replaced by $(h-1)$), the function $f_{\eufm{norm}, h-1}^{n-h}$ depends only on $${\rm t}^{(h-1)}_{n-h}=\big(\Theta^{(h-1)}_{n-h},\ \vartheta^{(h-1)}_{n-h}\big)\quad {\rm and}\quad \hat{\rm y}^{(h-1)}_{n-h}$$ therefore is integrable. In Section \ref{local coordinates***}, we shall construct a canonical, real-analytic change of coordinates 
 \beqano {\phi_{\eufm{int}}^{n-h}}:\ &&\eufm{D}_{\eufm{int}}^h\to \eufm{D}_{\eufm{sec}}^{h-1}\nonumber\\
&& ({\rm t}^\pph_*, {\rm z}^\pph_*,{\rm y}^\pph_*,{\rm x}^\pph_* )\to({\rm t}^{(h-1)}, {\rm z}^{(h-1)},{\rm y}^{(h-1)},{\rm x}^{(h-1)} ) \eeqano \beqa{domain Dsec***} {\eufm D}^{h}_{\eufm{int}}:= {\rm T}^h_{\hat{\eufm c}_h\theta}\times B^2_{\hat{\eufm c}_h\varepsilon_{n-h}} \times B^{*,h-1}_{\hat{\eufm c}_h\varepsilon}\times{\rm X}^{*,h-1}_{\hat{\eufm c}_h\theta, \hat{\eufm c}_h\ovl r} \times{\rm A}_{\hat{\eufm c}_h\widetilde r}\times \torus^n_{\hat{\eufm c}_h s}\times \torus^n_{\hat{\eufm c}_h s} \eeqa such that \beq{integration h} f_{\eufm{norm}, h-1}^{n-h}\circ\phi_{\eufm{ int}}^{n-h}={\rm h}_{\eufm{sec}}^{n-h}(\hat{\rm y}^\pph_{*, n-h}) \eeq depends only on $\hat{\rm y}^\pph_{*,n-h}$, where $\hat{\rm y}^\pph_{*,i}$ is defined analogously to $\hat{\rm y}^\pph_i$ in \equ{good transformation}. Here, \beqano \arr{ {\rm t}^\pph_*:=\big(\Theta^\pph_*,\ \vartheta^\pph_*\big)\\
{\rm z}^\pph_*:=\big(p^\pph_*,\ q^\pph_*\big)\\
{\rm y}^\pph_*:=\big(\chi^\pph_*,\ \L^\pph_*\big)\\
{\rm x}^\pph_*:=\big(\k^\pph_*,\ \ell^\pph_*\big)}\qquad \qquad \arr{ {\rm t}^{(h-1)}:=\big(\Theta^{(h-1)},\ \vartheta^{(h-1)}\big)\\
{\rm z}^{(h-1)}:=\big(p^{(h-1)},\ q^{(h-1)}\big)\\
{\rm y}^{(h-1)}:=\big(\chi^{(h-1)},\ \L^{(h-1)}\big)\\
{\rm x}^{(h-1)}:=\big(\k^{(h-1)},\ \ell^{(h-1)}\big)} \eeqano are defined analogously to \equ{new coordinates*}.

\nl We shall construct ${\phi_{\eufm{int}}^{n-h}}$ {in such a way that} it involves only the coordinates \beqano\phi_{\eufm{int}}^{n-h}:\ ({\rm z}_{*, n-h}^\pph, \ {\rm y}_{*, n-h}^\pph,\ {\rm x}_{*,n-h}^\pph)\to ({\rm t}_{n-h}^{(h-1)},\ {\rm z}_{n-h+1}^{(h-1)}, \ {\rm y}_{n-h}^{(h-1)},\ {\rm x}_{n-h}^{(h-1)})\eeqano

\nl with \beqa{true dependence} &&{\rm z}^{(h)}_{*,n-h}:= \big(p^\pph_{*,n-h},\ \cdots,\ p^\pph_{n-1},\ q^\pph_{*,n-h}, \ \cdots,\ q^\pph_{n-1}\big) \nonumber\\
&&{\rm y}_{*,n-h}^\pph:=\big(\chi^\pph_{*,n-h-1},\ \cdots,\ \chi_{*,n-1},\ \L^\pph_{*,n-h},\ \cdots,\ \L^\pph_n\big)\nonumber\\
&&{\rm x}_{*,n-h}^\pph:=\big(\k^\pph_{*,n-h-1},\ \cdots,\ \k^\pph_{*,n-1},\ \ell^\pph_{*,n-h},\ \cdots,\ \ell^\pph_n\big)\nonumber\\
\nonumber\\
&&{\rm t}_{n-h}^{(h-1)}:=\big(\Theta^{(h-1)}_{n-h},\ \vartheta^{(h-1)}_{n-h}\big)\nonumber\\
&&{\rm z}_{n-h+1}^{(h-1)}:=\big(p^{(h-1)}_{n-h+1},\ \cdots,\ p^{(h-1)}_{n-1},\ q^{(h-1)}_{n-h+1}, \ \cdots,\ q^{(h-1)}_{n-1}\big) \nonumber\\
&&{\rm y}_{n-h}^{(h-1)}:=\big(\chi^{(h-1)}_{n-h-1},\ \cdots,\ \chi_{n-1},\ \L^{(h-1)}_{n-h},\ \cdots,\ \L^{(h-1)}_n\big)\nonumber\\
&&{\rm x}_{n-h}^{(h-1)}:=\big(\k^{(h-1)}_{n-h-1},\ \cdots,\ \k^{(h-1)}_{n-1},\ \ell^{(h-1)}_{n-h},\ \cdots,\ \ell^{(h-1)}_n\big)\eeqa and has the form \beqa{integration} \phi_{\eufm{int}}^{n-h}:\quad \left\{
\begin{array}
    {lll} \dst\Theta^{(h-1)}_{n-h}={\rm F}_{\eufm{int}}^\pph(p^\pph_{*,n-h},\ q^\pph_{*,n-h},\ \tilde{\rm y}^{(h)}_*) \\
    \\ \dst\vartheta^{(h-1)}_{n-h}-\p={\rm G}_{\eufm{int}}^\pph(p^\pph_{*,n-h},\ q^\pph_{*,n-h},\ \tilde{\rm y}^{(h)}_*) \\
    \\ \dst \hat {\rm z}_j^{(h-1)}=\hat {\rm z}_{*,j}^{(h)}e^{{\rm i}\psi_{\eufm{int},j}^\pph(p^\pph_{*,n-h},\ q^\pph_{*,n-h},\ \tilde{\rm y}^{(h)}_*)}\\
    \\ \dst{\rm y}_{n-h}^{(h-1)}={\rm y}^\pph_{*,n-h}\\
    \\ \dst{\rm x}_{n-h}^{(h-1)}={\rm x}^\pph_{*,n-h}+\f_{\eufm{int}}^\pph(p^\pph_{*,n-h},\ q^\pph_{*,n-h},\ \tilde{\rm y}^{(h)}_*)
\end{array}
\right. \eeqa with ${\rm F}_{\eufm{int}}^\pph$, ${\rm G}_{\eufm{int}}^\pph $ odd, $\psi_{\eufm{int},j}^\pph$, $\f_{\eufm{int}}^\pph $ even in $(p^\pph_{*,n-h},\ q^\pph_{*,n-h})$, \beqa{dependence on variables} &&\tilde{\rm y}_*^{(h)}:=\big(\frac{(p^\pph_{*,n-h+1})^2+(q^\pph_{*,n-h+1})^2}{2},\cdots \frac{(p^\pph_{*,1})^2+(q^\pph_{*,1})^2}{2},\ { \rm y}^\pph_{*, n-h}\big)\nonumber\\
&&\hat{\rm z}^{(h-1)}_j:=\big(p^{(h-1)}_{j},\ q^{(h-1)}_{j}\big):=p^{(h-1)}_{j}+{\rm i}q^{(h-1)}_{j}\nonumber\\
&& \hat{\rm z}^{(h)}_{*,j}:=\big(p^{(h)}_{*,j},\ q^{(h)}_{*,j}\big):=p^{(h)}_{*,j}+{\rm i}q^{(h)}_{*,j} \eeqa with $j=n-h+1$, $\cdots$, $n-1$, for $ n\ge 3$, $h\ge 2$ and ${\rm y}^\pph_{*, n-h}$ as in \equ{true dependence}.

\nl In particular, observe that $\phi_{\eufm{int}}^{n-h}$ enjoys the following properties:
\begin{enumerate}
    \item[\tiny\textbullet] it is co-variant with the symmetry: if \beqano \phi_{\eufm{int}}^{n-h}({\rm t}_*^\pph, {\rm z}_*^\pph, \ {\rm y}_*^\pph,\ {\rm x}_*^\pph)=({\rm t}^{(h-1)},\ {\rm z}^{(h-1)}, \ {\rm y}^{(h-1)},\ {\rm x}^{(h-1)}),\eeqano then \beqano \phi_{\eufm{int}}^{n-h}(-{\rm t}_*^\pph, -{\rm z}_*^\pph, \ {\rm y}_*^\pph,\ {\rm x}_*^\pph)=(-{\rm t}^{(h-1)},\ -{\rm z}^{(h-1)}, \ {\rm y}^{(h-1)},\ {\rm x}^{(h-1)});\eeqano \item[\tiny\textbullet] leaves the ``actions'' $$\tilde{\rm y}^\pph_*=\tilde{\rm y}^{(h-1)}$$ unvaried, where $\tilde{\rm y}^\pph_*$ is as in \equ{dependence on variables}, and \beqano \tilde{\rm y}^{(h-1)}&:=&\big(\frac{(p^{(h-1)}_{n-h+1})^2+(q^{(h-1)}_{n-h+1})^2}{2},\cdots \frac{(p^{(h-1)}_{1})^2+(q^{(h-1)}_{1})^2}{2},\ {\rm y}_{n-h}^{(h-1)}\big)\ \eeqano is defined analogously;

\item[\tiny\textbullet] leaves the averages with respect to the ${\rm x}$-coordinates unvaried. Namely, for any real-analytic function $g$ on ${\eufm D}_{\eufm{sec}}^{h-1}$, $$\P_{{\rm x}_*^\pph}\big(g\circ\phi_{\eufm{int}}^{n-h}\big)=\big(\P_{{\rm x}^{(h-1)}}g\big) \circ\phi_{\eufm{int}}^{n-h}.$$
\end{enumerate}
\nl Applying ${\phi_{\eufm{int}}^{n-h}}$ to $\widehat{f_{\eufm{sec}, h-1}}$ in \equ{beginning}, we obtain \beqano {f_{\eufm{sec}, \eufm{int}, h-1}}&:=&\widehat{f_{\eufm{sec}, h-1}}\circ{\phi_{\eufm{int}}^{n-h}}={\rm h}_{\eufm{sec}, h-1}+{{\rm h}_{\eufm{sec}}^{n-h}}+\sum_{i=1}^{n-h-1}f_{\eufm{norm},\eufm{int}, h-1}^i ({\rm t}_{*, i}^\pph,\tilde{\hat{\rm y}}_{*, i}^\pph,\tilde{\hat{\rm x}}_{*, i}^\pph)\nonumber\\
&=& \sum_{i=n-h}^{n-1}{\rm h}_{\eufm{sec}, h}^i(\hat{\rm y}^\pph_{*,i}) +\sum_{i=1}^{n-h-1}f_{\eufm{norm},\eufm{int}, h-1}^i ({\rm t}_{*, i}^\pph,\tilde{\hat{\rm y}}_{*, i}^\pph,\tilde{\hat{\rm x}}_{*, i}^\pph) \eeqano with \beqa{norm int h-1} {\rm h}_{\eufm{sec}, h}&:=&{\rm h}_{\eufm{sec}, h-1}+{{\rm h}_{\eufm{sec}}^{n-h}},\qquad f_{\eufm{norm},\eufm{int}, h-1}^i:=f_{\eufm{norm}, h-1}^i\circ\phi^{n-h}_{\eufm{ int}}
\eeqa and (as it follows from \equ{good transformation} with $h-1$ replacing $h$ and \equ{integration}) $f_{\eufm{norm},\eufm{int}, h-1}^i$ depends only on the arguments 
\begin{equation}\label{1st step***}
\begin{split} {\rm t}_{*, i}^\pph:= &\big(\Theta^\pph_{*, i},\cdots,\Theta^\pph_{*, n-h-1}, \vartheta^\pph_{*, i},\cdots,\vartheta^\pph_{*, n-h-1} \big) \\
\tilde{\hat{\rm y}}_{*, i}^\pph:=& \big(p^\pph_{*, n-h},\ q^\pph_{*, n-h},\ \frac{(p^\pph_{*, n-h+1})^2+(q^\pph_{*, n-h+1})^2}{2},\ \cdots ,\ \frac{(p^\pph_{*, n-1})^2+(q^\pph_{*, n-1})^2}{2},\\
&\chi^\pph_{*, i-1},\ \cdots,\ \chi_{*, n-1},\L^\pph_{*, i},\ \cdots,\ \L^\pph_{*, n}\big) \\		
\tilde{\hat{\rm x}}_{*, i}^\pph:=&\arr{\big(\k^\pph_{*, i},\ \cdots,\ \k^\pph_{*, n-h-1}\big)\quad n\ge4\ \&\ 1\le h-1\le n-3 \\
\\ \emptyset\qquad\qquad{\rm otherwise}} \end{split}\end{equation} The next step will be to retain the dependence on $(p^\pph_{n-h}, q^\pph_{n-h})$ only via $\frac{(p^\pph_{n-h})^2+(q^\pph_{n-h})^2}{2}$ and, for $h<n-1$, to to eliminate from ${f_{\eufm{sec}, \eufm{int}, h-1}}$ the dependence upon the angle $k^\pph_{*, n-h-1}$, up to an exponential remainder. Namely, we look for another canonical, real-analytic change of coordinates \beqa{normalization} {\phi_{\eufm{norm}}^{n-h}}:\ &&\eufm{D}_{\eufm{sec}}^h\to \eufm{D}_{\eufm{int}}^h\nonumber\\
&& ({\rm t}^\pph,\ {\rm z}^\pph,\ {\rm y}^\pph,\ {\rm x}^\pph )\to({\rm t}^\pph_*,\ {\rm z}^\pph_*,\ {\rm y}^\pph_*,\ {\rm x}^\pph_* )\eeqa so as to conjugate ${{f_{\eufm{sec},\eufm{int}, h-1}}}$ to a new Hamiltonian
 \begin{equation}\label{normalization***}
 \begin{split}\widehat{f_{\eufm{sec}, h}}:={{f_{\eufm{sec},\eufm{int}, h-1}}}\circ\phi_{\eufm{norm}}^{n-h}={\rm h}_{\eufm{sec}, h}+\sum_{i=1}^{n-h-1}f_{\eufm{norm}, h}^i ({\rm t}_{ i}^\pph,\hat{\rm y}_{i}^\pph,\hat{\rm x}_{i}^\pph)+\widehat{f_{\eufm{ exp}, \eufm{sec}, h}}\end{split} \end{equation} where $f_{\eufm{norm}, h}^i$ and $\widehat{f_{\eufm{ exp}, \eufm{sec}, h}}$ satisfy \equ{split step h}-\equ{h-1 bound}. We choose $\eufm{D}_{\eufm{sec}}^h$ as the subset of $\eufm{D}_{\eufm{int}}^h$ where the map $$\omega_{\eufm{ sec}, h}:=\arr{\dst\partial_{\frac{(p^\pph_{n-h})^2+(q^\pph_{n-h})^2}{2}, \chi^\pph_{n-h-1}}{\rm h}_{\eufm{sec}, h} \quad h=2,\ \cdots, n-2\ \&\ n\ge 4\\
\\
\dst\partial_{\frac{(p^{(n-1)}_{1})^2+(q^{(n-1)}_{1})^2}{2}}{\rm h}_{\eufm{sec}, n-1} \qquad h=n-1 }$$ does not {verify} resonances up to order $\bar K$, and next we apply a suitable normal form theory (Proposition \ref{iterative lemma}). We shall choose $\phi^{n-h}_{\eufm{ norm}}$ {in such a way that}
\begin{enumerate}
    \label{item} \item[\tiny\textbullet] it is co-variant with the symmetry: if \beqano \phi_{\eufm{norm}}^{n-h}({\rm t}^{(h)},\ {\rm z}^{(h)}, \ {\rm y}^{(h)},\ {\rm x}^{(h)})=({\rm t}_*^\pph, {\rm z}_*^\pph, \ {\rm y}_*^\pph,\ {\rm x}_*^\pph),\eeqano then \beqa{item1} \phi_{\eufm{norm}}^{n-h}(-{\rm t}^{(h)},\ -{\rm z}^{(h)}, \ {\rm y}^{(h)},\ {\rm x}^{(h)})=(-{\rm t}_*^\pph, -{\rm z}_*^\pph, \ {\rm y}_*^\pph,\ {\rm x}_*^\pph),\eeqa \item[\tiny\textbullet] leaves the ``actions'' $${\rm y}^\pph_{*, n-h}={\rm y}^{(h)}_{n-h}$$ unvaried, where \beqa{item2} &&{\rm y}^{(h)}_{n-h}:=\big(\frac{(p^{(h)}_{n-h+1})^2+(q^{(h)}_{n-h+1})^2}{2},\cdots \frac{(p^{(h)}_{1})^2+(q^{(h)}_{1})^2}{2},\ \chi^\pph_{n-h},\ \cdots,\chi^\pph_{n-1},\nonumber\\
    &&\qquad\qquad \L_1^\pph,\ \cdots,\ \L_n^\pph \big);\nonumber\\
    &&{\rm y}^{(h)}_{*,n-h}:=\big(\frac{(p^{(h)}_{*,n-h+1})^2+(q^{(h)}_{*,n-h+1})^2}{2},\cdots \frac{(p^{(h)}_{*,1})^2+(q^{(h)}_{*,1})^2}{2},\ \chi^\pph_{*,n-h},\ \cdots,\chi^\pph_{*,n-1},\nonumber\\
    &&\qquad\qquad \L_{*,1}^\pph,\ \cdots,\ \L_{*,n}^\pph \big); \eeqa

 \item[\tiny\textbullet] verifies \beq{item3}\P_{{\rm z}_{*, n-h+1}^\pph, {\rm x}_{*, n-h+1}^\pph}\big(g\circ\phi_{\eufm{norm}}^{n-h}\big)=\big(\P_{{\rm z}^{(h)}_{n-h+1},{\rm x}^{(h)}_{n-h+1}}g\big) \circ\phi_{\eufm{norm}}^{n-h}.\eeq
\end{enumerate}
The thesis of Proposition \ref{claim} at rank $h$ follows, with $$f_{\eufm{ sec}, h}:=\widehat{f_{\eufm{ sec}, h}}+{f_{\eufm{ exp}, \eufm{sec}, h-1}}\circ\phi_{\eufm{ sec}}^{n-h},\ {f_{\eufm{ exp}, \eufm{sec}, h}}:=\widehat{f_{\eufm{ exp}, \eufm{sec}, h}}+{f_{\eufm{ exp}, \eufm{sec}, h-1}}\circ\phi_{\eufm{ sec}}^{n-h}.$$ \vskip.1in \noi {\sl Case $h=1$.} The proof of this case uses similar ideas as the proof of the case $2\le h\le n-1$ for $n\ge3$. However, due to subtle differences between the two cases (compare, \eg, the inductive assumption on $f^{n-h}_{\eufm{ norm}, h-1}$ in \equ{good transformation} for $h\ge 2$ with Eq. \equ{split1}; the definition of ${\rm h}^{n-h}_{\eufm{sec}}$, ${{\phi}^{n-h}_{\eufm{int}}}$ for $h\ge 2$ in \equ{integration h}, with the definition of ${\rm h}^{n-1}_{\eufm{sec}}$, ${{\phi}^{n-1}_{\eufm{int}}}$ in \equ{integration 1} and \equ{integration 1*}), for sake of precision, we briefly discuss also this case.

\nl Let ${f_{\eufm{sec}, 0}}$ be as in \equ{first step}. In view of \equ{fast split} and \equ{averaged bounds*}, we can split\beq{split1}{f_{\eufm{sec}, 0}}=\ovl{f_{\cP}^{n-1, n}}^\ppd+\ovl{f_{\cP}^{n-1}}^{(\ge 3)}+\widetilde{f_{\eufm{fast}}^{n-1}}+\sum_{i=1}^{n-2}{ {f}^{i}_{\eufm{fast}}} \eeq where $$\ovl{f_{\cP}^{n-1}}^{(\ge 3)}:=\ovl{f_{\cP}^{n-1}}^{(\ge 2)}-\ovl{f_{\cP}^{n-1}}^{(2)}$$ and the summand appears only when $n\ge 3$. As for $\ovl{f_{\cP}^{n-1, n}}^\ppd$, by Lemmata \ref{integration***} (see also Lemma \ref{base integration}), we find a domain $\ovl{{\eufm B}_{n-1}}$ (defined in Eq. \equ{cond on gamma} below), a real-analytic and canonical transformation 
\begin{equation}\label{first integration***} \begin{split}{\phi^{n-1}_{\eufm{int}}}:\ \big({\rm z}^\ppu_{*,n-1}, {\rm y}^\ppu_{*,n-1}, {\rm x}^\ppu_{*,n-1}\big)\in \ovl{{\eufm B}_{n-1}} \to \big({\rm z}_{n-1}^\ppo, {\rm y}_{n-1}^\ppo, {\rm x}_{n-1}^\ppo\big)\in \ovl{{\eufm D}_{n-1}}:={\phi^{n-1}_{\eufm{int}}}\big(\ovl{{\eufm B}_{n-1}}\big)\end{split} \end{equation} of the form \equ{integration}, with $h=1$ (but neglecting the coordinates $\hat{\rm z}^\ppo_j$, $\hat{\rm z}^\ppu_{j, *}$) such that \beq{integration 1}\ovl{f_{\cP}^{n-1}}^{(2)}\circ {\phi^{n-1}_{\eufm{int}}}=\ovl{{\rm h}_{\eufm{sec}}^{n-1}}( \hat{\rm y}_{*,n-1}^\ppu)\eeq depends only on \beq{hat y n-1}\hat{\rm y}^\ppu_{*,n-1}=\Big( \frac{(p^\ppu_{*,n-1})^2+(q^\ppu_{*,n-1})^2}{2},\ \chi^\ppu_{*,n-2},\ \chi^\ppu_{*,n-1},\ \L^\ppu_{*, n-1},\ \L^\ppu_{*, n}\Big). \eeq In \equ{first integration***}, we have let \beqano \left\{
\begin{array}
    {llll} {\rm z}^\ppu_{*,n-1}:=\big(p^\ppu_{*,n-1},\ q^\ppu_{*, n-1}\big)\\
    \\ {\rm y}^\ppu_{*,n-1}:=\big(\chi^\ppu_{*,n-2},\ \chi^\ppu_{*,n-1},\ \L^\ppu_{*,n-1},\ \L^\ppu_{*,n}\big)\\
    \\ {\rm x}^\ppu_{*,n-1}:=\big(\k^\ppu_{*,n-2},\ \k^\ppu_{*,n-1}\big)
\end{array}
\right. \qquad \left\{
\begin{array}
    {llll} {\rm t}_{n-1}^\ppo:=\big(\Theta^\ppo_{n-1},\ \vartheta^\ppo_{n-1}\big)\\
    \\ {\rm y}_{n-1}^\ppo:=\big(\chi^\ppo_{n-2},\ \chi^\ppo_{n-1},\ \L^\ppo_{n-1},\ \L^\ppo_{n}\big)\\
    \\ {\rm x}_{n-1}^\ppo:=\big(\k^\ppo_{n-2},\ \k^\ppo_{n-1}\big).
\end{array}
\right. \eeqano We let \beqano \left\{
\begin{array}
    {llll} {\rm t}^\ppo:=\big(\Theta^\ppo,\ \vartheta^\ppo\big)\\
    \\ {\rm y}^\ppo:=\big(\chi^\ppo, \L^\ppo\big)\\
    \\ {\rm x}^\ppo:=\big(\k^\ppo, \ell^\ppo\big)
\end{array}
\right.\qquad \left\{
\begin{array}
    {llll} {\rm t}^\ppu_*:=\big(\Theta^\ppu_*,\ \vartheta^\ppu_*\big)\\
    {\rm z}^\ppu_*:=\big(p^\ppu_*,\ q^\ppu_*\big)\\
    {\rm y}^\ppu_*:=\big(\chi^\ppu_*, \L_*^\ppu\big)\\
    {\rm x}^\ppu_*:=\big(\k^\ppu_*, \ell_*^\ppu\big)
\end{array}
\right. \eeqano analogously to \equ{new coordinates*}, with $h=0,1$, and then we regard the map in \equ{first integration***} as a map $$ {\phi^{n-1}_{\eufm{int}}}:\ \big({\rm t}^\ppu_*,\ {\rm z}^\ppu_*,\ {\rm y}^\ppu_*,\ {\rm x}^\ppu_*\big)\in{\eufm{D}_{\eufm{int}}^1}\to \big({\rm t}^\ppo, \ {\rm y}^\ppo,\ {\rm x}^\ppo\big)$$ on the set \beqano {\eufm{D}_{\eufm{int}}^1}&:=& \Big\{\big({\rm t}^\ppu_*, \ {\rm z}^\ppu_*,\ {\rm y}^\ppu_*,\ {\rm x}^\ppu_*\big) :\quad \big({\rm z}^\ppu_{*,n-1},\ {\rm y}^\ppu_{*,n-1},\ {\rm x}^\ppu_{*,n-1} \big)\in\ovl{\eufm{B}_{n-1}}\Big\} \eeqano where ${\phi^{n-1}_{\eufm{int}}}$ is defined on the extra-coordinates via the identity. ${\eufm{D}_{\eufm{int}}^1}$ has the form in \equ{domain Dsec***}, with $h=1$. Applying this extension to ${f_{\eufm{sec}, 0}}$ in \equ{split1} we obtain \beqano f_{\eufm{sec}, \eufm{int}, 0}:= {f_{\eufm{sec}, 0}}\circ\ovl{\phi^{n-1}_{\eufm{int}}} =\ovl{{\rm h}_{\eufm{sec}}^{n-1}}( \hat{\rm y}_{*,n-1}^\ppu)+\sum_{i=1}^{n-1}f^i_{\eufm{norm}, \eufm{int}, 0}({\rm t}_{*, i}^\ppu,\tilde{\hat{\rm y}}_{*, i}^\ppu,\tilde{\hat{\rm x}}_{*, i}^\ppu) \eeqano where $$f^{n-1}_{\eufm{norm}, \eufm{int}, 0}:=\big(\ovl{f_{\cP}^{n-1}}^{(\ge 3)}+\widetilde{f_{\eufm{fast}}^{n-1}}\big) \circ\ovl{\phi^{n-1}_{\eufm{int}}} , \qquad f^i_{\eufm{norm}, \eufm{int}, 0}:=\widehat{f^i_{\eufm{ fast}}}\circ\ovl{\phi^{n-1}_{\eufm{int}}}$$ and, as a consequence of \equ{averaged bounds*} and of \equ{integration}, with $h=1$, $f^i_{\eufm{norm}, \eufm{int}, 0}$ depends only on the arguments \beqano &&{\rm t}_{*, i}^\ppu:=\big(\Theta^\ppu_{*, i},\cdots,\Theta^\ppu_{*, n-2}, \vartheta^\ppu_{*, i},\cdots,\vartheta^\ppu_{*, n-2}, \big) \nonumber\\
&&\tilde{\hat{\rm y}}_{*, i}^\ppu:= \big( p^\ppu_{*, n-1},\ q^\ppu_{*, n-1},\ \chi^\ppu_{*, i-1},\ \cdots,\ \chi_{*, n-1},\ \L^\ppu_{*, i},\ \cdots,\ \L^\ppu_{*, n}\big)\nonumber\\
&&\tilde{\hat{\rm x}}_{*, i}^\ppu:=\big(\k^\ppu_{*, i},\ \cdots,\ \k^\ppu_{*, n-1}\big).\eeqano Note, in particular, that $f^{n-1}_{\eufm{norm}, \eufm{int}, 0}$ is a function of \beq{term n-1} ({\rm t}_{*,n-1}, {\rm y}_{*,n-1}, {\rm x}_{*,n-1})=(p_{*n-1}^\ppu,\ q_{*n-1}^\ppu, \chi_{*,n-2},\ \chi_{*,n-1},\ \L_{*,n-1},\ \L_{*,n},\ \k_{*n-1}).\eeq In view of the fact that $\ovl{{\rm h}^{n-1}_{\eufm{sec}}}$ depends on the actions in \equ{hat y n-1}, we aim to eliminate from $f_{\eufm{sec}, \eufm{int}, 0}$ the dependence on the following angles $$\left\{
\begin{array}
    {lll} \dst \k_{*,1}\qquad \qquad\qquad&{\rm if}\ &n=2\\
    \\ \dst\k_{*,n-2},\ \k_{*,n-1}&{\rm if}&n\ge 3
\end{array}
\right. $$ and to retain the dependence on $(p^\ppu_{*, n-1},\ q^\ppu_{*, n-1})$ only via $\frac{(p^\ppu_{*, n-1})^2+(q^\ppu_{*, n-1})^2}{2}$. Then we choose a domain ${\eufm D}_{\eufm{sec}}^1\subset {\eufm D}_{\eufm{int}}^1$ as in \equ{domain Dsec} where the frequency $$\omega_{\eufm{ sec}, 1}:=\arr{ \dst \partial_{\frac{p^\ppu_{*, n-1})^2+(q^\ppu_{*, n-1})^2}{2}, \chi^\ppu_{*, n-1}}\ovl{{\rm h}_{\eufm{sec}}^{n-1}}\qquad\qquad\qquad \qquad \ n=2\\
\\
\dst \partial_{\frac{p^\ppu_{*, n-1})^2+(q^\ppu_{*, n-1})^2}{2}, \chi^\ppu_{*, n-2}, \chi^\ppu_{*, n-1}}\ovl{{\rm h}_{\eufm{sec}}^{n-1}}\qquad \qquad\quad \ \ n\ge 3 } $$ is non-resonant up to the order $\bar K$ and on this domain we construct a real-analytic transformation ${\phi_{\eufm{norm}}^{n-1}}$ as in \equ{normalization} which conjugates $\ovl{{{\rm f}_{ \eufm{sec}, 1}}}$ to a Hamiltonian \beqano f_{\eufm{sec}, 1}:=f_{\eufm{sec}, \eufm{int}, 0}\circ\phi^{n-1}_{\eufm{ norm}}=\ovl{{\rm h}_{\eufm{sec}}^{n-1}}( \hat{\rm y}_{n-1}^\ppu)+\sum_{i=1}^{n-1}f^i_{\eufm{norm}, 1}({\rm t}_{ i}^\ppu,\hat{\rm y}_{ i}^\ppu,\hat{\rm x}_{ i}^\ppu)+f_{\eufm{exp}, \eufm{sec}, 1} \eeqano Now, since (as it follows from \equ{term n-1}), $f^{n-1}_{\eufm{norm}, 1}$ is actually a function of $\hat{\rm y}^\ppu_{n-1}$ only, this step is proved, with \beq{integration 1*}{\rm h}_{\eufm{sec}}^{n-1}(\hat{\rm y}^\ppu_{n-1}):=\ovl{{\rm h}_{\eufm{sec}}^{n-1}}( \hat{\rm y}_{n-1}^\ppu)+f^{n-1}_{\eufm{norm}, 1}(\hat{\rm y}^\ppu_{n-1}). \eeq\qed

\subsection{Construction of ${\phi_{\eufm{int}}^{n-1}}$}\label{Local coordinates}

 The following lemma completes Lemma \ref{integration***}. In particular, it provides the transformation ${\phi_{\eufm{int}}^{n-1}}=\ovl{\phi_{\eufm{int}}^{n-1}}$ in \equ{integration 1}.
\begin{lemma}
    \label{base integration} Let $i=1$, $\cdots$, $n-1$. Let $\cA$, $\rm X$, $\theta$ in \equ{domain***} be chosen {in such a way that}
\beqa{cond on gamma} && \inf_{{\eufm D}_\cP}|g|>0,\quad \sup_{{\eufm D}_\cP}|\arg g|<\frac{\p}{4}\nonumber\\
    && \forall\ g\in \Big\{\chi_{i-1},\quad \chi_i,\quad \chi_{i-1}+\chi_i,\quad{5\chi_{i-1}\L^2_{i}-(\chi_{i-1}-\chi_i)^2(4\chi_{i-1}-\chi_i)}\Big\}.\eeqa Then, the domains $\ovl{{\eufm B}_i}$ in \equ{tr*}, the functions $\ovl{{\rm h}_{\eufm{sec}}^{i}}$ and the transformations $\ovl{\phi_{\eufm{int}}^{i}}$ can be taken as follows
    
    \beqa{Nn} \ovl{{\eufm B}_i}&=&\arr{ \dst B^2_{\varepsilon_i}\times \cA^i_{\bar\theta^i}\times \chi^i_{\bar\theta^i}\times \torus_{\bar s^i}^4\quad\quad i=n-1\\
    \\
    \dst B^2_{\varepsilon _i}\times \cA^i_{\bar\theta^i}\times \chi^i_{\bar\theta^i}\times \torus_{\bar s^i}^5\qquad i=1,\cdots, n-2\ \ \&\ \ n\ge 3 }\nonumber\\
    \ovl{\phi_{\eufm{int}}^{i}}&:&\quad \arr{\dst \Theta_{i}=\frac{p_{i}}{\b_i}+f_i(p_i, q_i, {\rm y}_i^*)\\
    \\
    \dst\vartheta_{i}-\p=\b_iq_{i}+g_i(p_i, q_i, {\rm y}_i^*)\\
    \\ \dst{\rm y}_i={\rm y}_i^*\\
    \\
    \dst {\rm x}_i={\rm x}_i^* + \f_i(p_i, q_i, {\rm y}_i^*) } \nonumber\\
    \ovl{{\rm h}_{\eufm{sec}}^{i}}&=&\cA_i\Big[{\rm E}_i+\Omega_i\frac{p_i^2+q_i^2}{2}+\t_i(\frac{p_i^2+q_i^2}{2})^2+{\rm O}(p_i,q_i)^6\Big]\eeqa where ${\rm X}_{\bar\theta^i}^i\times{\rm A}_{\bar\theta^i}^i$ denote the projection of the set ${\rm X}_{\bar\theta}\times {\rm A}_{\bar\theta}$ over the coordinates ${\rm y}_i$ in \equ{yx}, $\bar\theta:=\theta/2$, $\bar s:=s/2$, $f_i$, $g_i$ are ${\rm O}(p_i, q_i)^3$, odd in $(p_i, q_i)$, $\varphi_i$ is ${\rm O}(p_i, q_i)^2$, and \beqa{coefficients} \varepsilon_{i}&=&{\eufm c}_i\,\sqrt{\theta_i}\nonumber\\
    \b_i&:=& \sqrt[4]{\frac{{5\chi_{i-1}\L^2_{i}-(\chi_{i-1}-\chi_i)^2(4\chi_{i-1}-\chi_i)}}{\chi_{i-1}^2\chi_i^2(\chi_{i-1}+\chi_i)} }\qquad \nonumber\\
    \cA_i&:=& m_{i} m_{i+1} \frac{a_{i}^2}{4a_{i+1}^3} \nonumber\\
    {\rm E}_i&:=&-\frac{\L_{i+1}^3}{2(\chi_i-\chi_{i+1})^3} \Big( 5-3\frac{(\chi_{i-1}-\chi_{i})^2}{\L_{i}^2} \Big)\quad \quad \nonumber\\
    \Omega_i&:=&\frac{3\L_{i+1}^3}{\chi_i\L_i^2(\chi_i-\chi_{i+1})^3} \sqrt{ ({5\chi_{i-1}\L^2_{i}-(\chi_{i-1}-\chi_i)^2(4\chi_{i-1}-\chi_i)})(\chi_{i-1}+\chi_i) }\quad \quad \nonumber\\
    \t_i&:=&\frac{\L_{i+1}^3}{\chi_i^2(\chi_i-\chi_{i+1})^3}\Big[ -{\frac{9}{16}}\frac{(\chi_{i-1}-\chi_{i})^{2}(3\chi_{i-1}-\chi_{i})(5\chi_{i-1}+\chi_{i})}{\chi_{i-1}^3\chi_i\L_{i}^2\b_i^4}\nonumber\\
    &&-{\frac{3}{8}}\frac{2\chi_{i-1}^3+9\chi_{i-1}^2\chi_i+2\chi_{i-1}\chi_i^2+\chi_i^3}{\chi_{i-1}\L_{i}^2}-{\frac{3}{16}}\frac{\chi_{i-1}\chi_i^2}{\L_{i}^2}(4 \chi_{i-1}+ \chi_{i})\b_i^4\Big]\eeqa with $\chi_n\equiv 0$, $\bar{\eufm c}_{i}$ depending at most on the ratios $a_i^+/a_i^-$, the masses $m_1$, $\cdots$, $m_n$ and, as usual, $\sqrt[m]{z}$ denoting the principal determination of the $m^{\rm th}$ root of a complex number $z$.
\end{lemma}
\begin{proof} Since the formula for $\ovl{\ovl{f^{n-1, n}_\cP}}$ coincides with the one for $\ovl{f^{n-1, n}_\cP}$ taking $\chi_n\equiv 0$, we shall only work on the terms $\ovl{\ovl{f^{i, i+1}_\cP}}$'s. 

\nl Let ${\rm y}_i$ be as in \equ{yx}, and let \beq{old domain}{\eufm D}_i:\quad (\Theta_i, \vartheta_i)\in{\rm T}^i_{\Theta_i^+, \vartheta_i^+}\qquad {\rm y}_i\in\cA^i_{\theta_i}\times {\rm X}^i_{\theta_i}\qquad {\rm x}_i\in \torus^{m_i}_{s} \eeq where ${\rm T}^i_{\Theta_i^+, \vartheta_i^+}$ is the projection of ${\rm T}_{\Theta^+, \vartheta^+}$ over the coordinates $(\Theta_i, \vartheta_i)$, while $m_i$ is 4 or 5, accordingly to \equ{Nn}. We shall obtain the transformation $\ovl{\phi^i_{\eufm{int}}}$ in \equ{tr*} as a product $\ovl{\phi^i_{\eufm{int}}}=\ovl{\phi^i_{\eufm{diag}}}\circ\ovl{\phi^i_{\eufm{bir}}}$, where $\ovl{\phi^i_{\eufm{diag}}}$ and $\ovl{\phi^i_{\eufm{bir}}}$ are described below.

\nl A Taylor expansion of $\ovl{\ovl{f^{i, i+1}_\cP}}$ around $(\Theta_i,\vartheta_i)=(0,\p)$ gives \beqa{exp} \ovl{\ovl{f^{i, i+1}_\cP}}&=&\cA_i\Big[{\rm E}_i+\O_i\frac{\b_i^2\Theta_i^2+\frac{(\vartheta_i-\p)^2}{\b_i^2}}{2}+\cR_i\Big] \eeqa where $\dst \cA_i$, ${\rm E}_i$, $\b_i$, $\O_i$ are as in \equ{coefficients}. Note that $\b_i$, $\O_i$ are well defined under the assumption \equ{cond on gamma}. The expansion in \equ{exp} shows that $(\Theta_i, \vartheta_i)=(0,\p)$ is an {\sl elliptic} equilibrium point for $\ovl{\ovl{f^{i, i+1}_\cP}}$. The remainder $\cR_i$ is given by \beqano \cR_i&=& {\rm F}\Big[-\frac{3}{2}\frac{4\Theta_i^2-\chi_i^2}{\L_{i}^2}\Big(\frac{(\chi_i^2-\chi_{i-1}^2)^2}{(\sqrt{\chi_i^2-\Theta_i^2}+\sqrt{\chi_{i-1}^2-\Theta_i^2})^2}\nonumber\\
&+&2\sqrt{(\chi_i^2-\Theta_i^2)(\chi_{i-1}^2-\Theta_i^2)}(1+\cos{\vartheta_i})\Big)+\frac{1}{2}\frac{(\chi_i^2-\Theta_i^2)(\chi_{i-1}^2-\Theta_i^2)}{\L_{i}^2}\sin^2{\vartheta_i}\Big] \eeqano where the symbol ${\rm F}$ on the left means that only terms of the fourth order in $(\Theta_i, \vartheta_i-\p)$ have to be included. The lower order expansion of $\cR_i$ is $$\cR_i=\t_{1, i}\Theta_i^4+\t_{2, i}(\vartheta_i-\p)^2\Theta_i^2+\t_{3, i}{(\vartheta_i-\p)^4}+{\rm O}(\Theta_i,\vartheta_i-\p)^6$$ with \beqano \t_{1, i}&:=&\t_1({\rm y}_i):=-\frac{3(\chi_{i-1}-\chi_{i})^{2}(3\chi_{i-1}-\chi_{i})(5\chi_{i-1}+\chi_{i})}{8\chi_{i-1}^3\chi_i\L^2_i}\nonumber\\
\t_{2, i}&:=&\t_2({\rm y}_i):=-\frac{3(2\chi_{i-1}^3+9\chi_{i-1}^2\chi_i+2\chi_{i-1}\chi_i^2+\chi_i^3)}{4\chi_{i-1}\L^2_i}\nonumber\\
\t_{3, i}&:=&\t_3({\rm y}_i):=-\frac{\chi_{i-1}\chi_i^2}{8\L^2_i}(4 \chi_{i-1}+ \chi_{i}).\eeqano We introduce the generating function $$S_{\eufm{diag}, i}(\tilde p_i, \tilde{\rm y}_i, \vartheta_i, {\rm x}_i)=\frac{\tilde p_i(\vartheta_i-\p)}{\tilde\b_i}+\tilde{\rm y}_i{\rm x}_i.$$ It generates the canonical transformation \[\ovl{\phi_{\eufm{diag}}^i}:\quad \Theta_i=\frac{\tilde p_i}{\tilde\b_i}\qquad \vartheta_i-\p= \tilde\b_i \tilde q_i,\quad {\rm y}_i=\tilde {\rm y}_i,\quad {\rm x}_i=\tilde {\rm x}_i+\frac{\partial_{{\rm y}_i}\b_i(\tilde {\rm y}_i)}{\b_i(\tilde {\rm y}_i)} \tilde p_i\tilde q_i\] which transforms $\ovl{\ovl{f^{i, i+1}_\cP}}$ into \beqa{diag} \ovl {\ovl{f_{\eufm{diag}, i}}}=\ovl{\ovl{f^{i, i+1}_\cP}}\circ\ovl{\phi_{\eufm{diag}}^i}&=&\tilde\cA_i\Big[\tilde {\rm E}_i+\tilde \O_i\frac{\tilde p_i^2+\tilde q_i}{2}+\tilde\cR_i\Big] \eeqa with \beqano \tilde\b_i&:=&\b(\tilde{\rm y}_i),\quad \tilde\cA_i:=\cA_i(\tilde{\rm y}_i),\quad \tilde {\rm E}_i:=C(\tilde{\rm y}_i),\quad \tilde \O_i:=\O(\tilde{\rm y}_i),\nonumber\\
\nonumber\\
\tilde\cR_i&:=&\cR_i\circ\phi_{\eufm{diag}}^i=\tilde\t_{1, i}\tilde p_i^4+\tilde\t_{2, i}\tilde p_i^2\tilde q_i^2+\tilde\t_{3, i}\tilde q_i^4+{\rm O}(\tilde p_i, \tilde q_i)^6\nonumber\\
\nonumber\\
\tilde \t_{1, i}&:=&\frac{\t_{1}(\tilde{\rm y}_i)}{\tilde\b_i^4},\qquad \tilde \t_{2, i}:={\t_{2}(\tilde{\rm y}_i)},\qquad \tilde \t_{3, i}:={\t_{3}(\tilde{\rm y}_i)}{\tilde\b_i^4} \eeqano To compute the domain of $\ovl{\phi_{\eufm{diag}}^i}$, we use the following inequalities, which readily follow from the definitions: $$\hat{\eufm c}_i\frac{\sqrt{\theta_i}}{{\rm G}_n^+}\le |\b_i|\le \frac{1}{\hat{\eufm c}_i}\frac{\sqrt{\theta_i}}{{\rm G}_n^-}$$ and $$|\frac{\partial_{{\rm y}_i}\b_i(\tilde {\rm y}_i)}{\b_i(\tilde {\rm y}_i)} |\le \frac{1}{\hat{\eufm c}_i{\theta_i}}.$$ We then see that, choosing a suitable $ \tilde{\eufm c}_i\le \hat{\eufm c}_i$, and the domain $$\widetilde{\ovl{{\eufm B}_i}}:\qquad |(\tilde p_i, \tilde q_i)|\le \tilde\varepsilon_i=\tilde{\eufm c}_i \sqrt{\theta_i}\qquad \tilde{\rm y}_i\in\cA^i_{\theta_i}\times {\rm X}^i_{\theta_i}\qquad \tilde{\rm x}_i\in \torus^{m_i}_{\frac{3}{4} s}$$ inequalities\footnote{Compare \equ{Choice of parameters}.} \equ{old domain} are verified, as desired. Now we look for another canonical transformation $$\ovl{\phi^i_{\eufm{bir}}}:\quad (p^*_i, q^*_i, {\rm y}^*_i, {\rm x}^*_i)\to (\tilde p_i, \tilde q_i, \tilde{\rm y}_i, \tilde{\rm x}_i)\qquad ({\rm y}_i^*=\tilde{\rm y}_i)$$ defined in a analogous domain $$\ovl{{\eufm B}_i}^*:=\ovl{{\eufm B}_i}:\qquad |(p^*_i, q^*_i)|\le\varepsilon_i= {\eufm c}_i^* \sqrt{\theta_i}\qquad {\rm y}^*_i\in\cA^i_{\theta_i}\times {\rm X}^i_{\theta_i}\qquad {\rm x}^*_i\in \torus^{m_i}_{\frac{s}{2}}$$ with ${\eufm c}_i^*=:{\eufm c}_i\le\tilde{\eufm c}_i/2$, such that $$\ovl{\ovl{f_{\eufm{diag}, i}}}\circ\ovl{\phi^i_{\eufm{bir}}}=\ovl{{\rm h}_{\eufm{sec}}^{i}} $$ satisfies the thesis of the lemma. We aim to apply Theorem \ref{perturbation argument}, with $${\rm h}=\tilde {\rm E}_i+\tilde \O_i\frac{\tilde p_i^2+\tilde q_i}{2},\qquad f= \cR_i,\quad \varepsilon=2{\eufm c}_i^*\sqrt{\theta_i},\quad \bar\varepsilon={\eufm c}_i^*\sqrt{\theta_i} .$$ We have to check that inequalities \equ{smallness integration} are satisfied. We can take ${\eufm a}$ and ${\eufm e}$ as it follows from the following inequalities, which, in turn, are easily implied by the definitions \beqano \inf_{{\eufm B}^*_i}|\partial {\rm h}|&=&\inf_{{\eufm B}^*_i}|\tilde\O_i|\ge \frac{\check{\eufm c}_i|{\rm G}_n^-|^2}{\theta_i}=:{\eufm a}\nonumber\\
\nonumber\\
\sup_{{\eufm B}^*_i}|\tilde\cR_i|&\le& \sup_{{\eufm D}^*_i}|\cR_i|\le \frac{1}{\check{\eufm c}_i}\max \sup_{{\eufm D}^*_i}\big\{\frac{(\Theta^*_i)^4}{({\rm G}_n^-)^2},\ {(\vartheta^*_i-\p)^2\Theta_i^2},\ ({\rm G}_n^+)^2(\vartheta^*_i-\p)^4\big\}\nonumber\\
&\le& \frac{({\eufm c}^*)^4({\rm G}_n^+)^2}{\bar{\eufm c}_i}=:{\eufm e} \eeqano Here, we have used that, for $|(p^*_i, q^*_i)|\le 2{\eufm c}^*\sqrt{\theta_i}$, $(\Theta_i^*, \vartheta_i^*):=(\ovl{\phi_{\eufm{diag}}^i})^{-1}(p_i^*, q_i^*)$ verifies \beqano &&|\Theta_i^*|=\frac{|p^*_i|}{|\b_i|}\le 2{\eufm c}^*\sqrt{\theta_i}\frac{{\rm G}_n^+}{\hat{\eufm c}_i\sqrt{\theta_i}}=2\frac{{\eufm c}^*{\rm G}_n^+}{\hat{\eufm c}_i}\nonumber\\
&&|\vartheta^*_i-\p|={|q^*_i|}{|\b_i|}\le \frac{2{\eufm c}^*\sqrt{\theta_i}}{{\eufm c_1}}\frac{\sqrt{\theta_i}}{{\rm G}_n^-}=2\frac{{\eufm c}^*}{\hat{\eufm c}_i}\frac{\theta_i}{{\rm G}_n^-}. \eeqano We then have that condition \equ{smallness integration} holds, provided one takes $${\eufm c}^*:=\min\big\{\frac{{\rm G}_n^-}{{\rm G}_n^+}\sqrt{\check{\eufm c}_i\bar{\eufm c}_i},\ \frac{\tilde{\eufm c}_i}{2} \big\} .$$ From \equ{diag}, one easily computes that the fourth order term of $\ovl{{\rm h}_{\eufm{sec}}^{i}} $ corresponds to be as in \equ{Nn}, with $$\t_i=\frac{3}{2}\t^*_{1, i}+\frac{1}{2}\t^*_{2, i}+\frac{3}{2}\t^*_{3, i}\qquad \t^*_{j, i}:=\tilde\t_{j,i}({\rm y}_i^*). $$ Finally, properties \equ{symmetries} easily follow from the construction.\end{proof}

\subsection{Construction of ${\phi_{\eufm{int}}^{1}}$, $\cdots$, ${\phi_{\eufm{int}}^{n-2}}$ $(n\ge 3)$}\label{local coordinates***} We have to solve \equ{integration h}, assuming that Proposition \ref{claim} holds, up to rank $h-1$. Accordingly to \equ{split step h}, \equ{zero approximation} and letting $$\ovl{\Phi_{\eufm{int}}^{n-h+1}}:=\ovl{\phi_{\eufm{int}}^{n-h+1}}\circ\cdots\circ\ovl{\phi_{\eufm{int}}^{n-1}}$$ we may split \beqano f_{\eufm{norm}, h-1}^{n-h} &=&\sum_{j=n-h+1}^n\P_{h-1}\big(\ovl{f_\cP^{n-h, j}}^{(\ge 2)}\circ\ovl{\Phi_{\eufm{int}}^{n-h+1}}\big) +\widetilde{f_{\eufm{sec}, h-1}^{n-h}}\nonumber\\
&=&\P_{h-1}\big(\ovl{f_{\cP}^{n-h, n-h+1}}^{(2)}\circ\ovl{\Phi_{\eufm{int}}^{n-h+1}}\big)+\P_{h-1}\big(\ovl{f_{\cP}^{n-h, n-h+1}}^{(\ge3)}\circ\ovl{\Phi_{\eufm{int}}^{n-h+1}}\big)\nonumber\\
&+&\sum_{j=n-h+2}^n\P_{h-1}\big(\ovl{f_{\cP}^{n-h, j}}^{(\ge2)}\circ\ovl{\Phi_{\eufm{int}}^{n-h+1}}\big)+\widetilde{f_{\eufm{sec}, h-1}^{n-h}}\nonumber\\
&=&\ovl{\ovl{f_{\cP}^{n-h, n-h+1}}}^{(2)}+\P_{h-1}\big(\widetilde{\ovl{f_{\cP}^{n-h, n-h+1}}}^{(2)} \circ\ovl{\Phi_{\eufm{int}}^{n-h+1}}\big) \nonumber\\
&+&\P_{h-1}\big(\ovl{f_{\cP}^{n-h, n-h+1}}^{(\ge3)}\circ\ovl{\Phi_{\eufm{int}}^{n-h+1}}\big)+\sum_{j=n-h+2}^n\P_{h-1}\big(\ovl{f_{\cP}^{n-h, j}}^{(\ge2)}\circ\ovl{\Phi_{\eufm{int}}^{n-h+1}}\big)\nonumber\\
&+&\widetilde{f_{\eufm{sec}, h-1}^{n-h}} \eeqano where \beqano \ovl{f_{\cP}^{n-h, n-h+1}}^{(\ge3)}&:=&\ovl{f_{\cP}^{n-h, n-h+1}}^{(\ge2)}-\ovl{f_{\cP}^{n-h, n-h+1}}^{(2)}\nonumber\\
\widetilde{\ovl{f_{\cP}^{n-h, n-h+1}}}^{(2)}&:=&{\ovl{f_{\cP}^{n-h, n-h+1}}}^{(2)}-\ovl{\ovl{f_{\cP}^{n-h, n-h+1}}}^{(2)} \eeqano and $\ovl{\ovl{f_{\cP}^{n-h, n-h+1}}}^{(2)}$ as in Lemma \ref{integration***}. Note that we have used that $\ovl{\ovl{f_{\cP}^{n-h, n-h+1}}}^{(2)}$ is left unvaried by $\ovl{\Phi_{\eufm{int}}^{n-h+1}}$. Let $\ovl{\eufm{B}_{n-h}}$, $\ovl{\phi^{n-h}_{\eufm{int}}}$ be as in Lemmata \ref{integration***}, with the symbols $(\Theta_{n-h}$, $\vartheta_{n-h})$, ${\rm y}_{n-h}$, ${\rm x}_{n-h}$ of that lemma corresponding to \beqano &&\ovl{\rm t}^{{(h-1)}}_{n-h}:=\big(\Theta_{n-h}^{(h-1)},\ \vartheta_{n-h}^{(h-1)}\big)\nonumber\\
&&\ovl{\rm y}_{n-h}^{(h-1)}:=\big(\chi^{(h-1)}_{n-h-1},\ \chi^{(h-1)}_{n-h},\ \chi^{(h-1)}_{n-h+1},\ \L^{(h-1)}_{n-h},\ \L^{(h-1)}_{n-h+1}\big)\nonumber\\
&&\ovl{\rm x}_{n-h}^{(h-1)}:=\big(\k^{(h-1)}_{n-h-1},\ \k^{(h-1)}_{n-h},\ \k^{(h-1)}_{n-h+1},\ \ell^{(h-1)}_{n-h},\ \ell^{(h-1)}_{n-h+1}\big) \eeqano and the symbols $(p_{n-h}$, $q_{n-h})$, ${\rm y}_{n-h}^*$, ${\rm x}_{n-h}^*$ to \beqano &&\ovl{\rm z}^{*\pph}_{n-h}:=\big(p_{n-h}^{*\pph},\ q_{n-h}^{*\pph}\big)\nonumber\\
&&\ovl{\rm y}_{n-h}^{*\pph}:=\big(\chi^{*(h)}_{n-h-1},\ \chi^{*(h)}_{n-h},\ \chi^{*(h)}_{n-h+1},\ \L^{*(h)}_{n-h},\ \L^{*(h)}_{n-h+1}\big)\nonumber\\
&&\ovl{\rm x}_{n-h}^{*\pph}:=\big(\k^{*(h)}_{n-h-1},\ \k^{*(h)}_{n-h},\ \k^{*(h)}_{n-h+1},\ \ell^{*(h)}_{n-h},\ \ell^{*(h)}_{n-h+1}\big). \eeqano Defining \beqano &&{\rm t}^{*(h)}:=\big(\Theta^{*(h)},\ \vartheta^{*(h)}\big)\nonumber\\
&&{\rm z}^{*(h)}:=\big(p^{*(h)},\ q^{*(h)}\big)\nonumber\\
&&{\rm y}^{*(h)}:=\big(\chi^{*(h)}, \L^{*(h)}\big)\nonumber\\
&&{\rm x}^{*(h)}:=\big(\k^{*(h)}, \ell^{*(h)}\big) \eeqano in an alagous way as in \equ{new coordinates*}, we regard $\ovl{\phi^{n-h}_{\eufm{int}}}$ as a map on the set \beqano \ovl{\eufm{D}_{\eufm{int}}^h}&:=&\Big\{\big({\rm t}^{*\pph},\ {\rm z}^{*\pph},\ {\rm y}^{*\pph},\ {\rm x}^{*\pph}\big) :\ \big(\ovl{\rm z}_{n-h}^{*\pph},\ \ovl{\rm y}_{n-h}^{*\pph},\ \ovl{\rm x}_{n-h}^{*\pph} \big)\in\ovl{\eufm{B}_{n-h}}\Big\} \eeqano extended via the identity on the extra-coordinates. We then have that $\ovl{\phi^{n-h}_{\eufm{int}}}$ transforms $f_{\eufm{sec}, h-1}^{n-h}$ into \beqano \ovl{f_{\eufm{sec}, h-1}^{n-h}}:=f_{\eufm{sec}, h-1}^{n-h}\circ\ovl{\phi^{n-h}_{\eufm{int}}}&=& \ovl{{\rm h}^{n-h}_{\eufm{sec}}}+\ovl{{\rm f}^{n-h}_{\eufm{sec}}} \eeqano where \beqa{perturbation} \ovl{{\rm f}^{n-h}_{\eufm{sec}}}&:=&\P_{h-1}\big(\widetilde{\ovl{f_{\cP}^{n-h, n-h+1}}}^{(2)} \circ\ovl{\Phi_{\eufm{int}}^{n-h}}\big) \nonumber\\
&&+\P_{h-1}\big(\ovl{f_{\cP}^{n-h, n-h+1}}^{(\ge3)}\circ\ovl{\Phi_{\eufm{int}}^{n-h}}\big)+\sum_{j=n-h+2}^n\P_{h-1}\big(\ovl{f_{\cP}^{n-h, j}}^{(\ge2)}\circ\ovl{\Phi_{\eufm{int}}^{n-h}}\big)\nonumber\\
&&+\widetilde{f_{\eufm{sec}, h-1}^{n-h}}\circ\ovl{\phi_{\eufm{int}}^{n-h}}. \eeqa Here, we have used $\ovl{\Phi_{{\eufm{int}}}^{n-h}}=\ovl{\Phi_{\eufm{int}}^{n-h+1}}\circ\ovl{\phi_{\eufm{int}}^{n-h}}$; that $\P_{h-1}$ and $\ovl{\phi_{\eufm{int}}^{n-h}}$ commute and observe that $\ovl{\eufm{D}_{\eufm{int}}^h}$ has the form of ${\eufm{D}_{\eufm{int}}^h}$ in \equ{domain Dsec***}, with $\hat{\eufm c}_h$ replaced by a suitable $\hat{\eufm c}'_h$ of the same form. The function $\ovl{{\rm f}_{\eufm{sec}}^{n-h}}$ satisfies the following two properties:
\begin{enumerate}
    \item[\tiny\textbullet] It depends only on $$\big(p^{*(h)}_{n-h},\ q^{*(h)}_{n-h},\ \tilde{\rm y}^{*(h)}\big)$$ where $\tilde{\rm y}^{*(h)}$ is defined analogously to \equ{dependence on variables};

 \item[\tiny\textbullet] is uniformly bounded by the right hand side of the first inequality in \equ{h-1 bound} (this follows from the definition in \equ{perturbation}); \item[\tiny\textbullet] is even for $$(p^{*(h)}_{n-h}, q^{*(h)}_{n-h})\to -(p^{*(h)}_{n-h}, q^{*(h)}_{n-h}).$$
\end{enumerate}
Proceeding in a similar way as we did for the construction of $\ovl{\phi^i_{\eufm{bir}}}$ in the proof of Lemma \ref{base integration}, we may apply Theorem \ref{perturbation argument}, with \beqano &&{\rm h}=\ovl{{\rm h}^{n-h}_{\eufm{sec}}},\quad f=\ovl{{\rm f}^{n-h}_{\eufm{sec}}},\quad (P, Q)=(p^{*(h)}_{n-h}, q^{*(h)}_{n-h})\nonumber\\
&&(P',Q')=\hat{\rm z}^{*\pph}_{n-h},\quad {\rm y}={\rm y}^{*\pph}_{n-h},\quad {\rm x}={\rm x}^{*\pph}_{n-h}. \eeqano with ${\rm y}^{*\pph}_{n-h}$, ${\rm x}^{*\pph}_{n-h}$ defined analogously to ${\rm y}^{\pph}_{*,n-h}$, ${\rm x}^{\pph}_{*,n-h}$ in \equ{true dependence} and $\hat{\rm z}^{*\pph}_{n-h}$ defined analogously to $\hat{\rm z}^{\pph}_{*,n-h}$ in \equ{dependence on variables}. We then find another domain $ {{\eufm D}_{\eufm{int}}^h}$ as in \equ{domain Dsec***} and another real-analytic transformation \beqano {\phi^{n-h}_{*,\eufm{int}}}:\quad&& \big( {\rm t}_*^\pph, {\rm z}_*^\pph, {\rm y}_*^\pph, {\rm x}_*^\pph\big)\in{{\eufm D}_{\eufm{int}}^h}\to \big( {\rm t}^{*\pph}, {\rm z}^{*\pph}, {\rm y}^{*\pph}, {\rm x}^{*\pph}\big) \in\ovl{{\eufm D}_{\eufm{ int}, h}}\eeqano such that \beqano \widetilde{f_{\eufm{sec}, h-1}^{n-h}}:=\ovl{f_{\eufm{sec}, h-1}^{n-h}}\circ{\phi^{n-h}_{*,\eufm{int}}}=f_{\eufm{sec}, h-1}^{n-h}\circ\ovl{\phi^{n-h}_{\eufm{int}}}\circ{\phi^{n-h}_{*,\eufm{int}}}&=& {{\rm h}^{n-h}_{\eufm{sec}}} \eeqano as desired, depends only on $\hat{\rm y}^\pph_{*, n-h}$ in \equ{good transformation}, and hence \equ{integration h} is satisfied. That ${\phi^{n-h}_{*,\eufm{int}}}$ may be also chosen of a form {analogous} to \equ{integration}, with $\Theta^{(h-1)}_{n-h}$, $\vartheta^{(h-1)}_{n-h}$, $\hat{\rm z}^{(h-1)}$, ${\rm y}^{(h-1)}_{n-h}$, ${\rm x}^{(h-1)}_{n-h}$ replaced by $p^{*\pph}$, $q^{*\pph}$, $\hat{\rm z}^{*\pph}$, ${\rm y}^{*\pph}$, ${\rm x}^{*\pph}$ also easily follows from the properties bove. Therefore the composition $$\phi^{n-h}_{\eufm{int}}:=\ovl{\phi^{n-h}_{\eufm{int}}}\circ{\phi^{n-h}_{*,\eufm{int}}}$$ has again the form in \equ{integration} and satisfies \equ{integration h}, as wanted. \qed

\subsection{Construction of $\phi_{\eufm{norm}}^1$, $\cdots$, $\phi_{\eufm{norm}}^{n-2}$ ($n\ge 3$)}\label{construction of normalizations}

In this section we aim to determine, for $n\ge 3$ and $1\le h\le n-2$, a transformation $\phi_{\eufm{norm}}^{n-h}$ solving \equ{normalization}-\equ{normalization***}, assuming the Proposition \ref{claim} holds up to rank $(h-1)$ and that $\phi_{\eufm{int}}^{n-h}$ has been constructed.

\vskip.1in \noi We switch from the coordinates $(\chi_*^\pph, \k_*^\pph)$ defined implicitly via the right hand side of \equ{normalization} to the auxiliary coordinates $${\rm G}_\eufm{aux}^\pph=({\rm G}_{\eufm{aux},1}^\pph, \cdots, {\rm G}_{\eufm{aux},n}^\pph),\quad {\rm g}_\eufm{aux}^\pph=({\rm g}_{\eufm{aux}, 1}^\pph, \cdots, {\rm g}_{\eufm{aux}, n}^\pph)$$ defined via the linear transformation \beq{transf*} \phi_{\eufm{aux}}^{n-h}:\quad \arr{\chi_{*,i-1}^\pph={\rm G}_{\eufm{aux},i}^\pph+\cdots+{\rm G}_{\eufm{aux},n}^\pph\\
\\
\k^\pph_{*, i-1}={\rm g}_{\eufm{aux}, i}^\pph-{\rm g}_{\eufm{aux}, i-1}^\pph}\eeq with $1\le i\le n$ and ${\rm g}_{\eufm{aux}, 0}:=0$. We regard $\phi_{\eufm{aux}}^{n-h}$ as a transformation on all the coordinates, extending it as the identity on the remaining ones. We denote the new coordinates as \beqano {\rm t}_{\eufm{aux}}^\pph&:=&\arr{\big({\Theta}_{\eufm{aux}, 1}^\pph,\ \cdots,\ \Theta^\pph_{\eufm{aux}, n-h-1},\ {\vartheta}_{\eufm{aux}, 1}^\pph,\ \cdots,\ \vartheta^\pph_{\eufm{aux}, n-h-1}\big)\\
\qquad n\ge 4,\ 2\le h\le n-2\\
\emptyset\qquad{\rm otherwise} } \nonumber\\
{{\rm z}}_{\eufm{aux}}^\pph&:=&(p_{\eufm{aux}, n-h}^\pph,\ \cdots,\ p_{\eufm{aux}, n-1}^\pph,\ q_{\eufm{aux}, n-h}^\pph,\ \cdots,\ q_{\eufm{aux}, n-1}^\pph)\nonumber\\
{{\rm y}}_{\eufm{aux}}^\pph&:=& \big( {\rm G}^\pph_{\eufm{aux}, 1},\ \cdots,\ {\rm G}_{\eufm{aux}, n}^\pph,\ \L^\pph_{\eufm{aux}, {1}},\ \cdots,\ \L^\pph_{\eufm{aux}, n}\big)\nonumber\\
{{\rm x}}_{\eufm{aux}}^\pph&:=& \big( {\rm g}^\pph_{\eufm{aux},1},\ \cdots,\ {\rm g}_{\eufm{aux}, n}^\pph,\ \ell^\pph_{\eufm{aux}, 1},\ \cdots,\ \ell^\pph_{\eufm{aux}, n}\big) \eeqano the new Hamiltonian as 
\begin{equation}\label{new f***}
\begin{split}{{f_{\eufm{sec},\eufm{int}, \eufm{aux}, h-1}}}({\rm t}_\eufm{aux}^\pph,{\rm z}_{\eufm{aux}}^\pph, {\rm y}_\eufm{aux}^\pph,{\rm x}_\eufm{aux}^\pph) :={{f_{\eufm{sec},\eufm{int}, h-1}}}\circ\phi_{\eufm{aux}}^{n-h}({\rm t}_\eufm{aux}^\pph,{\rm z}_{\eufm{aux}}^\pph, {\rm y}_\eufm{aux}^\pph, {\rm x}_\eufm{aux}^\pph).\end{split} \end{equation} Now we define the domain where we want to consider ${{f_{\eufm{sec},\eufm{int}, \eufm{aux}, h-1}}}$. Firstly, we let \[{\eufm D}_{\eufm {int},\eufm{aux}}^h:=\Big\{ ({\rm t}_{\eufm{aux}}^\pph, {\rm z}_{\eufm{aux}}^\pph, {\rm y}_{\eufm{aux}}^\pph, {\rm x}_{\eufm{aux}}^\pph):\quad ({\rm t}_{*}^\pph, {\rm z}_{*}^\pph, {\rm y}_{*}^\pph, {\rm x}_{*}^\pph)\in{\eufm D}_{\eufm {int}}^h \Big\}\] where ${\eufm D}_{\eufm {int}}^h$ is defined in \equ{domain Dsec***}. Then ${\eufm D}_{\eufm {int},\eufm{aux}}^h$ is given by \beqano{\eufm D}_{\eufm {int},\eufm{aux}}^h={\rm T}^h_{\hat{\eufm c}_h\theta}\times B^2_{\hat{\eufm c}_h\varepsilon_{n-h}} \times B^{*,h-1}_{\hat{\eufm c}_h\varepsilon}\times({\rm G}_{*})_{ \ovl{\eufm c}_h\theta, \ovl{\eufm c}_h\bar r}\times{\rm A}_{\ovl{{\eufm c}}_h\widetilde r^\pph}\times \torus^n_{\ovl{\eufm c}_h s}\times \torus^n_{\ovl{\eufm c}_h s}, \eeqano with $$({\rm G}_*)_{ \ovl{\eufm c}_h\theta, \ovl{\eufm c}_h\bar r}:=({\rm G}_1)_{\hat{\eufm c}_1\theta_1}\times\cdots\times({\rm G}_{n-h})_{\hat{\eufm c}_{n-h}\theta_{n-h}}\times ({\rm G}^*_{n-h+1})_{\hat{\eufm c}_{n-h+1}\ovl r_{n-h+1}}\times\cdots\times ({\rm G}^*_{n-1})_{\hat{\eufm c}_{n-1}\ovl r_{n-1}}.$$ Next, for $1\le h'\le h$ and any fixed $\bar\g$, $\bar K>0$ and $\bar\t>2$, we define 
\begin{equation}\label{omega sec n-h} \begin{split}
&\o_{\eufm{ sec}}^{n-h'}(\hat{\rm y}_{\eufm{aux}, n-h'}^{(h)}):=\\
&\left\{
\begin{array}
    {lllll} \dst \partial_{\frac{(p_{\eufm{aux}, n-1}^\ppu)^2+(q_{\eufm{aux}, n-1}^\ppu)^2}{2}, {\rm G}_{\eufm{aux}, n-1}^\ppu, {\rm G}_{\eufm{aux}, n}^\ppu}\,{{{\rm h}_{\eufm{sec}}^{n-1}}}(\hat{\rm y}_{\eufm{aux}, n-1}^{(h)})\ &n\ge 3,\ h'=1,\ 2\le h\le n-1\\
    \\ \dst \partial_{\frac{(p_{\eufm{aux}, n-h'}^{(h')})^2+(q_{\eufm{aux}, n-h'}^{(h')})^2}{2}, {\rm G}_{\eufm{aux}, n-h'}^{(h')}}\,{{{\rm h}_{\eufm{sec}}^{n-h'}}}(\hat{\rm y}_{\eufm{aux}, n-h}^{(h)})\qquad &n\ge 3,\quad 2\le h'\le h\le n-1,\\
    &(h',\ h)\ne (n-1, n-1)\\
    \\ \dst \partial_{\frac{(p_{\eufm{aux}, 1}^{(n-1)})^2+(q_{\eufm{aux}, 1}^{(n-1)})^2}{2}}\,{{{\rm h}_{\eufm{sec}}^{n-h}}}(\hat{\rm y}_{\eufm{aux}, 1}^{(n-1)})\qquad & h'=h=n-1.
\end{array}
\right.
\end{split}\end{equation} We then choose the following sub-domain of ${\eufm D}_{\eufm{int}, \eufm{aux}}^h$ \beqa{normalized auxiliary domain} {\eufm D}_{\eufm{sec}, \eufm{aux}}^h&:=&\Big\{({\rm t}_{\eufm{norm}, \eufm{aux}}^\pph,{\rm z}_{\eufm{norm}, \eufm{aux}}^\pph, {\rm y}_{\eufm{norm}, \eufm{aux}}^\pph,{\rm x}_{\eufm{norm}, \eufm{aux}}^\pph)\in {\eufm D}_{\eufm {int},\eufm{aux}}^h:\nonumber\\
&&|\o_{\eufm{ sec}}^{n-h'}\cdot {{\rm k}}|\ge \frac{(a_{n-h'}^+)^2}{(a_{n-h'+1}^-)^3\theta_{n-h}}\frac{\bar\g}{\bar K^{\bar\t}},\nonumber\\
&& \forall\ {{\rm k}}\in \integer^j\setminus\{0\},\ |{{\rm k}}|_1\le \bar K,\quad \forall\ 2\le h'\le h \Big\}.\eeqa Here $j$ is chosen to be $3$, $2$ or $1$ accordingly to the three cases above. The set ${\eufm D}_{\eufm {int},\eufm{aux}}^h$ is non-empty, if $\bar\g$ is chosen suitably small. Indeed, if we put $$\hat{\rm y}_{\eufm{aux}, n-h}^{(h)}:=\Big(\frac{(p_{\eufm{aux}, n-h}^\pph)^2+(q_{\eufm{aux}, n-h}^\pph)^2}{2}, {\rm G}_{\eufm{aux}, n-h}^\pph, \L_{\eufm{aux}, n-h}^\pph,\ \hat{\rm y}_{\eufm{aux}, n-h+1}^{(h)} \Big)$$ then standard quantitative arguments show that, for any fixed value $$\big( \bar\L_{\eufm{aux}, n-h}^\pph,\ \hat{\bar{\rm y}}_{\eufm{aux}, n-h+1}^{(h)}\big)\in \P_{ \L_{\eufm{aux}, n-h}^\pph, \hat{\rm y}_{\eufm{aux}, n-h+1}^{(h)}}{\eufm D}_{\eufm {int},\eufm{aux}}^h,$$ the measure of the set $\cN _{n-h}\subset B^2_{\ovl{\eufm c}_h\varepsilon_{n-h}}\times {\rm G}_{n-h} $ of $(p_{\eufm{aux}, n-h}^\pph$, $q_{\eufm{aux}, n-h}^\pph$, ${\rm G}_{\eufm{aux}, n-h}^\pph\big)$ where the inequality in \equ{normalized auxiliary domain} does not hold may be bounded as $$\meas\cN_{n-h}\le \frac{\bar\g}{\eufm c}\meas\big( B^2_{\ovl{\eufm c}_h\varepsilon_{n-h}}\times {\rm G}_{n-h} \big),$$ (where ${\eufm c}$ depends only on the semi-axes ratio and the masses), hence \equ{good set} follows. This is because $\o_{\eufm{ sec}}^{n-h'}(\hat{\rm y}_{\eufm{aux}, n-h}^{(h)})$ is a diffeomorphism (Compare Appendix \ref{Checking the non-degeneracy assumption}).

\nl Now we inspect the form of ${{f_{\eufm{sec},\eufm{int}, \eufm{aux}, h-1}}}$ in \equ{new f***}. Introducing the following symbols 
\beqano {\rm t}_{\eufm{aux}, i}^\pph&:=&\arr{\big({\Theta}_{\eufm{aux}, i}^\pph,\ \cdots,\ \Theta^\pph_{\eufm{aux}, n-h-1},\ {\vartheta}_{\eufm{aux}, i}^\pph,\ \cdots,\ \vartheta^\pph_{\eufm{aux}, n-h-1}\big)\\
\qquad n\ge 4,\ 2\le h\le n-2,\ 1\le i\le n-h-1\\
\\
\emptyset\qquad {\rm otherwise} } \nonumber\\
{{\rm y}}_{\eufm{aux}, i}^\pph&:=& \big( {\rm G}^\pph_{\eufm{aux}, i},\ \cdots,\ {\rm G}_{\eufm{aux}, n}^\pph,\ \L^\pph_{\eufm{aux}, {1}},\ \cdots,\ \L^\pph_{\eufm{aux}, n}\big)\nonumber\\
{{\rm x}}_{\eufm{aux}, i}^\pph&:=& \big( {\rm g}^\pph_{\eufm{aux}, i},\ \cdots,\ {\rm g}_{\eufm{aux}, n}^\pph,\ \ell^\pph_{\eufm{aux}, {1}},\ \cdots,\ \ell^\pph_{\eufm{aux}, n}\big) \nonumber\\
{\hat{\rm y}}_{\eufm{aux}, i}^\pph&:=& \big(\frac{(p^\pph_{\eufm{aux}, i})^2+(q^\pph_{\eufm{aux}, i})^2}{2},\ \cdots ,\ \frac{(p^\pph_{\eufm{aux}, n-1})^2+(q^\pph_{\eufm{aux}, n-1})^2}{2},\ {\rm G}^\pph_{\eufm{aux}, i},\ \cdots,\ {\rm G}_{\eufm{aux}, n}^\pph,\nonumber\\
&&\L^\pph_{\eufm{aux}, i},\ \cdots,\ \L^\pph_{\eufm{aux}, n}\big)\nonumber\\
{\hat{\rm x}}_{\eufm{aux}, i}^\pph&:=&\arr{\big({\rm g}^\pph_{\eufm{aux}, i+1}-{\rm g}^\pph_{\eufm{aux}, i},\ \cdots,{\rm g}^\pph_{\eufm{aux}, n-h}-{\rm g}^\pph_{\eufm{aux}, n-h-1}\big)\quad n\ge4\ \&\ 1\le h-1\le n-3 \\
\\ \emptyset\qquad\qquad{\rm otherwise}}\nonumber\\
{\hat{\rm X}}_{\eufm{aux}, i}^\pph&:=&\arr{\big({\rm G}^\pph_{\eufm{aux}, i},\ \cdots,{\rm G}^\pph_{\eufm{aux}, n-h}\big)\quad n\ge4\ \&\ 1\le h-1\le n-3 \\
\\ \emptyset\qquad\qquad{\rm otherwise}}\nonumber\\
{\rm z}^\pph_{\eufm{aux}, n-h}&:=&\big(p^\pph_{\eufm{aux}, n-h},\ q^\pph_{\eufm{aux}, n-h}\big),\qquad \hat {\rm z}_{\eufm{norm}, j}^{(h)}:=p_{\eufm{norm}, j}^{(h)}+{\rm i}q_{\eufm{norm}, j}^{(h)}\nonumber\\
{\tilde{\rm y}}_{\eufm{aux}, i}^\pph&:=& \big(\frac{(p^\pph_{\eufm{aux}, n-h+1})^2+(q^\pph_{\eufm{aux}, n-h+1})^2}{2},\ \cdots ,\ \frac{(p^\pph_{\eufm{aux}, n-1})^2+(q^\pph_{\eufm{aux}, n-1})^2}{2},\nonumber\\
&&{\rm G}^\pph_{\eufm{aux}, i},\ \cdots,\ {\rm G}_{\eufm{aux}, n}^\pph,\L^\pph_{\eufm{aux}, i},\ \cdots,\ \L^\pph_{\eufm{aux}, n}\big) \nonumber\\
\tilde{\rm y}_\eufm{aux}^\pph&:=&\tilde{\rm y}_{\eufm{aux}, 1}^\pph,\quad \hat{\rm x}_\eufm{aux}^\pph:=\hat{\rm x}_{\eufm{aux}, 1}^\pph,\quad \hat{\rm X}_\eufm{aux}^\pph:=\hat{\rm X}_{\eufm{aux}, 1}^\pph, \eeqano by means of \equ{1st step***}, we have 
\begin{equation}\label{before normalization}\begin{split} {{f_{\eufm{sec},\eufm{int}, \eufm{aux}, h-1}}}({\rm t}_\eufm{aux}^\pph,{\rm z}_{\eufm{aux}}^\pph, {\rm y}_\eufm{aux}^\pph,{\rm x}_\eufm{aux}^\pph) =&{\rm h}_{\eufm{sec}, h}({\hat{\rm y}}^\pph)+ {{f_{\eufm{norm},\eufm{int}, \eufm{aux}, h-1}}}({\rm t}_\eufm{aux}^\pph,{\rm z}_{\eufm{aux},n-h}^\pph,\tilde{\rm y}_\eufm{aux}^\pph,\hat{\rm x}_\eufm{aux}^\pph)\\
=&\sum_{i=n-h}^{n-1}{{{\rm h}_{\eufm{sec}}^i}}({\hat{\rm y}}_{ i}^\pph)\\
+& \sum_{i=1}^{n-h-1}{{f_{\eufm{norm},\eufm{int}, \eufm{aux}, h-1}^i}}({\rm t}_{\eufm{aux}, i}^\pph,{\rm z}_{\eufm{aux},n-h}^\pph,\tilde{\rm y}_{\eufm{aux}, i}^\pph,\hat{\rm x}_{\eufm{aux}, i}^\pph)\\
\end{split}
\end{equation} where we have let \begin{equation}\label{fi norm int aux} \begin{split}&{{f_{\eufm{norm},\eufm{int}, \eufm{aux}, h-1}}} :={{f_{\eufm{norm},\eufm{int}, h-1}}}\circ\phi_{\eufm{aux}}^{n-h} ,\qquad {{f_{\eufm{norm},\eufm{int}, \eufm{aux}, h-1}^i}} :={{f_{\eufm{norm},\eufm{int}, h-1}^i}}\circ\phi_{\eufm{aux}}^{n-h} .
\end{split}
\end{equation} On the domain ${\eufm D}_{\eufm {sec},\eufm{aux}}^h$ specified in \equ{normalized auxiliary domain}, we aim to construct and real-analytic and canonical transformation \beq{auxiliary normalization}\phi^{n-h}_{\eufm{norm}, \eufm{aux}}:\ ({\rm t}_{\eufm{norm}, \eufm{aux}}^\pph,{\rm z}_{\eufm{norm}, \eufm{aux}}^\pph, {\rm y}_{\eufm{norm}, \eufm{aux}}^\pph,{\rm x}_{\eufm{norm}, \eufm{aux}}^\pph)\in {\eufm D}_{\eufm {sec},\eufm{aux}}^h\to ({\rm t}_{\eufm{aux}}^\pph, {\rm z}_{\eufm{aux}}^\pph {\rm y}_{\eufm{aux}}^\pph,{\rm x}_{\eufm{aux}}^\pph)\in {\eufm D}_{\eufm{int}, \eufm{aux}}^h\eeq such that the transformed Hamiltonian \[{{f_{\eufm{sec},\eufm{aux}, h}}} :={{f_{\eufm{sec},\eufm{int}, \eufm{aux}, h-1}}}\circ \phi^{n-h}_{\eufm{norm}, \eufm{aux}} \] has the form \beqano {{f_{\eufm{sec},\eufm{aux}, h}}} &=&{\rm h}_{\eufm{sec}, h}({\hat{\rm y}}_{\eufm{norm}, \eufm{aux}}^\pph)+{{f_{\eufm{norm}, \eufm{aux}, h}}}({\rm t}_{\eufm{norm}, \eufm{aux}}^\pph,\hat{\rm y}_{\eufm{norm}, \eufm{aux}}^\pph,\hat{\rm x}_{\eufm{norm}, \eufm{aux}}^\pph)\nonumber\\
&=&\sum_{i=n-h}^{n-1}{{{\rm h}_{\eufm{sec}}^i}}({\hat{\rm y}}_{\eufm{norm}, \eufm{aux}, i}^\pph)+\sum_{i=1}^{n-h-1}{{f_{\eufm{norm},\eufm{aux}, h}^i}}({\rm t}_{\eufm{norm}, \eufm{aux}, i}^\pph,\hat{\rm y}_{\eufm{norm}, \eufm{aux}, i}^\pph,\hat{\rm x}_{\eufm{norm}, \eufm{aux}, i}^\pph)\nonumber\\
&+&f_{\eufm{exp}, \eufm{sec}, \eufm{aux}, h}({\rm t}_{\eufm{norm}, \eufm{aux}}^\pph,{\rm z}_{\eufm{norm}, \eufm{aux}}^\pph, {\rm y}_{\eufm{norm}, \eufm{aux}}^\pph,{\rm x}_{\eufm{norm}, \eufm{aux}}^\pph) \eeqano where \beqano {\hat{\rm x}}_{\eufm{norm},\eufm{aux}, i}^\pph&:=&\arr{\big({\rm g}^\pph_{\eufm{norm},\eufm{aux}, i+1}-{\rm g}^\pph_{\eufm{norm},\eufm{aux}, i},\ \cdots,{\rm g}^\pph_{\eufm{norm},\eufm{aux}, n-h-1}-{\rm g}^\pph_{\eufm{norm},\eufm{aux}, n-h-2}\big)\\
\\\
qquad{\rm if}\qquad n\ge4\ \&\ 1\le h-1\le n-3 \\
\\ \emptyset\qquad\qquad{\rm otherwise}}\nonumber\\
{\hat{\rm y}}_{\eufm{norm},\eufm{aux}, i}^\pph&:=& \big( \frac{(p^\pph_{\eufm{norm},\eufm{aux}, n-h})^2+(q^\pph_{\eufm{norm},\eufm{aux}, n-h})^2}{2},\ \cdots ,\ \frac{(p^\pph_{\eufm{norm},\eufm{aux}, n-1})^2+(q^\pph_{\eufm{norm},\eufm{aux}, n-1})^2}{2},\nonumber\\
&&{\rm G}^\pph_{\eufm{norm},\eufm{aux}, i},\ \cdots,\ {\rm G}_{\eufm{norm},\eufm{aux}, n}^\pph,\L^\pph_{\eufm{norm},\eufm{aux}, i},\ \cdots,\ \L^\pph_{\eufm{norm},\eufm{aux}, n}\big) \eeqano and $f_{\eufm{exp}, \eufm{sec}, \eufm{aux}, h}$ satisfies the bound for $f_{\eufm{exp}, \eufm{sec}, h}$ in \equ{h-1 bound}. This will conclude the proof, up to apply the inverse transformation of \equ{transf*}, with ${\rm G}_{\eufm{aux},i}^\pph$, ${\rm g}_{\eufm{aux},i}^\pph$, ${\rm \chi}_{*,i}^\pph$, ${\rm \k}_{*,i}^\pph$ replaced by ${\rm G}_{\eufm{norm},\eufm{aux}, i}^\pph$, ${\rm g}_{\eufm{norm},\eufm{aux}, i}^\pph$, ${\rm \chi}_{i}^\pph$, ${\rm \k}_{i}^\pph$, and to take $${\eufm D}_{\eufm {sec}}^h:=\phi_{\eufm{aux}}^{n-h}\big({\eufm D}_{\eufm {sec},\eufm{aux}}^h\big).$$

\nl We shall obtain the transformation $\phi^{n-h}_{\eufm{norm}, \eufm{aux}}$ in \equ{auxiliary normalization} via an application of Proposition \ref{iterative lemma}. Before doing it, we just remark that, since, in our particular case, $f_{\eufm{norm},\eufm{int},\eufm{aux}, h-1}$ depends on ${\rm z}_{\eufm{aux}}^\pph$, ${\rm y}_{ \eufm{aux}}^\pph$, ${\rm x}_{ \eufm{aux}}^\pph$ only via ${\rm z}_{\eufm{aux}, n-h}^\pph$, $\tilde{\rm y}_\eufm{aux}^\pph$, $\hat{\rm x}_\eufm{aux}^\pph$ and is even in $({\rm t}_{ \eufm{aux}}^\pph$, ${\rm z}_{\eufm{aux},n-h}^\pph)$, the proof of Proposition \ref{iterative lemma} can be easily handled to show that $\phi^{n-h}_{\eufm{norm}, \eufm{aux}}$ can be chosen of the form \beqano \phi_{\eufm{norm}, \eufm{aux}}^{n-h}:\ \left\{
\begin{array}
    {lllll} \dst\Theta^{(h)}_{\eufm{aux}, j}=&{\rm F}_{\eufm{norm}, \eufm{aux}, j}^\pph({\rm t}_{\eufm{norm}, \eufm{aux}}^\pph,{\rm z}_{\eufm{norm},\eufm{aux}, n-h}^\pph,\tilde{\rm y}_{\eufm{norm}, \eufm{aux}}^\pph,\hat{\rm x}_{\eufm{norm}, \eufm{aux}}^\pph) \\
    \\ \dst\vartheta^{(h)}_{\eufm{aux}, j}-\p=&{\rm G}_{\eufm{norm}, \eufm{aux}, j}^\pph({\rm t}_{\eufm{norm}, \eufm{aux}}^\pph,{\rm z}_{\eufm{norm},\eufm{aux}, n-h}^\pph,\tilde{\rm y}_{\eufm{norm}, \eufm{aux}}^\pph,\hat{\rm x}_{\eufm{norm}, \eufm{aux}}^\pph)\\
    &j=1,\ \cdots,\ n-h-1 \\
    \\ \dst{{\rm z}}_{\eufm{aux},n-h}^\pph=&{\rm Z}_{\eufm{norm}, \eufm{aux}}^\pph({\rm t}_{\eufm{norm}, \eufm{aux}}^\pph,{\rm z}_{\eufm{norm},\eufm{aux}, n-h}^\pph,\tilde{\rm y}_{\eufm{norm}, \eufm{aux}}^\pph,\hat{\rm x}_{\eufm{norm}, \eufm{aux}}^\pph) \\
    \\ \dst\big({\hat{\rm X}}_{\eufm{aux}}^\pph,\ {\hat{\rm x}}_{\eufm{aux}}^\pph\big)=&{\rm X}_{\eufm{norm}, \eufm{aux}}^\pph({\rm t}_{\eufm{norm}, \eufm{aux}}^\pph,{\rm z}_{\eufm{norm},\eufm{aux}, n-h}^\pph,\tilde{\rm y}_{\eufm{norm}, \eufm{aux}}^\pph,\hat{\rm x}_{\eufm{norm}, \eufm{aux}}^\pph)\\
    \\ \dst \hat {\rm z}_{\eufm{norm}, j}^{(h)}=&\hat {\rm z}_{\eufm{norm}, \eufm{aux},j}^{(h)}e^{{\rm i}\psi_{\eufm{norm}, \eufm{aux},j}^\pph({\rm t}_{\eufm{norm}, \eufm{aux}}^\pph,{\rm z}_{\eufm{norm},\eufm{aux}, n-h}^\pph,\tilde{\rm y}_{\eufm{norm}, \eufm{aux}}^\pph,\hat{\rm x}_{\eufm{norm}, \eufm{aux}}^\pph)}\\
    &j=n-h+1.\ \cdots, n-1\\
    \\ \dst{\rm y}_{\eufm{aux}, n-h+1}^{(h)}=&{\rm y}^\pph_{\eufm{norm}, \eufm{aux},n-h+1}\\
    \\ \dst{{\rm x}}_{\eufm{aux}, n-h+1}^\pph=&{{\rm x}}_{\eufm{norm}, \eufm{aux}, n-h+1}^\pph\\
    \qquad\qquad&+\f_{\eufm{norm}, \eufm{aux}}^\pph({\rm t}_{\eufm{norm}, \eufm{aux}}^\pph,{\rm z}_{\eufm{norm},\eufm{aux}, n-h}^\pph,\tilde{\rm y}_{\eufm{norm}, \eufm{aux}}^\pph,\hat{\rm x}_{\eufm{norm}, \eufm{aux}}^\pph)
\end{array}
\right. \eeqano where ${\rm F}_{\eufm{norm}, \eufm{aux}}^\pph$, ${\rm G}_{\eufm{norm}, \eufm{aux}}^\pph$ and ${\rm Z}_{\eufm{norm}, \eufm{aux}}^\pph$ are odd; ${\rm X}_{\eufm{norm}, \eufm{aux}}^\pph$, $\psi_{\eufm{norm}, \eufm{aux},j}^\pph$ and $\f_{\eufm{norm}, \eufm{aux}}^\pph$ are even under the change $$({\rm t}_{\eufm{norm}, \eufm{aux}}^\pph,{\rm z}_{\eufm{norm}, \eufm{aux}, n-h}^\pph)\to -({\rm t}_{\eufm{norm}, \eufm{aux}}^\pph,{\rm z}_{\eufm{norm}, \eufm{aux}, n-h}^\pph).$$ Then \equ{item1}-\equ{item3} follow.

\nl Now we proceed with proving the existence of $\phi^{n-h}_{\eufm{norm}, \eufm{aux}}$. We can choose, in \equ{abstract system},\equ{f***} and \equ{form of perturbation}, \begin{equation}\label{choices in averaging}\begin{split}&\n_i=2(h+1),\qquad \ell_i= h,\qquad m_i=3i,\quad i=1,\ \cdots,\ n-h-1=N\\\
&{\rm h}(p,q,I)=\sum_{i=n-h}^{n-1}{{{\rm h}_{\eufm{sec}}^i}}(\hat{\rm y}^\pph_i),\qquad f(p,q,I, \f, \eta,\xi)=\sum_{i=1}^{n-h-1} f^i(u_i, p,q,\f)\\
& f^i(u_i, p,q,\f):=f^{n-h-i}_{\eufm{norm},\eufm{int},\eufm{aux}, h-1}({\rm t}_{\eufm{aux},n-h- i}^\pph,\tilde{\rm y}_{\eufm{aux},n-h- i}^\pph,\hat{\rm x}_{\eufm{aux},n-h- i}^\pph)\\
&{\eufm Z}:={\eufm Z}_i:=\Big\{(k', k'', k''')\in \integer^{h}\times\integer^{h+1}\times\integer^{h+1}:\quad k'_{n-h+1}= \cdots=k'_{n-1}=0\\
&\qquad\qquad\ \ k''_{n-h+1}=\cdots=k''_{n}=0,\ k'''_{n-h}=\cdots=k'''_{n}=0,\ k''_1+\cdots+k''_{n-h}=0\Big\}\\
&{\eufm L}:=\Big\{(k', k'', k''')\in{\eufm Z}:\qquad k'_{n-h}=k''_{n-h}=0\Big\} \end{split}\end{equation} where we have re-named \beqano (p,q)&:=&(p_\eufm{aux}^\pph,q^\pph_\eufm{aux})=(p_{\eufm{aux},n-h}^\pph,\cdots, p^\pph_{\eufm{aux},n-1},\ q_{\eufm{aux},n-h}^\pph,\cdots, q^\pph_{\eufm{aux},n-1}, )\nonumber\\
I&:=&\big({\rm G}^\pph_{\eufm{aux}, n-h},\ \cdots,\ {\rm G}_{\eufm{aux}, n}^\pph,\L^\pph_{\eufm{aux}, n-h},\ \cdots,\ \L^\pph_{\eufm{aux}, n}\big)\nonumber\\
\f&:=&\big({\rm g}^\pph_{\eufm{aux}, n-h},\ \cdots,\ {\rm g}_{\eufm{aux}, n}^\pph,\ell^\pph_{\eufm{aux}, n-h},\ \cdots,\ \ell^\pph_{\eufm{aux}, n}\big)\nonumber\\
u_i&:=&(I,\ \eta^i,\ \xi^i),\qquad \eta:=\eta^1,\qquad \xi:=\xi^1\nonumber\\
\eeqano with \beqano \eta^{i}&:=&\big(\Theta_{\eufm{aux},n-h-i},\ \cdots,\ \Theta_{\eufm{aux},n-1},\ {\rm G}_{\eufm{aux}, n-h-i},\ \cdots,\ {\rm G}_{\eufm{aux}, n-h-1},\nonumber\\
&& \L_{\eufm{aux},n-h-i},\ \cdots,\ \L_{\eufm{aux},n-h-1}\big)\nonumber\\
\xi^{i}&:=&\big(\vartheta_{\eufm{aux},n-h-i},\ \cdots,\ \vartheta_{\eufm{aux},n-1},\ {\rm g}_{\eufm{aux}, n-h-i},\ \cdots,\ {\rm g}_{\eufm{aux}, n-h-1},\nonumber\\
&&\ell_{\eufm{aux},n-h-i},\ \cdots,\ \ell_{\eufm{aux},n-h-1}\big). \eeqano In order to verify that Proposition \ref{iterative lemma} can be applied, we have to check conditions \equ{non res} and \equ{new smallness cond}. Due to the choices of ${\eufm Z}$, ${\eufm L}$ and to the fact that only the function ${{{\rm h}_{\eufm{sec}}^{n-h}}}$ in the summand for ${{{\rm h}_{\eufm{sec}}}}$ in \equ{before normalization} depends on $(p_{\eufm{aux}, n-h}^\pph, q_{\eufm{aux}, n-h}^\pph, {\rm G}_{\eufm{aux}, n-h}^\pph)$, it is sufficient to check that condition \equ{non res} holds with $$\o=\o_{\eufm{ sec}}^{n-h},\quad (k',k)\in \integer^2\setminus\{0\},\quad K=\bar K.$$ But due to the choice of ${\eufm D}_{\eufm {int},\eufm{aux}}^h$ in \equ{normalized auxiliary domain}, we have that \equ{non res} is verified, with $$ {\eufm a}=\frac{(a_{n-h}^+)^2}{(a_{n-h+1}^-)^3\theta_{n-h}}\frac{\bar\g}{\bar K^{\bar\t}},\quad r={\eufm c}_{h}\frac{\theta_{n-h}\bar\g}{\bar K^{\bar\t+1}},\quad \varepsilon={\eufm c}_h\sqrt{\theta_{n-h}}. $$ It remains to check the inequalities in \equ{new smallness cond}. In view of the definition of $f^i$ following from the formulae \equ{norm int h-1}, \equ{fi norm int aux} and \equ{choices in averaging}, of the definition of ${f_{\eufm{norm}, h-1}^i}$ in \equ{split step h}, the definition of $\ovl{f_{\eufm{norm}, h-1}^i}$, the bound for $\widetilde{f_{\eufm{norm}, h-1}^i}$ in \equ{h-1 bound}, and first inequality in \equ{small secular}, we see that the former of the inequalities in \equ{new smallness cond} is satisfied with \beq{Ei}E_i=\frac{1}{{\eufm c}_h}\max\Big\{\frac{(a_{n-h-i}^+)^2}{(a_{n-h-i+1}^+)^3},\ \m \bar K\big(\frac{a_n^+}{a_1^-}\big)^{\frac{3}{2}}\frac{1}{a_{n-h-i+1}^-}\Big\}\qquad i=1,\ \cdots,\ n-h-1.\eeq In order to check that also the second inequality in \equ{new smallness cond} is satisfied, we previously note that the number $d_i$ in \equ{c and d} can be taken to be \beqano d_i={\eufm c}_{h}\min\big\{\frac{\theta_{n-h}\bar\g}{\bar K^{\bar\t+1}},\ \theta_{n-h-i}\big\},\qquad i=1,\ \cdots,\ n-h-1. \eeqano Inserting then the above values for $K$, ${\eufm a}$, $E_i$ and $d_i$ into the left hand side of the second inequality in \equ{new smallness cond}, we find that this can be bounded by

 \beqano \frac{1}{\tilde{\eufm c}_h}\max\Big\{&& \frac{\bar K^{2{\bar\t}+2}}{\bar\g^2}\frac{(a_{n-h-i}^+)^2}{(a_{n-h}^+)^2}\frac{(a_{n-h+1}^-)^3}{(a_{n-h-i+1}^-)^3},\quad \frac{\bar K^{\bar\t+1}}{\bar\g}\frac{(a_{n-h-i}^+)^2}{(a_{n-h}^+)^2}\frac{(a_{n-h+1}^-)^3}{(a_{n-h-i+1}^-)^3}\frac{\theta_{n-h}}{\theta_{n-h-i}} \nonumber\\
&&\frac{\bar K^{2{\bar\t}+2}}{\bar\g^2}\frac{ \m \bar K\big(\frac{a_n^+}{a_1^-}\big)^{\frac{3}{2}}}{(a_{n-h}^+)^2}\frac{(a_{n-h+1}^-)^3}{a_{n-h-i+1}^- },\quad \frac{\bar K^{\bar\t+1}}{\bar\g}\frac{ \m \bar K\big(\frac{a_n^+}{a_1^-}\big)^{\frac{3}{2}}}{(a_{n-h}^+)^2}\frac{(a_{n-h+1}^-)^3}{a_{n-h-i+1}^-}\frac{\theta_{n-h}}{\theta_{n-h-i}} \Big\} \eeqano Using \equ{Choice of parameters}, one easily finds that this quantity does not exceed

\beqa{small parameter} \frac{1}{\hat{\eufm c}_h}\max\Big\{ \m(\frac{a_n}{a_1})^{5}\frac{\bar K^{2\bar\t+2}}{\bar\g^2},\quad \frac{\bar K^{\bar\t+1}\sqrt\a}{\bar\g} \Big\}<1. \eeqa where $\hat{\eufm c}_h$ depends only on the ratio $a_n^-/a_n^+$ and the masses and the inequality follows from \equ{small secular}. This conclude the proof of this case. \qed

 \subsection{Construction of $\phi_{\eufm{norm}}^{n-1}$}\label{first normalization} 

 The arguments we have used in the previous section to construct $\phi_{\eufm{norm}}^{1}$, $\cdots$, $\phi_{\eufm{norm}}^{n-2}$ also fit for the case of $\phi_{\eufm{norm}}^{n-1}$, therefore we shall not repeat them. We only limit to remark that, for this case, Equations \equ{omega sec n-h}, \equ{choices in averaging}, \equ{Ei} and \equ{small parameter} have to be replaced with \beqano &&\o_{\eufm{ sec}}^{n-1}(\hat{\rm y}_{\eufm{aux}, n-1}^{(1)}):=\left\{
\begin{array}
    {lllll} \dst \partial_{\frac{(p_{\eufm{aux}, n-1}^\ppu)^2+(q_{\eufm{aux}, n-1}^\ppu)^2}{2}, {\rm G}_{\eufm{aux}, n-1}^\ppu, {\rm G}_{\eufm{aux}, n}^\ppu}\,{{{\rm h}_{\eufm{sec}}^{n-1}}}(\hat{\rm y}_{\eufm{aux}, n-1}^{(1)})\ &n\ge 3\\
    \\ \dst \partial_{{\rm G}_{\eufm{aux}, 2}^\ppu}\,{{{\rm h}_{\eufm{sec}}^{2}}}(\hat{\rm y}_{\eufm{aux}, 1}^{(1)}) &n=2,
\end{array}
\right. \nonumber\\
&&f^i=f^{n-i}_{\eufm{norm}, \eufm{int}, \eufm{aux}, 0}({\rm t}_{\eufm{aux}, i}^\ppu,\hat{\rm y}_{\eufm{aux}, i}^\ppu,\hat{\rm x}_{\eufm{aux}, i}^\ppu),\quad d_i={\eufm c}_{1}\min\big\{\frac{\theta_{n-1}\bar\g}{\bar K^{\bar\t+1}},\ \theta_{n-i-1}\big\}\nonumber\\
&&i=1,\ \cdots,\ n-1,\quad \theta_0:=\theta_1\nonumber\\
&&E_i=\frac{1}{\hat{\eufm c}_1}\left\{
\begin{array}
    {lll} \dst\max\Big\{\m \big(\frac{a_n^+}{a_1^-}\big)^{\frac{3}{2}}\frac{1}{a_{n}^-},\ \frac{(a_{n-1}^+)^3}{(a_n^-)^4}\Big\}\qquad &i=1\\
    \\ \dst\max\Big\{\m \bar K\big(\frac{a_n^+}{a_1^-}\big)^{\frac{3}{2}}\frac{1}{a_{n-i+1}^-},\ \frac{(a_{n-i}^+)^2}{(a_{n-i+1}^-)^3}\Big\}\qquad &n\ge 3,\ i=2,\cdots,\ n-1
\end{array}
\right.\nonumber\\
&&\frac{1}{\hat{\eufm c}_1}\max\Big\{ \m(\frac{a_n}{a_1})^{5}\frac{\bar K^{2\bar\t+2}}{\bar\g^2},\quad\frac{\bar K^{2(\bar\t+1)}\a}{\bar\g^2}\Big\}.  \eeqano\qed

\appendix \chapter{Computing the domain of holomorphy} \section{On the analyticity of the solution of Kepler equation}\label{Computing the domain of holomorphy} Here is a refinement of Proposition \ref{Kepler equation*}.
\begin{proposition}
    \label{Kepler equation} Let $\widehat e$ be as in \equ{hat e}. For any $0<\ovl e<\widehat e$ there exists $\ovl\eta=\ovl\eta(\ovl e)$ such that, for any $\ovl\eta<\eta<1$ and any ${\rm e}\in \complex$ with $|{\rm e}|\le \ovl e$, there exist two positive numbers $\bar\zeta=\bar\zeta(\eta, {\rm e})$, $\ovl\ell=\ovl\ell(\eta, \ovl e)$ such that the map \beq{Eq: Kepler map}\zeta\in \ovl\torus_{ \bar\zeta} \to K(\zeta,{\rm e}):=\zeta-{\rm e}\sin\zeta\eeq is injective, its image verifies $$K( \ovl\torus_{ \bar\zeta},{\rm e})\supset \ovl\torus_{\ovl\ell}\qquad \forall\ {\rm e}\in \complex:\ |{\rm e}|\le \ovl e.$$ The inverse function $$\ell\in\ovl\torus_{\ovl\ell} \to \zeta(\ell, {\rm e}):=K^{-1}(\ell, {\rm e})\in \ovl\torus_{ \bar\zeta_\eta({\rm e})}$$ verifies \beq{zeta is analytic}|1-{\rm e}\cos \zeta(\ell, {\rm e})|\ge 1-\eta\eeq Therefore, $ \zeta(\ell, {\rm e})$ is real-analytic for $\ell\in \ovl\torus_{\ovl\ell}$.
\end{proposition}
\nl The proof of Proposition \ref{Kepler equation} is elementary and goes along the same lines of \cite{leviCivita1904}. Therefore, we shall present it skipping some detail.
\begin{lemma}
    \label{bar l positive} Let $\widehat e$ be as in Proposition \ref{Kepler equation*}. For any $0<\ovl e<\widehat e$ there exists a unique $\ovl\eta=\ovl\eta(\ovl e)\in (\ovl e,1)$ such that $$\forall\ \eta\in[\ovl\eta,\ 1):\quad \ovl\ell_\eta(\ovl e):=\log\Big[\frac{\eta}{\ovl e}+\sqrt{1+\frac{\eta^2}{\ovl e^2}}\Big]-\sqrt{\eta^2+\ovl e^2}\ge 0,\qquad \ell_\eta(\ovl e)=0\quad \iff\quad \eta=\ovl\eta.$$
\end{lemma}
\begin{proof} By definition of $\widehat e$, and since the function $\r\in [0,1]\to \frac{\r\, e^{\sqrt{1+\r^2}}}{1+\sqrt{1+\r^2}}$ increases with $\r$, we have $$\dst\frac{\ovl e\, e^{\sqrt{1+\ovl e^2}}}{1+\sqrt{1+\ovl e^2}}<1.$$ Consider now the function $$\eta\in (0,1]\to g_\r(\eta):=\frac{\r\, e^{\sqrt{\eta^2+\r^2}}}{\eta+\sqrt{\eta^2+\r^2}}.$$ This function decreases with $\eta$ for any $\r\in (0,1]$. Since $$g_{\ovl e}(0)=e^{\ovl e}>1,\qquad g_{\ovl e}(1)=\frac{\ovl e\, e^{\sqrt{1+\ovl e^2}}}{1+\sqrt{1+\ovl e^2}}<1$$ we find a unique $\ovl\eta=\ovl\eta(\ovl e)\in [0,1]$ such that $$g_{\ovl e}(\eta)<1\qquad \forall\ \ovl\eta<\eta< 1,\quad g_{\ovl e}({\eta(\ovl e)})=1.$$ Since also $$g_{\ovl e}(\ovl e)=\frac{e^{\ovl e\sqrt2}}{1+\sqrt2}\ge \frac{e^{\sqrt2}}{1+\sqrt2}>1$$ we actually have $$\ovl e<\ovl\eta<1.$$\end{proof}

\vskip.in\textsc{Proof of  Proposition \ref{Kepler equation}.} We shall prove Proposition \ref{Kepler equation} with \beqa{zeta and ell} &&\bar\zeta(\eta, {\rm e}):=\log{\frac{\sqrt{\eta^2+{\rm e}_2^2}+\sqrt{\eta^2-{\rm e}_1^2}}{\sqrt{{\rm e}_1^2+{\rm e}_2^2}}}\nonumber\\
&&\ovl\ell(\eta, \ovl e):=\log\Big[\frac{\eta}{\ovl e}+\sqrt{1+\frac{\eta^2}{\ovl e^2}}\Big]-\sqrt{\eta^2+\ovl e^2} \eeqa where ${\rm e}={\rm e}_1+{\rm i}{\rm e}_2$. Observe that $\ovl\ell(\eta, \ovl e)>0$ by Lemma \ref{bar l positive}. Moreover, since $${\rm e}_1\le |{\rm e}|\le \ovl e<\ovl\eta<\eta$$ we have that $\bar\zeta(\eta, {\rm e})$ is well defined and positive\footnote{Actually, $\bar\zeta(\eta, {\rm e})$, as a function of $({\rm e}_1, {\rm e}_2)$, reaches its positive minimum $$\ovl\zeta_{\rm min}=\log\Big[\frac{\eta}{\ovl e}+\sqrt{1+\frac{\eta^2}{\ovl e^2}}\Big]>\log(1+\sqrt2)$$ for $({\rm e}_1, {\rm e}_2)=(0,\ovl e)$.}: $$\bar\zeta(\eta, {\rm e})\ge \log\frac{\eta}{\ovl e}>0.$$ 

\nl We split Equation \equ{Eq: Kepler map} into its real and imaginary part \beqano \arr{ K_1(\zeta_1,\zeta_2,{\rm e}_1, {\rm e}_2):=\zeta_1-({\rm e}_1\sin\zeta_1\cosh\zeta_2-{\rm e}_2\cos\zeta_1\sinh\zeta_2)=\ell_1\\
\\
K_2(\zeta_1,\zeta_2,{\rm e}_1, {\rm e}_2):=\zeta_2-({\rm e}_1\cos\zeta_1\sinh\zeta_2+{\rm e}_2\sin\zeta_1\cosh\zeta_2)=\ell_2 } \eeqano (with $\zeta=\zeta_1+{\rm i}\zeta_2$, $\ell=\ell_1+{\rm i}\ell_2$). The equation for the real part gives a unique solution $$\zeta_1={\rm Z}_1({\rm e}_1,{\rm e}_2, \zeta_2, \ell_1)$$ provided \beq{E<1}|{\rm e}_1|\le \eta,\qquad |\zeta_2|\le \ovl\zeta(\eta, {\rm e})\eeq since it reduces to an ordinary real Kepler equation $$\zeta_1-E_1({\rm e}_1, {\rm e}_2, \zeta_2)\sin(\zeta_1-\phi_1({\rm e}_1, {\rm e}_2, \zeta_2)=\ell_1\qquad {\rm if}\quad E_1({\rm e}_1, {\rm e}_2, \zeta_2)\ne 0$$ $$\zeta_1=\ell_1\qquad {\rm otherwise}$$ with \beqano &&E_1({\rm e}_1, {\rm e}_2, \zeta_2):=\sqrt{{\rm e}_1^2\cosh^2\zeta_2+{\rm e}_2^2\sinh^2\zeta_2}\nonumber\\
&&\phi_1({\rm e}_1, {\rm e}_2, \zeta_2):\qquad E_1\cos\phi_1={\rm e}_1\cosh\zeta_2,\quad E_1\sin\phi_1={\rm e}_2\sinh\zeta_2. \eeqano and, under condition \equ{E<1}, one has \beq{E1}E_1\le \eta<1.\eeq Observe that this solution ${\rm Z}_1({\rm e}_1,{\rm e}_2, \zeta_2, \ell_1)$ verifies \beq{Z1 is odd}{\rm Z}_1({\rm e}_1,{\rm e}_2, -\zeta_2, \ell_1)=-{\rm Z}_1({\rm e}_1,{\rm e}_2, \zeta_2, \ell_1)\qquad {\rm mod}\ 2\p.\eeq On the other hand, the function $$\zeta_2\to {\rm K}_2({\rm e}_1, {\rm e}_2, \zeta_2, \ell_1):=K_2({\rm Z}_1({\rm e}_1, {\rm e}_2, \zeta_2, \ell_1), \zeta_2, {\rm e}_1, {\rm e}_2 )$$ is strictly increasing, therefore, it maps the interval $ [-\ovl\zeta(\eta, {\rm e}), \ovl\zeta(\eta, {\rm e})]$, onto the interval $[-{\rm L}_2$ $(\eta$, ${\rm e}$, $\ell_1)$, ${\rm L}_2$ $(\eta$, ${\rm e}$, $\ell_1)]$, where ${\rm L}_2(\eta, {\rm e}, \ell_1):={\rm K}_2({\rm e}_1, {\rm e}_2,\ovl\zeta(\eta, {\rm e}), \ell_1)$ (note that ${\rm K}_2({\rm e}_1, {\rm e}_2,-\ovl\zeta(\eta, {\rm e}), \ell_1)=-{\rm K}_2({\rm e}_1, {\rm e}_2,\ovl\zeta(\eta, {\rm e}), \ell_1)$ because of \equ{Z1 is odd}). We have thus proved that the map \equ{Eq: Kepler map} maps bijectively the strip $\ovl\torus_{\ovl\zeta(\eta, {\rm e})}$ onto the set $$\ell={\ell}_1+{\rm i}\ell_2\in \complex:\quad \ell_1\in \torus,\quad \ell_2\in [-{\rm L}_2(\eta, {\rm e}, \ell_1), {\rm L}_2(\eta, {\rm e}, \ell_1)].$$ But the curve $$\ell_2={\rm L}_2(\eta, {\rm e}, \ell_1)\qquad \ell_1\in [0,2\p)$$ is concave, its minimum points are cusps, where ${\rm L}_2$ attains the value $${\rm L}_{2, \rm min}(\eta, {\rm e})=\ovl\zeta(\eta, {\rm e})-\sqrt{\eta^2-{\rm e}_1^2+{\rm e}_2^2}.$$ The minimum of this quantity while $|{\rm e}|\le \ovl e$ is just $\ovl \ell(\eta, \ovl e)$ in \equ{zeta and ell}. Inequality in \equ{zeta is analytic} follows from $$|1-{\rm e}\cos\zeta|\ge |\Re(1-{\rm e}\cos\zeta)|\ge 1-|\Re\big({\rm e}\cos\zeta\big)|$$ and (by \equ{E1}) $$|\Re\big({\rm e}\cos\zeta\big)|=|E_1({\rm e}_1, {\rm e}_2, \zeta_2)\cos(\zeta_1-\phi_1({\rm e}_1, {\rm e}_2, \zeta_2)|\le E_1\le \eta.$$
\qed
\section{Proof of Proposition \ref{prop: domain}}\label{Computing the domain of holomorphy***}

\nl Define \beqano &&\underline\d_j:=\sqrt{1-\overline e_j^2},\qquad \ovl\d_j:=\sqrt{1-\underline e_j^2}.\eeqano Assume \equ{initial bounds}, with \beqano &&{\rm A}:=(1-\s^2)\sqrt{\frac{1}{(1+{\s})^3(1+{\s}^2)^4}},\quad \cB:=\sqrt{\frac{1}{(1-\s^2)(1+\s)^3(1+\s^2)}}\nonumber\\
&&\ovl{\rm C}_i:=\arr{{\rm C}_1(\s)\ovl\d_i\quad i=1,\cdots, n-1\\
\ovl\d_n\quad i=n } ,\quad \underline{\rm C}_i:=\arr{{\rm C}_2(\s)\sqrt{\underline\d_i^2+2g(\s)^2\ovl\d_i^2}\quad i=1,\cdots, n-1\\
\sqrt{\underline\d_i^2+2g(\s)^2\ovl\d_n^2}\quad i=n }\nonumber\\
&&s=\s(1-\s) \eeqano where \beqa{AC1C2} &&{\rm C}_1(\s):=\sqrt{1-\s^2},\quad {\rm C}_2(\s):=\sqrt{\frac{(1+\s^2)^3}{(1-\s^2)^2}}\eeqa and $\s$, $g$ are chosen as follows: $g(\s')$ is a suitable positive function, depending at most on the ratios $\frac{\L_j^+}{\L_j^-}$, $\frac{{\rm G}_i^+}{{\rm G}_i^-}$, such\footnote{Since, for $j=1$, $\cdots$, $n$, $\|{\rm C}^\ppj_{\cP}\|^2$ depends only on $\chi_{j-1}$, $\chi_j$, $\Theta_j$ and $\vartheta_j$ as in \equ{CP} and all such coordinates, together also with $\L_j$, have their anomalies bounded by $\s'$, we can always find such a function $g(\s')$.} that \beq{g(eps)}g(\s')\to 0\qquad{\rm as}\quad \s'\to 0,\quad {\rm and}\quad |\sin\arg\frac{\|{\rm C}^\ppj_{\cP}\|^2}{\L_j^2}|\le g(\s'),\quad j=1,\cdots, n, \eeq provided $$\max\Big\{|\arg (\L_i)|, |\arg(\chi_j)|, |\arg(\Theta_j)|,\ |\arg(\vartheta_j)|\Big\}\le \s'$$ while $\s$ is so small that, if $\ovl\ell_1$, $\cdots$, $\ovl\ell_n$ are as in Proposition \ref{Kepler equation*}, with $\ovl e$ replaced by $\ovl e_1$, $\cdots$, $\ovl e_n$, then $$\s\le \min\Big\{\frac{3}{4},\ \ovl\ell_1,\ \cdots,\ \ovl\ell_n\Big\}$$ and the following inequality is satisfied \[\frac{{\rm C}_1(\s)}{{\rm C}_2(\s)}\frac{\ovl\d_j}{\sqrt{\underline\d_j^2+\sqrt2g(\s)\ovl\d_j}}>1\qquad \forall\ i=1,\cdots,n.\]

\nl Note that this inequality is satisfied for $\s$ suitably small, since, by definition, $$ \ovl\d_j>\underline\d_j,\quad {\rm C}_1(\s')\uparrow 1,\quad {\rm C}_2(\s')\downarrow 1,\quad g(\s')\downarrow 0 \qquad {\rm as} \quad \s'\to 0.$$

\nl Definitions and assumptions in \equ{initial bounds} imply, since $\s(1-\s)<\s$, \beqa{list} &&(1-\s){\rm G}_n^-<|\chi_{i}|< {\rm G}_n^+(1+\s)\nonumber\\
&&|\tan\arg(\chi_{i-1}-\chi_{i})|\le\frac{\max|\Im(\chi_{i-1}-\chi_{i})|}{\min|\Re(\chi_{i-1}-\chi_{i})|}\le\frac{\theta_i}{{\rm G}_i^-}\le\s\le1\nonumber\\
&&|\arg\chi_{i}|\le |\arg\chi_{n-1}|+\sum_{j=i+1}^{n-1}|\sin^{-1}\frac{|\chi_{j-1}-\chi_{j}|}{|\chi_{j}-\chi_{j+1}|}|\le\s\le \frac{\p}{3} \eeqa The previous inequalities imply that, firstly $$|\frac{\Theta_j}{\chi_{j-1}}|\le \frac{\s(1-\s){\rm G}_n^-}{(1-\s){\rm G}_n^-}\le \s$$ and, similarly, $$|\frac{\Theta_j}{\chi_{j}}|\le \s$$ therefore, the inequality for ${\rm i}_j$, $\iota_i$ is \equ{eccentricities} follows. Secondly, the definitions of $\Theta^+_{i}$, $\vartheta^+_{i}$ imply that conditions \equ{ass L2} are met and hence Lemma \ref{L2} applies. By the thesis \equ{bounds on C}, we have\footnote{Beware that, if $z=(z_1,z_2,z_3)\in \complex^3$, we denote $$\|z\|^2:=z_1^2+z_2^2+z_3^3.$$ For a given $z\in \complex$, the symbol $|z|$ denotes the usual modulus of $z\in \complex$: $$|z|:=\sqrt{(\Re z)^2+(\Im z)^2}.$$}, for $j=1$, $\cdots$, $n-1$,

\beqa{upper C} \big|\|{\rm C}_\cP^\ppj\|^2\big|&\le&\frac{\big|\chi_{j-1}-\chi_{j}\big|^2}{1-\s^2}+ (1+\s)(1+\s^2)|\chi_{j-1}||\chi_{j}||\vartheta_{j}-\p|^2\nonumber\\
&\le&\frac{({\rm G}_i^+)^2}{\ovl{\rm C}_j^2}+ \frac{({\rm G}_n^+)^2}{\ovl{\rm C}_j^2{\rm B}^2}|\vartheta_{j}-\p|^2\nonumber\\
&\le&\ovl\d_j^2(\L^-_j)^2. \eeqa For $j=n$, $$\big|\|{\rm C}_\cP^\ppj\|^2\big|=\big|\chi_{n-1}\big|^2\le ({\rm G}_n^+)^2<\ovl\d_{n}^2(\L^-_n)^2.$$ We suddenly have the left bound in \equ{eccentricities}: $$1-|e_{i,\cP}^2|\le|1-e_{i,\cP}^2|=|\frac{\big|\|{\rm C}^\ppi_{\cP}\|^2\big|}{\L_i^2}|\le \ovl\d_i^2=1-\underline e_{i}^2,$$ for $i=1$, $\cdots$, $n$. Now we check the right bound. To this end, previously check the following inequality \beq{lower bound}\big||\chi_{j-1}|-|\chi_{j}|\big|\ge \frac{1-\s^2}{1+\s^2}{\rm G}_j^-.\eeq

\nl Because of the second inequality in \equ{list}, $$|\arg\big[(\chi_{j-1}-\chi_{j})(\ovl\chi_{m-1}-\ovl\chi_{m})\big]|\le 2\tan^{-1}\s.$$ Then we have $$\Re\big[(\chi_{j-1}-\chi_{j})(\ovl\chi_{m-1}-\ovl\chi_{m})\big]\ge\frac{1-\s^2}{1+\s^2}|\chi_{j-1}-\chi_{j}||\ovl\chi_{m-1}-\ovl\chi_{m}|.$$ Taking the sum for $m=j+1$, $\cdots$, $n$, gives \beqano \Re(\chi_{j-1}-\chi_{j})\ovl\chi_{j}&\ge&\frac{1-\s^2}{1+\s^2}|\chi_{j-1}-\chi_{j}|\sum_{m=j+1}^{n}|\ovl\chi_{m-1}-\ovl\chi_{m}|\ge\frac{1-s^2}{1+s^2}|\chi_{j-1}-\chi_{j}||\ovl\chi_{j}| \nonumber\\
&\ge& \frac{1-\s^2}{1+\s^2}{\rm G}_j^-|\ovl\chi_{j}| \eeqano So, Lemma \ref{L1} with $$A=\chi_{j-1},\quad B=\chi_{j},\quad \D={\rm G}_j^-,\quad a=\frac{1-\s^2}{1+\s^2}$$ gives \equ{lower bound}. Then the thesis \equ{below} of Lemma \ref{L2} and the definition of $\vartheta_{j}$ provide, for $j=1$, $\cdots$, $n-1$, \beq{lower C}\big|\|{\rm C}^\ppj_{\cP}\|^2\big|\ge\frac{1}{\cA^2\underline\cC_j^2}\Big[\cA^2({\rm G}_j^-)^2-({\rm G}_n^+)^2|\vartheta_{j}-\p|^2\ge \Big]\ge(\underline\d_j^2+\sqrt2g(\s)\ovl\d_i)(\L_j^+)^2\eeq where $g(\s)$ is as in \equ{g(eps)}. Again, this inequality is implied by the definition of $\vartheta_{j}^+$ in \equ{initial bounds} and the ones of $\cA$ and ${\rm C}_2$ in \equ{AC1C2}. By \equ{g(eps)}, \equ{upper C} and \equ{lower C}, for $j=1,\cdots, n$, we have \beqa{ej} |e_{j,\cP}|^2&=&\sqrt{ (1-\Re\frac{\|{\rm C}^\ppj_{\cP}\|^2}{\L_j^2})^2+(\Im\frac{\|{\rm C}^\ppj_{\cP}\|^2}{\L_j^2})^2} \nonumber\\
&\le&\sqrt{\Big(1-\big|\frac{\|{\rm C}^\ppj_{\cP}\|^2}{\L_j^2}\big|\Big)^2+2|\Im\frac{\|{\rm C}^\ppj_{\cP}\|^2}{\L_j^2}}|\nonumber\\
&\le&\sqrt{\Big(1-\underline\d_j^2-\sqrt2g(\s)\ovl\d_j\Big)^2+2\overline\d_j^2g(\s)^2}\le1-\underline\d_j^2=\ovl e_{j}^2. \eeqa For $j=n$, $$\big|\|{\rm C}^\ppn_\cP\|^2\big|=|\chi_{n-1}|^2\ge(\underline\d_n^2+\sqrt2g(\s)\ovl\d_n)(\L_n^+)^2$$ again implies \equ{ej} with $j=n$.

\nl The proof of the inequality on the right in \equ{eccentricities} proceeds in a similar way. Indeed, starting with \beqano |d_{i, \cP}|^2=\Big|\|x_\cP^{(i+1)}\|^2-2x_\cP^\ppi\cdot x_\cP^{(i+1)}+\|x_\cP^\ppi\|^2\Big|\ge\Big|\|x_\cP^{(i+1)}\|^2\Big|-2\Big|x_\cP^\ppi\cdot x_\cP^{(i+1)}\Big|-\Big|\|x_\cP^\ppi\|^2\Big| \eeqano and using (as it follows from Proposition \ref{Kepler equation}) $$\Big|\|x_\cP^{(i+1)}\|^2\Big|=|a_{i+1}^2(1-e_{i+1,\cP}\cos\zeta_{i+1})^2|\ge (1-\eta_{i+1})^2(a^-_{i+1})^2$$ and {analogous} arguments as above to evaluate $\Big|x_\cP^\ppi\cdot x_\cP^{(i+1)}\Big|$ and $\Big|\|x_\cP^\ppi\|^2\Big|$, one easily finds the ansatz. \qed

\vskip.1in\noi\paragraph{\bf Estimates}
\begin{lemma}
    \label{L2} Fix a number $\s>0$. Assume that, for $1\le j\le n-1$, \begin{equation}\label{ass L2} \begin{split}\Re\ovl\chi_{j}(\chi_{j-1}-\chi_{j})>0,\quad |\Theta_{j}|\le \s\min\{|\chi_{j-1}|,\ |\chi_{j}|\},\quad |\Im(\vartheta_{j}-\p)|\le \log(1+\s). \end{split}\end{equation} Then \beqa{bounds on C} &&\big|\|{\rm C}^\ppj_{\cP}\|^2\big|\le \frac{\big|\chi_{j-1}-\chi_{j}\big|^2}{1-\s^2}+ (1+\s)(1+\s^2)|\chi_{j-1}||\chi_{j}||\vartheta_{j}-\p|^2\\
    &&\label{below}\big|\|{\rm C}^\ppj_{\cP}\|^2\big|\ge \frac{\big||\chi_{j-1}|-|\chi_{j}|\big|^2}{1+\s^2}-(1+\s)(1+\s^2)|\chi_{j-1}||\chi_{j}||\vartheta_{j}-\p|^2 \eeqa
\end{lemma}
\begin{proof} We use the formula \equ{CPjsq}. By Taylor's, given $a$, $b$, $z\in \complex$, with $|z|\le \s \min_{t\in[0,1]}|a+t(b-a)|$ \beqano \Big|\sqrt{b^2-z^2}-\sqrt{a^2-z^2}\Big|&=&\Big|\int_0^1\frac{d}{dt}\sqrt{\big(a+t(b-a)\big)^2-z^2}dt\Big|\nonumber\\
&=&\Big|(b-a)\int_0^1\frac{a+t(b-a)}{\sqrt{\big(a+t(b-a)\big)^2-z^2}}dt\Big|\nonumber\\
&\le&|b-a| \int_0^1\frac{|a+t(b-a)|}{\sqrt{\big|a+t(b-a)\big|^2-|z|^2}}dt\nonumber\\
&\le&\frac{|b-a|}{\sqrt{1-\s^2}} \eeqano We use this formula with $b:=\chi_{j-1}$, $a:=\chi_{j}$, $z:=\Theta_{j}$, with the observation that, for $\Re\ovl\chi_{j}(\chi_{j-1}-\chi_{j})>0$, the function $$t\in [0,1]\to\big|\chi_{j}+t(\chi_{j-1}-\chi_{j})\big|^2=|\chi_{j}|^2+2t\Re\ovl\chi_{j}(\chi_{j-1}-\chi_{j})+t^2|\chi_{j-1}-\chi_{j}|^2$$ reaches its minimum, given by $\min\{|\chi_{j-1}|^2,\ |\chi_{j}|^2\}$, for $t=0$ or $t=1$. Developing also the function $w\in \complex\to \cos w$ around $w=\p$, with $\varrho:=w-\p=\varrho_1+{\rm i}\varrho_2$ and $|\r_2|\le \log(1+\s)$ \beqano \big|\cos w+1\big|&=&\big|\int_0^1(1-t)\frac{d^2}{dt^2}\cos(\p+t(w-\p))\big|=\frac{1}{2}|\varrho|^2\sup_{|\varrho'|\le \varrho}|\cos (\p+\varrho')|\nonumber\\
&\le& \frac{1}{2}|\varrho|^2e^{|\varrho_2|}\le \frac{1}{2}|\varrho|^2(1+\s)\eeqano and using again the second inequality in \equ{ass L2}, then inequality in \equ{bounds on C} follows. The inequality in \equ{below} is obtained via the second inequality in \equ{ass L2} and \beqano \big|\sqrt{\chi_{j}^2-\Theta_{j}^2}-\sqrt{\chi_{j-1}^2-\Theta_{j}^2}\big|&=&\frac{\big|\chi_{j-1}^2-\chi_{j}^2\big|}{\big|\sqrt{\chi_{j}^2-\Theta_{j}^2}+\sqrt{\chi_{j-1}^2-\Theta_{j}^2}\big|}\nonumber\\
&\ge&\frac{\big||\chi_{j-1}|^2-|\chi_{j}|^2\big|}{\big|\sqrt{\chi_{j}^2-\Theta_{j}^2}+\sqrt{\chi_{j-1}^2-\Theta_{j}^2}\big|}\nonumber\\
&\ge& \frac{\big||\chi_{j-1}|-|\chi_{j}|\big|}{ \sqrt{1+\s^2} }.  \eeqano
\end{proof}
\begin{lemma}
    \label{L1} If $A$, $B\in \complex$ and $a$, $\D\in \real_+$ verify $|A-B|\ge \D$ and $\Re\ovl B(A-B)\ge a|B|\D$, where $0<a<1$, then $\big||A|-|B|\big|> a\D$.
\end{lemma}
\begin{proof} Let $D:=A-B$. Then $\big||A|-|B|\big|=\big||B+D|-|B|\big|\le a\D$ implies $$ |B|^2+|D|^2+2\Re\ovl BD=|B+D|^2\le (|B|+a\D)^2=|B|^2+a^2(\D)^2+2a|B|\D. $$ This contradicts assumptions $|D|\ge \D>a\D$ and $\Re\ovl BD\ge a|B|\D$. \end{proof}

\chapter{Proof of Lemma 3.2}\label{Proof of Lemma {integrabilityyyy}} In this section, we prove the formulae \equ{fn-1n} and \equ{ovl f} given in Lemma \ref{integrabilityyyy}. 

\nl We recall the following result
\begin{proposition}
    [\cite{pinzari13}] Let ${\eufm X}={\eufm X}_1\times\cdots\times{\eufm X}_n\subset \real^5\times \cdots\times \real^5$ and let

$$(\ell_k, {\rm X}_k)\in \torus^1\times {\eufm X}_k\to(y_\phi^\ppk( \ell_k,{\rm X}_k), x_\phi^\ppk( \ell_k,{\rm X}_k))\in \real^3\times \real^3\qquad k=1,\ \cdots,\ n$$ be mappings such that, for $1\le i<j\le n$
    \begin{enumerate}
        \item[{\bf (A)}] the map \[\phi_{ij}:\quad ( \ell_i,\ell_j,{\rm X}_i, {\rm X}_j)\to (y_\phi^\ppi, y_\phi^\ppd, x_\phi^\ppj, x_\phi^\ppd)\] is symplectomorphism of $\torus^2\times {\eufm X}_i\times {\eufm X}_j$ into $\real^{12}$. \item[{\bf (B)}] The map $( \ell_j,{\rm X}_j)\to(y_\phi^\ppd( \ell_j,{\rm X}_j), x_\phi^\ppd( \ell_j,{\rm X}_j))$ verifies \[ \frac{\|y_\phi^\ppd( \ell_j,{\rm X}_j)\|^2}{2m_j}-\frac{{\eufm m}_j\eufm{M}_j}{\|x_\phi^\ppd( \ell_j,{\rm X}_j)\|}=-\frac{\eufm{m}_j^3\eufm{M}_j^2}{2\L_j^2}; \] where $\L_j$ is the variable conjugated to $\ell_j$ in this symplectomorphism.
    \end{enumerate}
    Then the function \beqano {\rm P}^\ppi(\ell_i,{\rm X})&:=&-\frac{1}{2\p}\int_{\torus}d\ell_j\nonumber\\
    &&\quad\frac{3(x_\phi^\ppi(\ell_i,{\rm X}_i)\cdot x_\phi^\ppj(\ell_j,{\rm X}_j))^2-\|x_\phi^\ppi(\ell_i,{\rm X}_i)\|^2\|x_\phi^\ppj(\ell_j,{\rm X}_j)\|^2}{2\|x_\phi^\ppj(\ell_j,{\rm X}_j)\|^5} \eeqano is given by \beq{good formula}{\rm P}^\ppi=\frac{{\eufm M}_j\eufm{m}_j^2}{4}\frac{3( x_\phi^\ppi\cdot{\rm C}_\phi^\ppj)^2-\|x_\phi^\ppi\|^2\|{\rm C}_\phi^\ppj\|^2}{\|{\rm C}_\phi^\ppj\|^4}\frac{1}{2\p}\int_{\torus}\frac{d\ell_j}{\|x^\ppj_\phi\|^2}.\eeq with ${\rm C}_\phi^\ppj( {\rm X}):=x_\phi^\ppj(\ell_j,{\rm X})\times y_\phi^\ppj(\ell_j,{\rm X})$.
\end{proposition}
\nl Even though the $(i,j)$ projections of the $\cP$-map do not verify assumption (A), one has
\begin{corollary}
    \label{extension to P} The formula \equ{good formula} applies also to the $\cP$-map, or, more in general, to any Kepler map $\cK$ related to {the} map $\cDel$ in Definition \ref{Def: Delaunay} via $${\rm X}_{\cDel}={\eufm F}({\rm X}).$$
\end{corollary}
\begin{proof} $\cDel$ verifies (A) and (B). \end{proof}

\nl In particular, we have an expression for the second-order term of the doubly averaged Newtonian potential \beqano \ovl{f_\cK^{ij}}^\ppd&&:=-\frac{ m_{i}m_j}{(2\p)^2}\int_{\torus^2}d\ell_{i}d\ell_j\nonumber\\
&&\quad\frac{3(x_\cK^{(i)}(\ell_{i},{\rm X}_\cK)\cdot x_\cK^\ppj(\ell_j,{\rm X}_\cK))^2-\|x^{(i)}_\cK(\ell_{i},{\rm X}_\cK)\|^2\|x^\ppj_\cK(\ell_j,{\rm X}_\cK)\|^2}{2\|x_\cK^\ppj(\ell_j,{\rm X}_\cK)\|^5}. \eeqano
\begin{corollary}
    \label{new prop**} For any $\cK$ as in Corollary \ref{extension to P}, \beqa{f12****} {\ovl {f_{\cK}^{i j}}}^{(2)}&=&m_{i} m_{j} \frac{a_{i}^2}{4a_{j}^3}\frac{\L_{j}^3}{\|{\rm C}^{(j)}_{\cK}\|^5}\Big[ -\big(\frac{5}{2}-\frac{3}{2}\frac{\|{\rm C}^{(i)}_{\cK}\|^2}{\L_{i}^2}\big)\|{\rm C}^{(j)}_{\cK}\|^2 \nonumber\\
    &+&\frac{3}{2}\big(5-4\frac{\|{\rm C}^{({i})}_{\cK}\|^2}{\L_i^2}\big)(P_\cK^{(i)}\cdot {\rm C}_{\cK}^{(j)})^2+\frac{3}{2} \frac{\|{\rm C}^{({i})}_{\cK}\|^2}{\L_{i}^2}(Q_\cK^{(i)}\cdot {\rm C}_{\cK}^{(j)})^2 \Big] \eeqa
\end{corollary}
\begin{proof} Lemma \ref{extension to P} implies that \beqano \ovl{f^{ij}_\cK}^\ppd&=& m_i m_j\frac{M_jm_j^2}{4}\frac{\frac{1}{2\p}\int_{\torus}\Big(3( x_\cK^\ppi\cdot{\rm C}_\cK^\ppj)^2-\|x_\cK^\ppi\|^2\|{\rm C}_\cK^\ppj\|^2\Big)d\ell_i}{\|{\rm C}_\cK^\ppj\|^4}\nonumber\\
&&\times\frac{1}{2\p}\int_{\torus}\frac{d\ell_j}{\|x_\cK^\ppj\|^2}.\eeqano

\nl By \equ{xjyj} \beqano x_{\cK}^{(i)}\cdot {\rm C}_{\cK}^{(j)}&=&\big({\rm a}_{i,\cK}P^{(i)}_{\cK}+{\rm b}_{i,\cK}Q^{(i)}_{\cK}\big)\cdot {\rm C}^\ppj_{\cK} \nonumber\\
&=&{\rm a}_{i,\cK}P^{(i)}_{\cK}\cdot{\rm C}^\ppj_{\cK}+{\rm b}_{i,\cK}Q^{(i)}_{\cK}\cdot{\rm C}^\ppj_{\cK} \eeqano Therefore, squaring, $\ell_i$-averaging and using \beqano \frac{1}{2\p}\int_\torus({\rm a}_{i,\cK})^2d\ell_{i}&=&\frac{a_{i}^2}{2}\big(5-4\frac{\| {\rm C}^{(i)}_{\cK}\|^2}{\L_{i}^2}\big)\nonumber\\
\frac{1}{2\p}\int_\torus({\rm b}_{i,\cK})^2d\ell_{i}&=&\frac{a_{i}^2}{2}\frac{\| {\rm C}^{(i)}_{\cK}\|^2}{\L_{i}^2}\nonumber\\
\frac{1}{2\p}\int_\torus{\rm a}_{i,\cK}{\rm b}_{i,\cK}d\ell_{i}&=&0 \eeqano we obtain \beqano \frac{1}{2\p}\int_\torus(x_{\cK}^{(i)}\cdot {\rm C}_{\cK}^{(j)})^2d\ell_{i}&=&\frac{a_{i}^2}{2} (5-4\frac{\| {\rm C}^{(i)}_{\cK}\|^2}{\L_{i}^2})(P_{\cK}^{(i)}\cdot {\rm C}_{\cK}^{(j)})^2\nonumber\\
&+&\frac{a_{i}^2}{2}\frac{\| {\rm C}^{(i)}_{\cK}\|^2}{\L_{i}^2}(Q_{\cK}^{(i)}\cdot {\rm C}_{\cK}^{(j)})^2. \eeqano Using finally \beqano \frac{1}{2\p}\int_{\torus}{\|x^{(i)}_{\cK}\|^2}{d\ell_{i}} = a_{i}^2\big(\frac{5}{2}-\frac{3}{2}\dst\frac{\|{\rm C}^{(i)}_{\cK}\|^2}{\L_{i}^2}\big)\label{first},\quad \frac{1}{2\p}\int_{\torus}\frac{d\ell_{j}}{\|x^{(j)}_{\cK}\|^2}=\frac{1}{a_{j}^2}\frac{\L_{j}}{\|{\rm C}^{(j)}_{\cK}\|} \eeqano we obtain \equ{f12****}.\end{proof}

\nl Now we may proceed with proving the formulae in \equ{fn-1n} and \equ{ovl f}.

\vskip.in\textsc{Proof of  of \equ{fn-1n}.} We apply Corollary \ref{new prop**} with $\cK=\cP$, $i=n-1$, $j=n$. Using $\|{\rm C}_\cP^\ppn\|=\chi_{n-1}$ (see \equ{CP}), ${\rm C}^\ppn_{\cP}={\rm S}^\ppn_{\cP}$ and Eq. \equ{a e b}, Proposition \ref{inversion}, and Remark \ref{C perp to P}, we have \beqano P_{\cP}^{(n-1)}\cdot {\rm S}_{\cP}^{(n)}&=&\Theta_{n-1}\nonumber\\
Q_{\cP}^{(n-1)}\cdot {\rm S}_{\cP}^{(n)}&=&\frac{1}{\| {\rm C}^{(n-1)}_{\cP}\|} \big(({\rm S}^{(n-1)}_{\cP}-{\rm S}^{(n)}_{\cP})\times P_{\cP}^{(n-1)}\big)\cdot {\rm S}^\ppn_{\cP}\nonumber\\
&=&\frac{1}{\| {\rm C}^{(n-1)}_{\cP}\|} {\rm S}^{(n-1)}_{\cP}\times {\rm P}^{(n-1)}_{\cP}\cdot{\rm S}^{(n)}_{\cP} \nonumber\\
&=&\frac{1}{\| {\rm C}^{(n-1)}_{\cP}\|}\sqrt{(\chi_{n-1}^2-\Theta_{n-1}^2)(\chi_{n-2}^2-\Theta_{n-1}^2)}\sin{\vartheta_{n-1}}. \eeqano
\qed

\vskip.in\textsc{Proof of \equ{ovl f}.} By Corollary \ref{new prop**} with $\cK=\cP$, $j=i+1$, we find, for ${\ovl {f_{\cP}^{i,i+1}}}^{(2)}$ an expression as in \equ{f12****}, replacing $(n-1, n)$ with $(i, i+1)$.

 \begin{equation}\label{new products}\begin{split} &P_{\cP}^{(i)}\cdot {\rm C}^{(i+1)}_{\cP}=P_{\cP}^{(i)}\cdot ({\rm S}^{(i+1)}_{\cP}-{\rm S}^{(i+2)}_{\cP}) =\Theta_{i}-P_{\cP}^{(i)}\cdot {\rm S}^{(i+2)}_{\cP}\\
 &Q_{\cP}^{(i)}\cdot {\rm C}^{(i+1)}_{\cP}=Q_{\cP}^{(i)}\cdot ({\rm S}^{(i+1)}_{\cP}-{\rm S}^{(i+2)}_{\cP}) =\frac{1}{\| {\rm C}^{(i)}_{\cP}\|} ( \sqrt{(\chi_{i}^2-\Theta_{i}^2)(\chi_{i-1}^2-\Theta_{i}^2)}\sin{\vartheta_{i}}\\
&-{\rm S}^{(i)}_{\cP}\times {\rm P}^{(i)}_{\cP}\cdot{\rm S}^{(i+2)}_{\cP}-{\rm P}^{(i)}_{\cP}\times{\rm S}^{(i+1)}_{\cP}\cdot{\rm S}^{(i+2)}_{\cP} ).\end{split}\end{equation}

\nl Now, when $(\Theta_{i+1}, \vartheta_{i+1})=(0,\p)$, $\|{\rm C}^{(i+1)}_\cP\|$ reduces to $$\|{\rm C}^{(i+1)}_\cP\|=\chi_{i}-\chi_{i+1},$$ (provided $\arg(\chi_i-\chi_{i+1})\in (-\frac{\p}{2},\ \frac{\p}{2}]$ mod $2\p$) and ${\rm S}^{(i+2)}_{\cP}\parallel {\rm S}^{(i+1)}_{\cP}$, so $${\rm S}^{(i+2)}_{\cP}=\frac{\chi_{i+1}}{\chi_i}{\rm S}^{(i+1)}_{\cP}$$ and hence, the extra-terms in \equ{new products} reduce to \beqano && P_{\cP}^{i)}\cdot {\rm S}^{(i+2)}_{\cP}= \Theta_i\frac{\chi_{i+1}}{\chi_i} \nonumber\\
&& {\rm S}^{(i)}_{\cP}\times {\rm P}^{(i)}_{\cP}\cdot{\rm S}^{(i+2)}_{\cP}=\frac{\chi_{i+1}}{\chi_i}\sqrt{\chi_{i-1}^2-\Theta_i^2}\sqrt{\chi_{i}^2-\Theta_i^2} \sin\vartheta_i\nonumber\\
&& {\rm P}^{(i)}_{\cP}\times{\rm S}^{(i+1)}_{\cP}\cdot{\rm S}^{(i+2)}_{\cP}=0.\eeqano Then \equ{ovl f} readily follows. \qed

 \chapter{Checking the non-degeneracy condition}\label{Checking the non-degeneracy assumption} In this section we prove statement {\bf 4} of Proposition \ref{claim}.

\nl Due to the form of ${\rm h}_{\eufm{ sec}}$ in \equ{3split}-\equ{split step h} and to the bound for $\widetilde{{\rm h}_{\eufm{sec}, h}^i}$ in \equ{h-1 bound}, it is sufficient to prove that the maps $$\zeta^\pph_i\to \ovl{\o^i_{\eufm{ sec}}}:=\partial_{\zeta^\pph_i}\ovl{{\rm h}^i_{\eufm{ sec}}}(\zeta^\pph_i, \L^\pph_{n-h}, \L^\pph_{n-h+1})$$ in \equ{split step h}, where $$\zeta^\pph_i=\left\{
\begin{array}
    {llll} \dst\big(\frac{(p_{1}^\pph)^2+(q_{1}^\pph)^2}{2},\ \chi^\pph_{1}\big)\quad & i=1\ &\&\ n= 2\\
    \\ \dst\big(\frac{(p_{n-1}^\pph)^2+(q_{n-1}^\pph)^2}{2},\ \chi^\pph_{n-2},\ \chi^\pph_{n-1}\big)\quad & i=n-1\ &\&\ n\ge 3\\
    \\ \dst\big(\frac{(p_{i}^\pph)^2+(q_{i}^\pph)^2}{2},\ \chi^\pph_{i-1}\big)\quad & i=2,\ \cdots,\ n-2\ &\&\ n\ge 4\\
    \\ \dst\frac{(p_{1}^\pph)^2+(q_{1}^\pph)^2}{2}\quad & i=1\ &\&\ n\ge 3
\end{array}
\right. $$ are diffeomorphisms, with non-vanishing Hessian matrices. We shall do this verifications for just one of the cases above, and we choose the second case in the list, $i=n-1$, for $n\ge 3$. The explicit expression of $\ovl{{\rm h}^{n-1}_{\eufm{ sec}}}$ is given in \equ{Nn}-\equ{coefficients}. We neglect the coefficient $\cA_{n-1}$ (which does not depend on $\zeta^\pph_{n-1}$) and we denote $$\widehat{{\rm h}_{\eufm{sec}}^{n-1}}={{\rm E}_{n-1}}+{\Omega_{n-1}}\frac{p_{n-1}^2+q_{n-1}^2}{2}+{\t_{n-1}}(\frac{p_{n-1}^2+q_{n-1}^2}{2})^2+{\rm O}(p_{n-1},q_{n-1})^6\Big]$$ the function $\ovl{{\rm h}_{\eufm{sec}}^{n-1}}$ thus rescaled, and $\widehat{\o_{\eufm{int}}^{n-1}}$ its gradient with respect to $(\frac{(p_{n-1}^\pph)^2+(q_{n-1}^\pph)^2}{2}$, $\chi_{n-2}$, $\chi_{n-1})$. A perturbative argument shows that, under the choices of Corollary \ref{cor: parameters}, the frequency-map with respect to $(\chi_{n-2},\chi_{n-1})$ associated to $${\rm E}_{n-1}=-\frac{\L_n^3}{2\chi_{n-1}^3} \Big( 5-3\frac{(\chi_{n-2}-\chi_{{n-1}})^2}{\L_{{n-1}}^2} \Big)$$ is an injection of its domain and hence, by another perturbative argument, so is the gradient of $\widehat{{\rm h}_{\eufm{sec}}^{n-1}}$ with respect to the same coordinates, for any fixed value of $\frac{p_{n-1}^2+q_{n-1}^2}{2}$. On the other hand, since $\t_{n-1}$ does not vanish under the same assumptions of Corollary \ref{cor: parameters}, $\widehat{\o_{\eufm{int}}^{n-1}}$ is an injection. The computation shows that the Jacobian of $\widehat{\o_{\eufm{int}}^{n-1}}$ does not vanish. \qed

 \chapter{Some results from perturbation theory} \section{A multi-scale normal form theorem}\label{normal form theory}

 The purpose of this section is to present a normal form result which takes into account different scale lengths. It is a particularization of \cite[Normal Form Lemma, p. 192]{poschel93} and uses the same techniques of that paper. 

\nl Following \cite{poschel93}, the notations are as follows.
\begin{enumerate}
    \item[\tiny\textbullet]If $A\subset \real^{\n}$ is open and connected, $\torus:=\real/(2\p\integer)$ is the usual flat torus, $r$, $s$ are positive numbers, we denote as $A_r:=\bigcup_{x\in A}\Big\{z\in \complex^{\n}:\quad z\in B^\n_r(x)\Big\}$ the complex $r$-neighborhood of $A$. $\torus^\n_s$ will denote the complex set $\torus+{\rm i}[-s,s]$. As usual, $B^{\n}_r(x)$ denotes the ball in $\complex^{\n}$ with radius $r$ centered at $x$, accordingly to a prefixed norm $|\cdot|$ of $\complex^\n$.

\item[\tiny\textbullet]If $f=f(u, p,q,\varphi)$ is real-analytic for $(u, p,q,\varphi)\in W_{v, s, \varepsilon}= U_v \times B^{2\ell}_{\varepsilon}\times\torus^{\n}_s$, and affords the Taylor-Fourier expansion $$\dst f=\sum_{k\in \integer^{m}}f_{k, \a,\b}(u)e^{ik\cdot\varphi}\prod_{j=1}^\ell(\frac{p_j-{\rm i}q_j}{\sqrt2})^{\a_j}(\frac{p_j+{\rm i}q_j}{{\rm i}\sqrt2})^{\b_j},$$ we denote as $\|f\|_{v, s, \varepsilon}$ its ``sup-(Taylor, Fourier) norm'': $$\|f\|_{v, s, \varepsilon}:=\sum_{(a,b)\in \natural^{2\ell}\atop k\in \integer^{\n}}\sup_{u\in U_v}|f_{\a,\b, k}(u)|e^{|k|s}\varepsilon^{|(\a,\b)|}$$ with $|k|:=|k|_1$, $|(\a,\b)|:=|\a|_1+|\b|_1$.

\item[\tiny\textbullet] If $f$ is as in the previous item, $K>0$ and ${\eufm L}={\eufm L}_1\times{\eufm L}_2$ is a sub-lattice of $\integer^\n\times \integer^\ell$, $T_Kf$ and $\P_{\eufm L} f$ denote, respectively, the $K$-truncation and the ${\eufm L}$-projection of $f$: \beqano &&T_Kf:=\sum_{|(\a,\b)|\leq K,\ |k|\le K}f_{\a,\b,k}(u)e^{ik\cdot\varphi}\prod_{j=1}^\ell(\frac{p_j-{\rm i}q_j}{\sqrt2})^{\a_j}(\frac{p_j+{\rm i}q_j}{{\rm i}\sqrt2})^{\b_j}\nonumber\\
    && \P_{\eufm L}f:=\sum_{k\in {\eufm L}_1\atop \a-\b\in{\eufm L}_2}f_{\a,\b, k}(u)e^{ik\cdot\varphi}\prod_{j=1}^\ell(\frac{p_j-{\rm i}q_j}{\sqrt2})^{\a_j}(\frac{p_j+{\rm i}q_j}{{\rm i}\sqrt2})^{\b_j}. \eeqano
\end{enumerate}

\begin{proposition}
    [Multi-scale normal form]\label{iterative lemma} Let $$\n,\qquad \ell,\qquad 1\le m_1<\cdots<m_N=m$$ be natural numbers; $$A\subset\real^{\n},\ B\subset\real^{2\ell},\ C_1,\ C'_1\subset \real^{m_1},\ C_2,\ C'_2\subset \real^{m_2-m_1},\ \cdots,\ C_N,\ C'_N\subset \real^{m_N-m_{N-1}},$$ be open and connected sets; $$r,\ s,\ \varepsilon,\ \r_1\ge \r_2\cdots\ge \r_N,\ \r'_1\ge \r'_2\cdots\ge \r'_N$$ positive numbers. Put \beqano
    \begin{array}
        {llll} \dst v_i:=(r, \r_1,\cdots,\r_i,\r'_1,\cdots,\r'_i),\qquad &v:=v_N\\
        \\ \dst U^\ppi_{v_i}:=A_r\times {C_1}_{\r_1}\times\cdots\times {C_i}_{\r_i}\times {C'_1}_{\r'_1}\times\cdots\times {C'_i}_{\r'_i},\qquad &U_{v}:=U^{(N)}_{v_N}\\
        \\ \dst W^\ppi_{v_i, s, \varepsilon}:=U^\ppi_{v_i}\times\torus^\n_s\times B_{\varepsilon}, &W_{v, s, \varepsilon}:=W^{(N)}_{v_N, s, \varepsilon},
    \end{array}
    \eeqano with $i=1$, $\cdots$, $N$.
    
\nl Let ${\eufm a}$, $K>0$ with $0<s<6\log5/6$ and $Ks\ge 12$; let also ${\eufm L}$ and ${\eufm Z}_1$, $\cdots$, ${\eufm Z}_N$ be sub-lattices of $\integer^\ell\times \integer^\n$ and let ${\eufm Z}:={\eufm Z}_1\cup\cdots\cup {\eufm Z}_N$. 

\nl Let \beq{abstract system}H(u,\f,p,q)={\rm h}(p,q,I)+f(u,\f,p,q)\eeq be real-analytic for $(u, \f, p,q)\in W_{v, s, \varepsilon} $, where $u:$ $=$ $(I$, $\eta$, $\xi)$ $=$ $(I_1$, $\cdots$, $I_\n$, $\eta_1$, $\cdots$, $\eta_{m}$, $\xi_1$, $\cdots$, $\xi_{m})$. Suppose that
    \begin{enumerate}
        \item[{\rm (i)}] $h$ depends on $(p,q)$ only via $\frac{p_i^2+q_i^2}{2}$, with the frequency map $\o=(\o_1$, $\cdots$, $\o_\ell$, $\o_{\ell+1}$, $\cdots$, $\o_{\ell+\n})$ defined via $$\o_i:=\arr{ \dst \partial_{\frac{p_i^2+q_i^2}{2}}{\rm h}\quad 1\le i\le \ell\\
        \\
        \dst \partial_{I_{i-\ell}}{\rm h}\quad \ell+1\le i\le \ell+\n }$$ verifying \beq{non res} |\o(p,q, I)\cdot ( k',k)|\geq{\eufm a}\quad \forall\ (k',k)\in{\eufm Z}\setminus{\eufm L},\ |(k', k)|\leq K \eeq and all $(p,q,I)\in B^{2\ell}_\varepsilon\times A_r$; \item[{\rm (ii)}] $f$ is a sum \beq{f***}f=\sum_{i=1}^N f_i(u_i, \f, p, q) \eeq where $f_i$ is real-analytic on $W^\ppi_{v_i, s, \varepsilon}$ and has the form \begin{equation}\label{form of perturbation} \begin{split}f_i(u_i,\f, p, q)=\sum_{(\a^-\a^+, k)\in {\eufm Z}_i}f^i_{k,\a^-,\a^+}(u_i)\prod_{j=1}^{\n}e^{{\rm i}k_j\varphi_j} \prod_{k=1}^{\ell}\big(\frac{p_k-{\rm i}q_k}{\sqrt2}\big)^{\a_k^-}\big(\frac{p_k+{\rm i}q_k}{\sqrt2{\rm i}}\big)^{\a_k^+}\end{split}\end{equation} with \beqa{ui} &&u_i:=(I,\ \eta^{i}, \ \xi^{i}):=(I_1,\ \cdots,\ I_{\n}, \ \eta_1,\ \cdots, \ \eta_{m_i},\ \xi_1,\ \cdots,\ \xi_{m_i}); \eeqa \item[{\rm (iii)}] the following ``smallness'' conditions hold. If \beq{c and d}c_i:=e(1+\ell_i e+m_i e)/2,\qquad d_i:= \min\{ r s,\ \varepsilon^2, \ \r_i\r_i'\} \eeq with $e$ denoting Neper number, then \beqa{new smallness cond} \|f_{i}\|_{W^\ppi_{v_i,s,\varepsilon}}\le E_i,\qquad \sum_{i=1}^N\frac{7}{6}\big(\frac{9}{8}\big)^{i-1}\frac{2^7c_{i}K s}{{\eufm a} d_i}E_i<1 .\eeqa
    \end{enumerate}
    Then, one can find a real-analytic and symplectic transformation \[ \Phi:\quad W_{v/6^N,s/6^N, \varepsilon/6^N}\to W_{v,\s, \varepsilon} \] which conjugates $H$ to
    \begin{eqnarray*}
        H_{*}(u, \f,p,q):=H\circ \Phi={\rm h}(I,p,q)+\sum_{i=1}^Ng_i(u_i,\f,p,q)+\sum_{i=1}^Nf_i^{*}(u,\f,p,q),
    \end{eqnarray*}
    where $g_i$, $f_i$ verify \beqano g_i&=&\P_{{\eufm Z}_i\cap{\eufm L} }T_Kg_i\nonumber\\
    \|g_i-\P_{{\eufm Z}_i\cap{\eufm L}}T_Kf_i\|_{v_i/6^N,\s/6^N, \varepsilon/6^N}&\leq& (\frac{9}{8})^{2(i-1)}\frac{2^7c_{i}\,\|f_{i}\|^2_{v_{i},s,\varepsilon}}{{\eufm a} d_{i}}\nonumber\\
    &+&\frac{7}{6}(\frac{9}{8})^{2(i-1)}\sum_{j=1}^{i-1}\frac{2^7c_j\,\|f_j\|_{v_j,s,\varepsilon}}{{\eufm a} d_j}\|f_{i}\|_{v_{i},s, \varepsilon}\nonumber\\
    &+&\sum_{k=1}^{i-1}(\frac{9}{8})^{i-1-k}\frac{2^4c_k\|f_k\|_{v_k, s,\varepsilon}Ks}{{\eufm a} d_k}\|f_{i}\|_{v_{i}, s,\varepsilon}\nonumber\\
    \|f_i^*\|_{v_i/6^N, s/6^N,\varepsilon/6^N}&\le&\big(\frac{9}{8}\big)^{N-1}e^{-Ks/6^{i}}\| f_{i}\|_{v_{i}, s, \varepsilon} \eeqano Finally, $\Phi$ is close to the identity in the following sense. Given $F$, real-analytic on $W^\ppi_{v_i/6^N,s/6^N,\varepsilon/6^N}$, \beqano \|F\circ\Phi-F\|_{v/6^N,s/6^N,\varepsilon/6^N}\le \sum_{k=1}^N(\frac{9}{8})^{N-k}\frac{2^4c_k\|f_k\|_{v_k, s,\varepsilon}Ks}{{\eufm a} d_{k,i}}\|F\|_{v_i/6^N,s/6^N,\varepsilon/6^N} \eeqano with $d_{k,i}:=\max\{d_k,d_i\}$.
\end{proposition}
\nl The proof of Proposition \ref{iterative lemma} is based on the following
\begin{lemma}
    \label{iterative lemma*} Let $\bar N\in \natural$, $\n$, $\ell$, $m_i$, $A$, $B$, $C_i$, $C'_{i}$, $r$, $s$, $\r_i$, $\r'_i$, $U^\ppi_{v_i}$, $W^\ppi_{v_i, s,\varepsilon}$, $c_i$, $d_i$, with $i=1$, $\cdots$, ${\bar N}+1$, be as in Proposition \ref{iterative lemma}; $v:=(r, \r_1,\cdots,\r_{\bar N+1}, \r_1',\cdots,\r'_{\bar N+1})$, $U_v:=U^{({\bar N}+1)}_{v_{{\bar N}+1}} $, $W_{v, s,\varepsilon}:=W^{({\bar N}+1)}_{v_{{\bar N}+1},s,\varepsilon}$. 
Let \beq{gg} H(p,q,I,\varphi, \eta, \xi)={\rm h}(p,q,I)+g(p,q,I,\varphi, \eta, \xi) +f(p,q,I,\varphi, \eta, \xi) \eeq be real-analytic for $(u, \f, p,q)\in W_{v, s, \varepsilon} $. Suppose assumption {\rm (i)} of Proposition \ref{iterative lemma} and, moreover, the following ones
    \begin{enumerate}
        \item[{\rm (ii)}] $g$ is a sum \beq{gg***}g=\sum_{i=1}^{\bar N} g_i(u_i,\f,p,q)\eeq where $g_i$ is real-analytic on $W^\ppi_{v_i, s, \varepsilon}$ and $u_i$ is as in \equ{ui}; \item[{\rm (iii)}] $g_1$, $\cdots$, $g_{\bar N}$ and $f$ satisfy $$g_i=\P_{\eufm L}g_i,\quad f=\P_{\eufm Z}f$$ and \beqa{exponential smallness} \sum_{i=1}^{\bar N}\frac{2^7c_{i}K s}{{\eufm a} d_i}\|g_i\|_{v_i,s_i, \varepsilon_i}<1,\quad \|f\|_{v,s, \varepsilon}<\frac{{\eufm a} d_{\bar N+1}}{2^7c_{\bar N+1}K s}.\eeqa
    \end{enumerate}
    Then, one can find a real-analytic and symplectic transformation \[ \Phi:\quad (u', \varphi', p', q')\in W_{v/6,s/6, \varepsilon/6}\to (u, \varphi,p,q)\in W_{v,\s, \varepsilon} \] such that
    \begin{eqnarray*}
        H_{*}:=H\circ \Phi=h+g+g_*+f_{*},
    \end{eqnarray*}
    where $g_*=\P_{{\eufm Z}\cap{\eufm L} }T_Kg_*$ is ${\eufm Z}\cap{\eufm L}$-resonant and the following bounds hold \beqano \|g_*-T_{K}\P_{{\eufm Z}\cap{\eufm L}}f\|_{v/6,\s/6, \varepsilon/6}&\leq&\big(\frac{2^7c_{\bar N+1}\,\|f\|_{v,s,\varepsilon}}{{\eufm a} d_{\bar N+1}}+\sum_{i=1}^n\frac{2^7c_i\,\|g_i\|_{v_i,s,\varepsilon}}{{\eufm a} d_i}\big)\|f\|_{v,s, \varepsilon}\nonumber\\
    &\leq& \frac{\|f\|_{v,s, \varepsilon}}{6}\nonumber\\
    \|f_{*}\|_{v/6,\s/6, \varepsilon/6}&\leq& e^{-K s/6}\|f\|_{v,s, \varepsilon}. \eeqano Finally, $\Phi$ is close to the identity in the following sense: for any $F$ which is real-analytic on $W^\ppi_{v, s,\varepsilon}$, \beq{close to id}\|F\circ\Phi-\Phi\|_{v/6, s/6,\varepsilon/6}\le \frac{2^4c_{\bar N+1}\|f\|_{v, s,\varepsilon}Ks}{{\eufm a} d_i}\|F\|_{v_i, s,\varepsilon}<\frac{1}{8}\|F\|_{v, s,\varepsilon}.\eeq
\end{lemma}
\nl The following Lemma is a trivial extension\footnote{In order to obtain the extension it is sufficient to replace $\phi$ of \cite[Appendix A]{poschel93} with $$\phi=\sum_{{ (\a-\b, k)\in{\eufm K}\setminus {\eufm L}\atop |(\a,\b)|\le K,\ |k|\le K} }\frac{f_{k, \a,\b}(u)}{{\rm i}(\a-\b,k)\cdot\o}e^{ik\cdot\varphi}\prod_{j=1}^\ell(\frac{p_j-{\rm i}q_j}{\sqrt2})^{\a_j}(\frac{p_j+{\rm i}q_j}{{\rm i}\sqrt2})^{\b_j}$$} of \cite[Iterative Lemma]{poschel93}. Its proof is omitted.
\begin{lemma}
    \label{itera1} Let $s=(s_1,\cdots, s_\n)$, $r=(r_1,\cdots, r_\n)$, $\varepsilon=(\varepsilon_1,\cdots, \varepsilon_{\ell})$, $\r=(\r_1,\cdots, \r_{m})$, $\r'=(\r'_1,\cdots, \r'_{m})$, $v:=(r, \r, \r')$, $\hat v:=(\hat r, \hat\r, \hat\r')<v/2$, $\hat s<s/2$, $\hat\varepsilon<\varepsilon/2$, $$\d:=\min_{{i=1,\cdots, \n\atop j=1,\cdots, \ell}\atop k=1,\cdots, m}\{ \hat r_i \hat s_i,\ \hat\varepsilon_j^2, \hat\r_k\hat\r_k'\} .$$ Let \[ H(u,\varphi,p,q)={\rm h}(I, p,q)+g(u, \f, p,q) +f(u, \f,p,q)\qquad g(u, \f, p,q)=\sum_{i=1}^mg_i(u, \f, p,q)\] be real-analytic on $W_{v,s,\varepsilon}$. Assume that inequality \equ{non res} and \beqa{itera smallness} \|f\|_{v,s, \varepsilon}<\frac{{\eufm a} \d}{c} \eeqa are satisfied. Then one can find a real-analytic and symplectic transformation $$\Phi:\ W_{v-2\hat v, s-2\hat s,\varepsilon-2\hat\varepsilon}\to W_{v, s,\varepsilon}$$ defined by the time-one flow\footnote{The time-one flow generated by $\phi$ is defined as the differential operator $$X_\phi^1:=\sum_{k=0}^\infty \frac{{\rm L}_\phi^k}{k!}$$ where ${\rm L}^0_\phi f:=f$ and ${\rm L}^k_\phi f:=\big\{\phi, {\rm L}^{k-1}_\phi f\big\}$, with $k=1,\ 2,\ \cdots$. } $X_\phi^1f:=f\circ\Phi$ of a suitable $\phi$ verifying $$\|\phi\|_{v,s,\varepsilon}\le \frac{\|f\|_{v,s,\varepsilon}}{{\eufm a}}$$ such that $$H_+:=H\circ\Phi=h+g+\P_{{\eufm L}\cap{\eufm Z}}f+f_+$$ and, moreover, the following bounds hold \beqano \|f_+\|_{v-2\hat v, s-2\hat s, \varepsilon-2\hat\varepsilon}&\le& \big(1-\frac{c}{{\eufm a} \d}\|f\|_{v, s,\varepsilon}\big)^{-1}\Big[\frac{c}{{\eufm a} \d}\|f\|_{v, s,\varepsilon}^2 \nonumber\\
    &&+e^{-K\hat s}\|f\|_{v, s,\varepsilon}+\big(\frac{\varepsilon-\hat\varepsilon}{\varepsilon}\big)^{K}\|f\|_{v, s,\varepsilon}+ \|\big\{\phi, g \big\}\|_{v-\hat v, s-\hat s, \varepsilon-\hat\varepsilon}\Big] \eeqano Finally, for any real-analytic function $F$ on $W _{v, s, \varepsilon}$, \beqano \|F\circ\Phi-F\|_{v-2\hat v, s-2\hat s, \varepsilon-2\hat\varepsilon}&\le&\frac{\|\{\phi, F\}\|_{v-\hat v, s-\hat s, \varepsilon-\hat\varepsilon}}{\dst1-\frac{c\|f\|_{v, s,\varepsilon}}{{\eufm a} \d}}. \eeqano
\end{lemma}
\vskip.in\textsc{Proof of  Lemma \ref{iterative lemma*}.} Following \cite{poschel93}, the proof is obtained via iterate applications of Lemma \ref{itera1}. 

\nl To avoid too many indices, we shall prove this lemma taking, in \equ{gg***}, $\bar N=1$; the extension to $\bar N\ge 1$ being straightforward. Namely, we take \beqa{rho} &&\r_1=\cdots=\r_{m_1}=\bar\r,\quad \r'_1=\cdots=\r'_{m_1}=\bar\r'\nonumber\\
&&\r_{m_1+1}=\cdots=\r_{m}=\r,\quad \r'_{m_1+1}=\cdots=\r'_{m}=\r \eeqa where $1\le m_1<m$. Letting \beqano &&v:=(r,\ \r,\ \r'),\qquad \bar v:=(r,\ \bar\r,\ \bar\r'),\quad E:=\|f\|_{v,s,\varepsilon},\quad G:=\|g\|_{\bar v,s,\varepsilon},\quad \bar c=c_1,\quad c=c_2,\nonumber\\
&&\bar d:=\min\{ { r} { s},\ { \varepsilon}^2,{\bar\r}{\bar\r'}\},\qquad d:=\{ r s,\ \varepsilon^2, \r\r'\}, \eeqano we rewrite the assumptions in \equ{exponential smallness} as \beq{new cond}\frac{2^7 \bar c GKs}{{\eufm a}\bar d }<1,\quad \frac{2^7 c EKs}{{\eufm a} d }<1.\eeq The inequality on the right clearly implies \equ{itera smallness}. So, we apply Lemma \ref{itera1} to the Hamiltonian \equ{gg}, taking $r_1=\cdots=r_{\n }= r$, $s_1=\cdots=s_{\n }= s$, $\varepsilon_1=\cdots=\varepsilon_{\ell }= \varepsilon$, $\r_k$, $\r'_k$ as in \equ{rho} and \beqano &&\hat v=\hat v_0:=v/6,\quad \hat s=\hat s_0:=s/6,\quad \hat\varepsilon=\hat\varepsilon_0:=\varepsilon/6\nonumber\\
&&\hat{\bar v}=\hat {\bar v}_0:=\bar v/6,\quad \hat{\bar s}:=\hat{\bar s}_0:=\bar s/6,\quad \hat{\bar\varepsilon}:=\hat{\bar\varepsilon}_0:=\bar\varepsilon/6\nonumber\\
&&\d:=\{ \hat r \hat s,\ \hat\varepsilon^2, \hat\r\hat\r\}=\frac{d}{36}. \eeqano Letting \beqano v_1:=v-2 \hat v_0=3/4 v,\qquad s_1:=s-2\hat s=2/3 s,\quad \varepsilon_1:=\varepsilon-2\hat\varepsilon=2/3\varepsilon \eeqano by Lemma \ref{itera1}, we find a canonical transformation $\Phi_0=X_{\phi_0}$ which is real-analytic on $W_{ v_1, s_1, \varepsilon_1}$ and conjugates $H$ to $H_1=h+g+g_1+f_1$, where $g_1=\P_{{\eufm L}\cap{\eufm Z}}T_K f$ and

\beqano \|f_1\|_{v_1,s_1,\varepsilon_1}&\le&(1-\frac{36c E}{{\eufm a}d})^{-1}\Big[\frac{36c E}{{\eufm a}d}+e^{-Ks/6}+\big(\frac{5}{6}\big)^K\Big]E\nonumber\\
&+&(1-\frac{36 c E}{{\eufm a} d})^{-1}\frac{36\bar c G}{{\eufm a}\bar d} E\nonumber\\
\eeqano where $$\bar\d:=\min\{ \hat{ r} \hat{ s},\ \hat{ \varepsilon}^2,\hat{\bar \r}\hat{\bar \r}'\}=\frac{\bar d}{36}.$$ Here, we have used \beqa{Poschel trick} \|\big\{\phi,\ g\big\}_{I,\varphi, \eta, \xi}\|_{v-\hat v. s-\hat s, \varepsilon-\hat\varepsilon}&=& \|\big\{\phi,\ g\big\}_{I,\varphi, \eta^1, \xi^1}\|_{\bar v-\hat{\bar v}, s-\hat s, \varepsilon-\hat\varepsilon}\nonumber\\
&\le&\frac{\bar c G}{{\eufm a}\bar\d}=36\frac{\bar c G}{{\eufm a}\bar d} \eeqa since $g$ depends on $\eta$, $\xi$ only via $\eta^1=(\eta_1,\cdots, \eta_{m_1})$, $\xi^1=(\xi_1,\cdots, \xi_{m_1})$. It is sufficient to consider the case $$e^{-Ks/6}+\big(\frac{5}{6}\big)^K\le \frac{18c E}{{\eufm a}d}$$ since otherwise the Lemma is proved. In such case, using \equ{new cond} we can write \beqa{intermediate step***} E_1=\|f_1\|_{v_1,s_1,\varepsilon_1} &\le&\frac{32}{23}\big(\frac{9}{32}\frac{2^7 cEKs}{{\eufm a}d}+\frac{9}{64}\frac{2^7 cEKs}{{\eufm a}d}+\frac{9}{32}\frac{2^7 \bar cGKs}{{\eufm a}\bar d}\big)\frac{E}{Ks}\nonumber\\
&<&\frac{E}{Ks}\max\Big\{\frac{2^7 cEKs}{{\eufm a}d},\ \frac{2^7 \bar cGKs}{{\eufm a}\bar d}\Big\}<\frac{E}{4} \eeqa Let $$L:=\big[\frac{Ks}{12\log 2}\big].$$ Note that \beq{L and K}L\ge 1,\qquad Ks>8 L,\eeq since we have assumed $Ks\ge 12$. We want to prove that Lemma \ref{itera1} can be applied $L$ times with parameters \beq{delta i***} \hat v_i=\frac{v}{4L},\quad \hat \varepsilon_i=\frac{\varepsilon}{4L},\quad \hat s_i=\frac{s}{4L},\quad \d_i= \frac{d}{16L^2},\quad i=1,\cdots, L.\eeq For $L=1$, this follows from \equ{intermediate step***}: \beqano E_1:=\|f_1\|_{v_1,s_1,\varepsilon_1}&\le&\frac{E}{Ks}\le2^{-7}\frac{{\eufm a} d }{c (Ks)^2}<2^{-13}\frac{{\eufm a} \d_1 }{c} \eeqano which is implied by the inequality in \equ{intermediate step***} and assumption \equ{exponential smallness}. We then assume $L\ge 2$. Suppose, by induction, that, for a certain $1\le i\le L-1$, and any $1\le j\le i$, we have conjugated $H$ to $$H_j={\rm h}+g+\bar g_j+f_j$$ where $\bar g_j=\sum_{k=0}^{j-1}\P_{{\eufm L}\cap{\eufm Z}}T_K f_k$ \beq{Ej}E_j:=\|f_j\|_{v_j, s_j, \varepsilon_j}\le \min\big\{\frac{E}{4^j},\ 2^{-6}\frac{{\eufm a}\d_j}{c}\big\} \eeq where $\hat v_0$, $\hat s_0$, $\hat \varepsilon_0$ are as above, $v_0:=v$, $s_0:=s$, $\varepsilon_0:=\varepsilon$ and $v_j=v_{j-1}-2\hat v_{j-1}$. Then by Lemma \ref{itera1}, on the domain $W_{v_{j+1}, s_{j+1}, \varepsilon_{j+1}}$, we fined a real-analytic transformation $\Phi_i=X_{\phi_i}$, which conjugates $H_i$ to $$H_{i+1}={\rm h}+g+\bar g_{i+1}+f_{i+1}$$ where $\bar g_{i+1}=\bar g_i+\P_{{\eufm L}\cap {\eufm K}}f_i=\sum_{k=0}^i\P_{{\eufm L}\cap{\eufm Z}}T_K f_k$. We prove that \equ{Ej} is satisfied for $j=i+1$. Using\footnote{For the proof of inequality $\|\big\{g_i,\phi_i\big\}\|_{v_i-\hat v_i,s_i-\hat s_i,\varepsilon_i-\hat\varepsilon_i}\le \frac{cE_i}{{\eufm a}\d_1}\big(E_1+\frac{E}{L}\big)$, compare \cite[Proof of the Normal Form Lemma]{poschel93}.} the assumption on the right in \equ{new cond}, \equ{intermediate step***}, the inequality for $Ks$ in \equ{L and K} and the definition of $\d_i$ in \equ{delta i***}, we have

\beqano \|\big\{\bar g_i,\phi_i\big\}\|_{v_i-\hat v_i,s_i-\hat s_i,\varepsilon_i-\hat\varepsilon_i}&\le& \big[\frac{c}{{\eufm a}\d_i}\big(E_1+\frac{E}{L}\big)\big]E_i\le\Big[\frac{c}{{\eufm a}\d_i}\frac{E}{Ks}+\frac{c}{{\eufm a}\d_i}\frac{E}{L}\Big]E_i<\frac{E_i}{32}. \eeqano Moreover, by a similar argument as in \equ{Poschel trick} and since $g$ is actually real-analytic in the larger domain $$W_{\bar v, s, \varepsilon}\supset W_{\bar v_i-\hat{\bar v}_i+\bar v, s_i-\hat s_i+\hat s, \varepsilon_i-\hat\varepsilon_i+\hat\varepsilon},$$ we have \beqano \|\big\{g,\phi_i\big\}\|_{ v_i-\hat{ v}_i, s_i-{\hat s}_i,\varepsilon_i-{\hat\varepsilon}_i}&=&\|\big\{g,\phi_i\big\}\|_{\bar v_i-\hat{\bar v}_i, s_i-{\hat s}_i,\varepsilon_i-{\hat\varepsilon}_i}\le\frac{\bar cE_i}{{\eufm a}\bar \d_i}\frac{G}{L}<\frac{E_i}{64}, \eeqano where \[ \bar\d_i:=\min\{\hat{r}_i\hat{s}_i, \bar{\hat\r}_i\hat{\bar \r}_i'\}= \frac{\bar d}{16L^2},\quad i=1,\cdots, L.\] Then we find\footnote{Since $K>8L$ and $L\ge 2$, one has $(1-\frac{3}{2L})^K\le \frac{1}{(1+\frac{3}{2L})^{8L}}$ with the r.h.s bounded above by $(4/7)^{16}$ (it decreases to $e^{-12}$ as $L\to+\infty$).} \beqano E_{i+1}=\|f_{i+1}\|_{v_{i+1},s_{i+1},\varepsilon_{i+1}}&\le&(1-\frac{c E_i}{{\eufm a}\d_1})^{-1}\Big[\frac{c E_i}{{\eufm a}\d_1}+e^{-K\hat s_i}+\big(\frac{\varepsilon_{i}-\hat\varepsilon_i}{\varepsilon_{i}}\big)^K\Big]E_i\nonumber\\
&+& (1-\frac{c E_i}{{\eufm a}\d_1})^{-1}\|\big\{\bar g_i,\phi_i\big\}\|_{v_i-\hat v_i,s_i-\hat s_i,\varepsilon_i-\hat\varepsilon_i}\nonumber\\
&+& (1-\frac{c E_i}{{\eufm a} \d_1})^{-1} \|\big\{g,\phi_i\big\}\|_{\bar v_i-\bar{\hat v}_i, s_i-{\hat s}_i,\varepsilon_i-{\hat\varepsilon}_i} \nonumber\\
&\le&\frac{64}{63}\big[\frac{1}{64}+\frac{1}{8}+(\frac{4}{7})^{16}+\frac{1}{32}+\frac{1}{64}\big]E_{i}\nonumber\\
&<&\frac{E_{i}}{4}< E_{1}<2^{-6}\frac{{\eufm a}\d_1}{c}.\eeqano since $i\ge 1$. Then we let $\Phi:=\Phi_0\circ\cdots\circ\Phi_L$, $H_*:=H\circ\Phi=h+g+\bar g_{L+1}+f_{L+1}$, $g_*:=g_{L+1}$, $f_*:=f_{L+1}$ and we have, by telescopic inequalities and \equ{intermediate step***}, \beqano \|g_{*}-\P_{{\eufm L}\cap{\eufm K}}T_K f\|_{v/6, s/6, \varepsilon/6}&=&\sum_{i=1}^L\|\P_{{\eufm L}\cap{\eufm K}}T_K f_i\|\le \sum_{i=1}^L E_i\le E_1\sum_{i=1}^L\frac{1}{4^{i-1}}\nonumber\\
&=&\frac{4}{3}E_1\le (\frac{2^7 c E}{{\eufm a} d}+\frac{2^7 \bar cG}{{\eufm a}\bar d})E\eeqano Now we prove \equ{close to id}. Let $F\in W_{\bar v,s,\varepsilon}$, $F_{-1}:=F$, $F_i:=F\circ\Phi_0\circ\cdots\circ\Phi_i$, $i=0$, $\cdots$, $L$. Then \beqano \|F\circ\Phi-F\|_{\bar v/6, s/6, \varepsilon/6}&=&\|F_L-F\|_{\bar v_{L+1}, s_{L+1}, \varepsilon_{L+1}}\le \sum_{i=0}^L\|F_{i-1}\circ\Phi_i-F_{i-1}\|_{\bar v_{i+1}, s_{i+1}, \varepsilon_{i+1}}\nonumber\\
&\le&\sum_{i=0}^L\frac{\frac{\bar cE_i}{{\eufm a} \bar\d_i}}{(1-\frac{\bar cE_i}{{\eufm a} \bar\d_i})}\|F\|_{\bar v_i, s_i, \varepsilon_i}\le \frac{\sum_{i=0}^L\frac{\bar cE_i}{{\eufm a} \bar\d_i}}{\prod_{i=0}^L(1-\frac{\bar cE_i}{{\eufm a} \bar\d_i})}\|F\|_{\bar v, s, \varepsilon}\nonumber\\
&\le&\sum_{i=0}^L\frac{\bar cE_i}{{\eufm a} \bar\d_i}e^{\frac{5}{4}\sum_{i=0}^L\frac{\bar cE_i}{{\eufm a} \bar\d_i}}\|F\|_{\bar v, s, \varepsilon}\le \frac{2^5\bar cE_0Ks}{{\eufm a} d}\|F\|_{\bar v, s, \varepsilon} \eeqano where we have used $\frac{\bar cE_i}{{\eufm a} \bar\d_i}<1/24$ that, for $0\le x\le 1/24$, $\log(1-x)^{-1}<\frac{5}{4}x$ and \beqano \sum_{i=0}^L\frac{\bar cE_i}{{\eufm a} \bar\d_i}&=&\frac{\bar cE_0}{{\eufm a} \bar\d_0}+\sum_{i=1}^L\frac{\bar cE_i}{{\eufm a} \bar\d_i}\le\frac{2^6\bar cE_0}{{\eufm a} d}+\frac{\bar cE_1}{{\eufm a} \bar\d_1}\sum_{i=1}^L\frac{1}{4^{i-1}}\nonumber\\
&\le&\frac{2^6\bar cE_0}{{\eufm a} d}+\frac{4}{3}\frac{\bar cE_1}{{\eufm a} \bar\d_1}<\frac{2^4\bar cE_0Ks}{{\eufm a} d}. \eeqano The proof for $F\in W_{v,s,\varepsilon}$ is similar. \qed

\vskip.in\textsc{Proof of Proposition \ref{iterative lemma}.} For simplicity of notations, we prove Proposition \ref{iterative lemma} in the case $\n=\ell=1$; the generalization to any $\n$, $\ell$ being straightforward. Consider the Hamiltonian $$H_0(u_1,\f,p,q):={\rm h}(I,p,q)+f_1(u_1,\f,p,q),\qquad (u_1,\f,p,q)\in W^\ppu_{v_1, s, \varepsilon}.$$ To this Hamiltonian let us apply Lemma \ref{iterative lemma*}, with $g\equiv 0$, so as to conjugate it to $$H_1:=H_0\circ\Phi_1=h+g_1+f_{*1}^\ppu,\qquad (u_1,\f,p,q)\in W^\ppu_{v_1/6, s/6, \varepsilon/6}$$ where $g_1$, $f_{*1}^\ppu$ correspond to $g_*$, $f_*$, hence satisfy \beqano \|f_{*1}^{(1)}\|_{v_{1}/6, s/6, \varepsilon/6}&\le& e^{-Ks/6}\|f^\ppi_{1}\|_{v_{1}, s, \varepsilon}\nonumber\\
\nonumber\\
\|g_{1}\|_{v_{1}/6, s/6, \varepsilon/6}&\le&\frac{7}{6}\|f_{1}\|_{v_{1}, s, \varepsilon}\nonumber\\
\|g_{1}-\P_{{\eufm L}\cap{\eufm Z}}T_K f_1\|_{v_{1}/6, s/6, \varepsilon/6}&\le&\frac{2^7c_{1}\,\|f_{1}\|^2_{v_{1},s,\varepsilon}}{{\eufm a} d_{1}} \eeqano Then we have $$H^\ppu(u, \f, p,q):=H\circ\Phi_1=H_0\circ\Phi_1+\sum_{j=2}^N f_j\circ\Phi_1={\rm h}+g_1+f_{1*}^\ppu+\sum_{j=2}^N f_j^\ppu$$ where $f_j^\ppu:=f_j\circ\Phi_1$. Assume, inductively, that, for some $1\le i\le N-1$ and any $1\le j\le i$ we have conjugated $H$ to $$H^\ppj(u, \f, p,q)=H\circ\Phi_1\circ\cdots\circ\Phi_j={\rm h}+\sum_{k=1}^j g_k+\sum_{k=1}^j f_{k*}^\ppj+\sum_{k=j+1}^Nf_k^\ppj$$ where $$\Phi_j:\ W^\ppj_{v/6^j, s/6^j, \varepsilon/6^j}\to W^{(j-1)}_{v/6^{j-1}, s/6^{j-1}, \varepsilon/6^{j-1}}$$ transforms $$H_{j-1}:={\rm h}+\sum_{k=1}^{j-1} g_k+f^{(j-1)}_j$$ into $$H_{j-1}\circ\Phi_j={\rm h}+\sum_{k=1}^{j} g_k+f_{*j}^\ppj.$$ The Hamiltonian $$H_i(u_{i+1}, \f,p,q):={\rm h}+\sum_{k=1}^i g_k(u_k, \f,p,q)+f^\ppi_{i+1}(u_{i+1}, \f,p,q)$$ is real-analytic for $(u_{i+1}, \f,p,q)\in W^{(i+1)}_{v_{i+1}/6^{i}, s/6^{i}, \varepsilon/6^{i}}$ and satisfies the assumptions of Lemma \ref{iterative lemma*}, with $\bar N=i$. Then one can find $\Phi_{i+1}:\ W^{(i+1)}_{v_{i+1}/6^{i+1}, s/6^{i+1}, \varepsilon/6^{i+1}}\to W^{(i+1)}_{v_{i+1}/6^{i}, s/6^{i}, \varepsilon/6^{i}}$ such that $H_{i}\circ\Phi_{i+1}={\rm h}+\sum_{k=1}^{i+1} g_k+f_{*i+1}^{(i+1)}$, where \beqano \|f_{*i+1}^{(i+1)}\|_{v_{i+1}/6^{i+1}, s/6^{i+1}, \varepsilon/6^{i+1}}&\le& e^{-Ks/6^{i+1}}\|f^\ppi_{i+1}\|_{v_{i+1}/6^{i}, s/6^{i}, \varepsilon/6^{i}}\nonumber\\
&\le&\big(\frac{9}{8}\big)^{i}e^{-Ks/6^{i+1}}\| f_{i+1}\|_{v_{i+1}, s, \varepsilon}\nonumber\\
\|g_{i+1}\|_{v_{i+1}/6^{i+1}, s/6^{i+1}, \varepsilon/6^{i+1}}&\le&\frac{7}{6}\|f^\ppi_{i+1}\|_{v_{i+1}/6^{i}, s/6^{i}, \varepsilon/6^{i}}\le \frac{7}{6}\big(\frac{9}{8}\big)^{i}\|f_{i+1}\|_{v_{i+1}, s, \varepsilon}\nonumber\\
\|g_{i+1}-\P_{{\eufm L}\cap{\eufm Z}}T_K f_{i+1}\|_{v_{i+1}/6^{i+1}, s/6^{i+1}, \varepsilon/6^{i+1}}&\le&\|g_{i+1}-\P_{{\eufm L}\cap{\eufm Z}}T_K f_{i+1}^\ppi\|_{v_{i+1}/6^{i+1}, s/6^{i+1}, \varepsilon/6^{i+1}}\nonumber\\
&+&\|\P_{{\eufm L}\cap{\eufm Z}}T_Kf^\ppi_{i+1}-\P_{{\eufm L}\cap{\eufm Z}}T_K f_{i+1}\|_{v_{i+1}/6^{i+1}, s/6^{i+1}, \varepsilon/6^{i+1}}\nonumber\\
&\le&\|g_{i+1}-\P_{{\eufm L}\cap{\eufm Z}}T_K f_{i+1}^\ppi\|_{v_{i+1}/6^{i+1}, s/6^{i+1}, \varepsilon/6^{i+1}}\nonumber\\
&+&\|f^\ppi_{i+1}-f_{i+1}\|_{v_{i+1}/6^{i+1}, s/6^{i+1}, \varepsilon/6^{i+1}}\nonumber\\
&\le&(\frac{9}{8})^{2i}\frac{2^7c_{i+1}\,\|f_{i+1}\|^2_{v_{i+1},s,\varepsilon}}{{\eufm a} d_{i+1}}\nonumber\\
&+&\frac{7}{6}(\frac{9}{8})^{2i}\sum_{j=1}^{i}\frac{2^7c_j\,\|f_j\|_{v_j,s,\varepsilon}}{{\eufm a} d_j}\|f_{i+1}\|_{v_{i+1},s, \varepsilon}\nonumber\\
&+&\sum_{k=1}^i(\frac{9}{8})^{i-k}\frac{2^4c_k\|f_k\|_{v_k, s,\varepsilon}Ks}{{\eufm a} d_k}\|f_{i+1}\|_{v_{i+1}, s,\varepsilon}\eeqano with $f_{k*}^{(i+1)}:=f_{k*}^{(i)}\circ\Phi_{i+1}$ for $1\le k\le i+1$ and $f_k^{(i+1)}:=f_k^{(i)}\circ\Phi_{i+1}$ for $i+2\le k\le N$. Then we find \beqano H^{(i+1)}&:=&H^\ppi\circ\Phi_{i+1}=({\rm h}+\sum_{k=1}^i g_k+\sum_{k=1}^i f_{k*}^\ppi+\sum_{k=i+1}^Nf_k^\ppi)\circ\Phi_{i+1}\nonumber\\
&=&H_i\circ\Phi_{i+1}+(\sum_{k=1}^i f_{k*}^\ppi+\sum_{k=i+2}^Nf_k^\ppi)\circ\Phi_{i+1}\nonumber\\
&=&{\rm h}+\sum_{k=1}^{i+1} g_k+\sum_{k=1}^{i+1}f_{k*}^{(i+1)}+\sum_{k=i+2}^Nf_k^{(i+1)} \eeqano and hence, after $N$ steps, \beqano H^{(N)}&:=&H\circ\Phi_1\cdots\circ\Phi_N={\rm h}+\sum_{k=1}^{N} g_k+\sum_{k=1}^{i+1}f_{k*}^{(N)}\eeqano satisfies the thesis of Proposition \ref{iterative lemma}. \qed

\section{A slightly-perturbed integrable system}

The following result is well known in the literature of close-to be integrable systems, hence its proof is omitted. Note that it deals with an integrable system, close to another integrable one.
\begin{theorem}
    \label{perturbation argument} One can find a number ${\eufm c}_0$ such that, for any real-analytic, one-dimensional, system $${\rm H}(P,Q)={\rm h}(\frac{P^2+Q^2}{2})+f(P,Q)\qquad (P, Q)\in {\eufm B}=B^2_{\varepsilon}(0)\subset \complex^2$$ and any $0<\bar\varepsilon<\varepsilon$, such that \beq{smallness integration}\inf_{B^2_{\varepsilon}}|\partial {\rm h}|\ge {\eufm a},\qquad \sup_{B^2_{\varepsilon}}|f|\le \eufm{e},\qquad \frac{1}{\eufm{c}_0}\frac{\eufm{e}}{{\eufm a}\bar\varepsilon^2}<1,\eeq one can find a real-analytic transformation $$\phi_*:\quad (P_*, Q_*)\in B^2_{\varepsilon-\bar\varepsilon}\to (P,Q)\in B^2_{\varepsilon}$$ which conjugates ${\rm H}$ to a function ${\rm H}_*={\rm H}\circ \phi_*$ depending only on $\frac{P_*^2+Q_*^2}{2}$. The assertion can be extended to the case that ${\rm h}$, $f$ are functions of other canonical coordinates $(P',Q', {\rm y}, {\rm x})$, depending on them only via ${\rm Y}=({\rm y}, \frac{{P'_1}^2+{Q'_1}^2}{2}, \cdots, \frac{{P'_m}^2+{Q'_m}^2}{2})$, with ${\rm y}\in{\rm Y}_\r$, $(P'_j, Q'_j)\in B^2_{\varepsilon'_j}$. In this case, letting $(P_*,Q_*)\to\phi_*(P_*,Q_*; {\rm Y})$ the transformation obtained for any fixed value of ${\rm Y}$, there exists a canonical, real-analytic, transformation $\Phi_*$ of the form 
    $$\Phi_*:\ (P,Q)=\phi_*(P_*,Q_*; {\rm Y}_*)\ {\rm y}={\rm y}_*,\ {\rm x}={\rm x}_*+\f ({\rm Y}_*),\ P_j'+{\rm i}Q_j'=e^{{\rm i}\psi_j ({\rm Y}_*)}(P_{*j}'+{\rm i} Q_{*j}')$$ which conjugates ${\rm H}$ to a function ${\rm H}_*={\rm H}\circ \Phi_*$ depending only on $\frac{P_*^2+Q_*^2}{2}$ and ${\rm Y}_*$. In this case, the functions $\varphi_j$, $\psi_j$ verify $$|\varphi_j|\le \frac{1}{{\eufm c}_0}\frac{{\eufm e}}{{\eufm a}\r_j},\qquad |\psi_j|\le \frac{1}{{\eufm c}_0}\frac{{\eufm e}}{{\eufm a}{\varepsilon'_j}^2}.$$
\end{theorem}
\chapter[The geometrical structure of the $\cP$-coordinates]{More on the geometrical structure of the $\cP$-coordinates, compared to Deprit's coordinates}\label{comparison} In this section we aim to point out differences and similarities between the $\cP$-coordinates and the coordinates denoted as $(\Psi, \G,\L,\psi,\g,\ell)$ in \cite{chierchiaPi11a, pinzari-th09, chierchiaPi11b}. 

\nl We recall that the ``planetary'' coordinates $(\Psi, \G,\L,\psi,\g,\ell)$ may be derived (after a canonical transformation) from a more general set of canonical coordinates studied by A. Deprit. In their planetary form, the coordinates $(\Psi, \G,\L,\psi,\g,\ell)$ have been rediscovered\footnote{The proof of their symplectic character found in \cite{pinzari-th09} has been published in \cite{chierchiaPi11a}. Another proof has been given in \cite{zhao14}.} by the author during her PhD, under the strong motivation of their application to the planetary problem \cite{pinzari-th09, chierchiaPi11b}.

\nl Let us recall their definition\footnote{For sake of uniformity, we use slightly different notations with respect to the ones in \cite{chierchiaPi11a}, actually closer to the ones of the paper \cite{deprit83}).}, in the spirit of {\sl Kepler maps} (Definition \ref{Kepler map}).

\nl Let ${\rm C}_\cE^\ppi$, ${\rm S}_\cE^\ppi$ be as in \equ{Cj Sj} of Section \ref{perihelia reduction} and define the {$\cD ep$-nodes} \begin{equation}\label{D nodes}\begin{split} n_i:=\left\{
\begin{array}
    {llll} k^\ppt\times {\rm S}_\cE^\ppu& i=0\\
    \\
    \dst{\rm S}_\cE^\ppi\times {\rm S}_\cE^{(i+1)}=-{\rm S}_\cE^\ppi\times {\rm C}_\cE^{(i)}\quad &i=1,\cdots, n-1.\\
    \\
    \dst-n_{n-1}& i=n
\end{array}
\right. \end{split}\end{equation} Then let $$\cE_{\cD ep}:=\big\{ (({\eufm E}_1,\cdots, {\eufm E}_n)\subset E^{3}\times \cdots\times E^3):\quad 0<e_i<1,\quad n_{i-1}\ne 0\quad \forall\ i=1,\cdots, n\big\}.$$ On $\cE_{\cD ep}$, define the map $$ \t_{\cD ep}^{-1}:\quad ({\eufm E}_1,\cdots, {\eufm E}_n)\in \cE_{\cD ep}\to {\rm X}_{\cD ep}\in{\eufm X}_{\cD ep}= \t_{\cD ep}^{-1}(\cE_{\cD ep}) $$ where $${\rm X}_{\cD ep}=(\Psi,\G,\L,\psi,\g)\in \real^n\times \real_+^n\times\real_+^n\times \torus^n\times \torus^n$$ where \beqano
\begin{array}
    {lll} &\dst\Psi=(\Psi_{-1},\Psi_0, \bar\Psi)\in \real_+\times\real_+\times \real_+^{n-2}& \psi=(\psi_{-1},\psi_0, \bar\Psi)\in \torus\times\torus\times \torus^{n-2}\\
    \\ &\dst\G=(\G_1,\cdots, \G_n)\in \real_+^{n}&\g=(\g_1,\cdots, \g_n)\in \torus^{n}\\
    \\ &\dst\L=(\L_1,\cdots,\L_n)\in \real_+^n\quad &
\end{array}
\eeqano with \beqano \bar\Psi=(\Psi_1,\cdots,\Psi_{n-2})\quad \bar\psi=(\psi_1,\cdots,\psi_{n-2}) \eeqano are defined as follows. The coordinates $\L_j$ are as in \equ{belle*}, while $(\Psi, \G,\psi, \g)$ are defined as \begin{equation}\label{planetary Deprit}\begin{split}
\begin{array}
    {llllrrr} \dst \Psi_{i-2}=\left\{
    \begin{array}
        {lrrr} \dst Z:={\rm S}^\ppu_\cE\cdot k^\ppt \\
        \\
        \dst |{\rm S}_\cE^{(i)}|
    \end{array}
    \right.& \psi_{i-2}=\left\{
    \begin{array}
        {lrrr} \dst\zeta:=\a_{k^\ppt}(k^\ppu, n_0)\\
        \\
        \dst \a_{{\rm S}_\cE^{(i-1)}}(n_{i-2}, n_{i-1})
    \end{array}
    \right.\qquad &\left.
    \begin{array}
        {lrrr} \dst i=1\\
        \\
        2\le i\le n
    \end{array}
    \right. \\
    \\
    \dst \G_{i}:=|{\rm C}_\cE^\ppi| & \g_i:=\a_{{\rm C}_\cE^\ppi}(n_i, P^\ppi)&\left.
    \begin{array}
        {lrrr} 1\le i\le n
    \end{array}
    \right.
\end{array}
\end{split}\end{equation}
Then $\t^{-1}_{\cD ep}$ is a bijection \cite{deprit83, pinzari-th09, chierchiaPi11a, zhao14}.
\begin{definition}
    \rm We call {\sl Deprit's map}, or $\cD ep$ map, the Kepler map \[ \cD ep:\quad {\rm Dep}= ({\rm X}_{\cD ep},\ell) \in\cD_{\cD ep}={\eufm X}_{\cD ep}\times \torus^n\to (y, x)\in\real^{3n}\times \real^{3n}\] associated to $\t_{\cD ep}$.
\end{definition}
 \vskip.1in\noi
{\bf \bf Comparing $\cP$ and $\cD ep$}
\vskip.1in
\noi
{\bf a)} Both the $\cP$ and ${\cD ep}$-coordinates reduce the system to $(3n-2)$ degrees of freedom. They share the following three coordinates (two actions and an angle) $$\Psi_{-1}=Z=\Theta_0,\quad \psi_{-1}=\zeta=\vartheta_0,\quad \Psi_0={\rm G}=\chi_0$$ which are integrals of the system. As a consequence, the coordinates $(Z, \zeta)$ and, respectively, $${\rm g}:=\psi_0,\qquad {\eufm g}:=\k_0$$ do not appear into the Hamiltonian. Note that ${\cD ep}$ and $\cP$ share also the fixed node $n_0=\n_1$. \vskip.1in\noi{\bf b)} The angle ${\rm g}$ for the set ${\cD ep}$ describes the motion of the node $n_1$ in \equ{D nodes} and, by the cyclic character of ${\rm g}$, this motion is negligible. Its counterpart in the set $\cP$ is {the} node ${\rm n}_1$ in \equ{good nodes}, the negligible motion of which is governed by ${\eufm g}$. 

\vskip.1in\noi{\bf c)} Compare the diagrams in \equ{chain} and \equ{P chain} with the two ones associated to the ${\cD ep}$-map, respectively:

 \beqano
\begin{array}
    {cccccccccccccccccc} \dst&&n_0&&n_1&&\vdots&&n_{n-2}&n_{n-1}&&\\
    \\
    \dst&&\Uparrow&&\Uparrow&&\vdots&&\Uparrow&\Uparrow&\\
    \\
    \dst k^\ppt&\to&{\rm S}_\cE^\ppu&\to& {\rm S}_\cE^\ppd&\to&\cdots&\to& {\rm S}_\cE^{(n-1)}&\to{\rm S}_\cE^\ppn= {\rm C}_\cE^\ppn\\
    \\
    \dst&&\downarrow&&\downarrow&&\vdots&&\downarrow&&\\
    \\
    \dst&&{\rm C}_\cE^\ppu&&{\rm C}_\cE^\ppd&&\vdots&&{\rm C}_\cE^{(n-1)}&&&\\
    \\
    \dst&&\Downarrow&&\Downarrow&&\vdots&&\Downarrow&&\\
    \\
    \dst&&-n_1&&-n_2&&\vdots&&-n_{n-1}&&&
\end{array}
\eeqano and \beqano
\begin{array}
    {llllllllllllllllllll} {\rm F}_0&\to&{\rm F}_1^*&\to&\cdots&\to& {\rm F}_i^*&\to& \cdots&\to&{\rm F}_n^*={\rm G}^*_n\\
    \\
    \dst&&\downarrow&&\vdots&&\downarrow&&\vdots&&\downarrow\\
    \\
    \dst&&{\rm G}_1^*&&&&{\rm G}_i^*&&&&{\rm G}_n^*
\end{array}
\eeqano where \beqano &&{\rm F}_i^*=(n_{i-1}, \ \cdot , {\rm S}_\cE^\ppi)\quad{\rm G}_i^*=(-n_i,\ \cdot,{\rm C}_\cE^\ppi)\qquad i=1,\cdots, n. \eeqano Note that, analogously to \equ{chain}, $n_i$ in \equ{D nodes} is the skew-product of its two previous vectors in the tree \equ{chain}.

 \vskip.1in\noi{\bf d)} While ${\cD ep}$ is not defined for the planar problem, $\cP$ is, and, in that case, the coordinates $(\Theta, \chi, \vartheta, \k)$ in \equ{belle*} reduce to\footnote{ Here by ``planar case'' we mean ${\rm C}_\cE^\ppu\parallel\cdots\parallel {\rm C}_\cE^\ppn\parallel k^\ppt$. Note that, to be more precise, $\vartheta_0$ and $\k_0$ would not exist in that case (since $\n_1=0$). However, since they are both cyclic angles, we can fix them to an arbitrary value. The choice above corresponds to replace $\n_1$ with $k^\ppu$. } \beqano
\begin{array}
    {llll} \dst \Theta_i=\arr{ \dst\chi_0 \\
    \\
    \dst0 }\qquad &\vartheta_i=\arr{ \dst0 \\
    \\
    \dst\p }\qquad \k_{i}=\arr{ \dst\arg P^\ppu-\frac{\p}{2} \\
    \\
    \dst\widehat{ P^{(i)}P^{(i+1)}}+\p } \qquad&
    \begin{array}
        {llll} \dst i=0\\
        \\
        \dst i=1,\cdots, n-1
    \end{array}
    \nonumber\\
    \nonumber\\
    \dst \chi_{i}=\sum_{j=i+1}^n\|{\rm C}_\cE^\ppj\| &
\end{array}
\eeqano while the $(\L, \ell)$ remain unchanged. \vskip.1in\noi{\bf e)} The $\cP$-map is singular when some eccentricity $e_i$ vanishes or some of the following relations hold $${\rm S}_\cE^\ppu\parallel k^\ppt\qquad P^\ppi\parallel {\rm S}_\cE^\ppi\qquad {\rm S}_\cE^{(i+1)}\parallel P^\ppi.$$ The former of such relations is negligible, while the other ones have no physical meaning. Therefore, the only physically relevant singularities of $\cP$ are for zero-eccentric motions.\\
The ${\cD ep}$-map is singular when some eccentricity $e_i$ vanishes or some of the following relations hold $${\rm S}_\cE^\ppu\parallel k^\ppt\qquad {\rm S}_\cE^{(i+1)}\parallel {\rm S}_\cE^\ppi\qquad i=1,\cdots, n-1.$$ The configurations ${\rm S}_\cE^\ppi\parallel {\rm S}_\cE^{(i+1)}$ have a relevant physical meaning, since the planar case corresponds to the intersection of all such configurations. A complete regularization of {\sl all} the singularities of the ${\cD ep}$-map has been obtained in \cite{pinzari-th09, chierchiaPi11b}, which allowed to overcome the problem of the {\sl rotational degeneracy} (see \cite{chierchiaPi11c} for information) of the planetary problem and to construct the Brkhoff normal form of it. It works at expenses of one extra-degree of freedom. \vskip.1in\noi{\bf f)} The Euclidean lengths $\|{\rm C}_\cE^\ppi\|$ of the planets' angular momenta are the actions $\G_i$ among ${\cD ep}$-coordinates: see \equ{planetary Deprit}. In terms of the $\cP$-coordinates they have more involved expressions in \equ{CP}. As mentioned in the previous item, this makes more difficult regularizing singular configurations with zero eccentricity. The formula simplifies in the planar case: $$\|{\rm C}_\cE^\ppi\|=\left\{
\begin{array}
    {llll} \dst|\chi_{i-1}-\chi_i|\quad &i=1,\cdots, n-1\\
    \\
    \chi_{n-1}& i=n
\end{array}
 \right. $$ where $|w|:=\sqrt{w^2}$, for a given $w\in \complex$. \vskip.1in\noi{\bf g)} Reflections are not well described in the framework of the ${\cD ep}$-reduction: Compare, \eg, \cite[Section 4.4]{pinzari14}. Instead, in the framework of the $\cP$-reduction, the transformation $$(\bar\Theta, \bar\vartheta)\to (-\bar\Theta, 2k\p-\bar\vartheta)\qquad k\in \integer^{n-1}$$ corresponds to changing the sign of the second component of any $y^\ppi$ and any $x^\ppi$. Therefore, any of the points $$(\bar\Theta, \bar\vartheta)=(0,k\p)\qquad k\in \integer^{n-1}$$ is an equilibrium point for the Hamiltonian, corresponding to a co-planar configuration. Compare Proposition \ref{integrability}.

\backmatter
\bibliographystyle{amsalpha}

\begin{thebibliography}{10}

\bibitem{arnold63c}
V.~I. Arnold.
\newblock Proof of a theorem by {A}. {N}. {K}olmogorov on the invariance of
  quasi-periodic motions under small perturbations of the {H}amiltonian.
\newblock {\em Russian Math. Survey}, 18:13-40, 1963.

\bibitem{arnold63}
V.I. Arnold.
\newblock {S}mall denominators and problems of stability of motion in classical
  and celestial mechanics.
\newblock {\em Russian Math. Surveys}, 18(6):85-191, 1963.

\bibitem{boigey82}
F.~Boigey.
\newblock {\'E}limination des n\oe uds dans le probl{\`e}me newtonien des
  quatre corps.
\newblock {\em Celestial Mech.}, 27(4):399-414, 1982.

\bibitem{cellettiPi05}
A.~Celletti and G.~Pinzari.
\newblock Four classical methods for determining planetary elliptic elements: a
  comparison.
\newblock {\em Celestial Mech. Dynam. Astronom.}, 93(1-4):1-52, 2005.

\bibitem{chierchia13}
L.~Chierchia.
\newblock {T}he {P}lanetary {N}-{B}ody {P}roblem.
\newblock {\em UNESCO Encyclopedia of Life Support Systems}, 6.119.55, 2012.

\bibitem{chierchiaPi10}
L.~Chierchia and G.~Pinzari.
\newblock Properly-degenerate {KAM} theory (following {V}.{I}. {A}rnold).
\newblock {\em Discrete Contin. Dyn. Syst. Ser. S}, 3(4):545-578, 2010.

\bibitem{chierchiaPi11a}
L.~Chierchia and G.~Pinzari.
\newblock {D}eprit's reduction of the nodes revised.
\newblock {\em Celestial Mech.}, 109(3):285-301, 2011.

\bibitem{chierchiaPi11c}
L.~Chierchia and G.~Pinzari.
\newblock Planetary {B}irkhoff normal forms.
\newblock {\em J. Mod. Dyn.}, 5(4):623-664, 2011.

\bibitem{chierchiaPi11b}
L.~Chierchia and G.~Pinzari.
\newblock The planetary {$N$}-body problem: symplectic foliation, reductions
  and invariant tori.
\newblock {\em Invent. Math.}, 186(1):1-77, 2011.

\bibitem{chierchiaPi14}
L.~Chierchia and G.~Pinzari.
\newblock {M}etric stability of the planetary n-body problem.
\newblock {\em {P}roceedings of the International Congress of Mathematicians},
  2014.

\bibitem{DKRS15}
A.~Delshams, V.~Kaloshin, A.~de~la Rosa, and T.~M. Seara.
\newblock {G}lobal instability in the elliptic restricted three body problem.
\newblock {\em ar{X}iv: 1501.01214}, 2015.

\bibitem{deprit83}
A.~Deprit.
\newblock Elimination of the nodes in problems of {$n$} bodies.
\newblock {\em Celestial Mech.}, 30(2):181-195, 1983.

\bibitem{fejoz15}
J.~F{\'e}joz.
\newblock Work in progress.

\bibitem{fejoz04}
J.~F{\'e}joz.
\newblock D{\'e}monstration du `th{\'e}or{\`e}me d'{A}rnold' sur la
  stabilit{\'e} du syst{\`e}me plan{\'e}taire (d'apr{\`e}s {H}erman).
\newblock {\em Ergodic Theory Dynam. Systems}, 24(5):1521-1582, 2004.

\bibitem{fejoz13b}
J.~F{\'e}joz.
\newblock On action-angle coordinates and the {P}oincar\'e coordinates.
\newblock {\em Regul. Chaotic Dyn.}, 18(6):703-718, 2013.

\bibitem{fejoz13}
J.~F{\'e}joz.
\newblock On "{A}rnold's theorem" in celestial mechanics -a summary with an
  appendix on the poincar{\'e} coordinates.
\newblock {\em Discrete and Continuous Dynamical Systems}, 33:3555-3565, 2013.

\bibitem{FGKR14}
J.~Fejoz, M.~Guardia, V.~Kaloshin, and P.~Roldan.
\newblock {K}irkwood gaps and diffusion along mean motion resonances in the
  restricted planar three body problem.
\newblock {\em J. Eur. Math. Soc.}, 2014.

\bibitem{ferrerO94}
S.~Ferrer and C.~Os{\'a}car.
\newblock Harrington's {H}amiltonian in the stellar problem of three bodies:
  reductions, relative equilibria and bifurcations.
\newblock {\em Celestial Mech. Dynam. Astronom.}, 58(3):245-275, 1994.

\bibitem{harrington69}
R.~S. Harrington.
\newblock The stellar three-body problem.
\newblock {\em Celestial Mech. and Dyn. Astrronom}, 1(2):200-209, 1969.

\bibitem{herman09}
M.~R. Herman.
\newblock Torsion du probl{\`e}me plan{\'e}taire, edited by J. F{\'e}joz in
  2009.
\newblock Available in the electronic `Archives Michel Herman' at
  {\small\url{http://www.college-de-france.fr/default/EN/all/equ_dif/archives_michel_herman.htm}}.

\bibitem{jacobi1842}
C.~G.~J. Jacobi.
\newblock Sur l'{\'e}limination des noeuds dans le probl{\`e}me des trois
  corps.
\newblock {\em Astronomische Nachrichten}, Bd XX:81-102, 1842.

\bibitem{kolmogorov54}
A.N. Kolmogorov.
\newblock On the {C}onservation of {C}onditionally {P}eriodic {M}otions under
  {S}mall {P}erturbation of the {H}amiltonian.
\newblock {\em Dokl. Akad. Nauk SSR}, 98:527-530, 1954.

\bibitem{laskarR95}
J.~Laskar and P.~Robutel.
\newblock Stability of the planetary three-body problem. {I}. {E}xpansion of
  the planetary {H}amiltonian.
\newblock {\em Celestial Mech. Dynam. Astronom.}, 62(3):193-217, 1995.

\bibitem{leviCivita1904}
T.~Levi-Civita.
\newblock Sopra la equazione di {K}epler.
\newblock {\em {A}stronomische {N}achrichten}, 165(20):313-314, 1904.

\bibitem{maligeRL02}
F.~Malige, P.~Robutel, and J.~Laskar.
\newblock Partial reduction in the {$n$}-body planetary problem using the
  angular momentum integral.
\newblock {\em Celestial Mech. Dynam. Astronom.}, 84(3):283-316, 2002.

\bibitem{moser1962}
J.~Moser.
\newblock On invariant curves of area-preserving mappings of an annulus.
\newblock {\em Nachr. Akad. Wiss. G\"ottingen Math.-Phys. Kl. II}, 1962:1-20,
  1962.

\bibitem{pinzari-th09}
G.~Pinzari.
\newblock {\em {O}n the {K}olmogorov set for many-body problems}.
\newblock PhD thesis, {U}niversit{\`a} {R}oma {T}re, April 2009.

\bibitem{pinzari13}
G.~Pinzari.
\newblock Aspects of the planetary {B}irkhoff normal form.
\newblock {\em Regul. Chaotic Dyn.}, 18(6):860-906, 2013.

\bibitem{pinzari14}
G.~Pinzari.
\newblock Canonical coordinates for the planetary problem.
\newblock {\em Acta Applicandae Mathematicae}, pages 1-28, 2014.

\bibitem{Poincare:1892}
H.~Poincar{\'e}.
\newblock {\em Les m{\'e}thodes nouvelles de la m{\'e}canique c{\'e}leste}.
\newblock Gauthier-Villars, Paris, 1892.

\bibitem{poschel93}
J.~P{\"o}schel.
\newblock Nekhoroshev estimates for quasi-convex {H}amiltonian systems.
\newblock {\em Math. Z.}, 213(2):187-216, 1993.

\bibitem{radau1868}
R.~Radau.
\newblock Sur une transformation des {\'e}quations diff{\'e}rentielles de la
  dynamique.
\newblock {\em Ann. Sci. Ec. Norm. Sup.}, 5:311-375, 1868.

\bibitem{robutel95}
P.~Robutel.
\newblock Stability of the planetary three-body problem. {II}. {KAM} theory and
  existence of quasiperiodic motions.
\newblock {\em Celestial Mech. Dynam. Astronom.}, 62(3):219-261, 1995.

\bibitem{Russmann:2001}
H.~R{\"u}ssmann.
\newblock Invariant tori in non-degenerate nearly integrable {H}amiltonian
  systems.
\newblock {\em Regul. Chaotic Dyn.}, 6(2):119-204, 2001.

\bibitem{tisserand1889}
F.~Tisserand.
\newblock Trait\'e de m\'ecanique c\'eleste.
\newblock {\em Gauthier-Villars}, I, 1889-1896

\bibitem{zhao14}
L.~Zhao.
\newblock Partial reduction and {D}elaunay/{D}eprit variables.
\newblock {\em Celestial Mechanics and Dynamical Astronomy}, 120(4):423-432,
  2014.

\end{thebibliography}

\printindex
\end{document}